\documentclass[12pt]{article}
\usepackage{amsmath}
\usepackage{graphicx}
\usepackage{enumerate}
\usepackage{natbib}
\usepackage{url} 

\newcommand{\blind}{0}

\def\spacingset#1{\renewcommand{\baselinestretch}%
{#1}\small\normalsize} \spacingset{0}

\usepackage[english]{babel}
\usepackage{amsthm}
\newtheoremstyle{colon}%
{}
{}
{\itshape}
{}
{\bfseries}
{:}
{ }
{}

\theoremstyle{colon}
\newtheorem{ass}{Assumption}
\newtheorem{theorem}{Theorem}
\newtheorem{example}{Example}
\newtheorem{definition}{Definition}
\newtheorem{prop}{Proposition}
\newtheorem{rem}{Remark}

\usepackage{scalerel,stackengine}			
\usepackage{amssymb}
\usepackage{authblk}
\usepackage{amsmath, amsfonts, amssymb, bm, bbm}
\usepackage{graphicx}
\usepackage{listings}
\usepackage{float}
\usepackage[titletoc]{appendix}
\usepackage{hvfloat}
\usepackage{longtable}
\usepackage{mathtools}
\usepackage{booktabs,xcolor,colortbl}
\usepackage{xparse}
\usepackage[bottom]{footmisc}
\usepackage[font=small,skip=0pt]{caption}
\usepackage{mathrsfs}
\usepackage{enumitem}
\usepackage{rotating}
\usepackage{comment}
\usepackage{accents}
\usepackage{csquotes}
\usepackage{booktabs}
\usepackage{threeparttable}
\usepackage{setspace}
\excludecomment{confidential}
\newcommand\sbullet[1][.5]{\mathbin{\vcenter{\hbox{\scalebox{#1}{$\bullet$}}}}}

\newcommand{\aggregate}[2]{\underset{#2}{\operatornamewithlimits{#1\ }}}
\newcommand\blfootnote[1]{%
  \begingroup
  \renewcommand\thefootnote{}\footnote{#1}%
  \addtocounter{footnote}{-1}%
  \endgroup
}

\let\originalleft\left
\let\originalright\right
\renewcommand{\left}{\mathopen{}\mathclose\bgroup\originalleft}
\renewcommand{\right}{\aftergroup\egroup\originalright}

\title{On the Impact of Serial Dependence on Penalized Regression Methods}
\author{}

\date{}

\begin{document}

\def\spacingset#1{\renewcommand{\baselinestretch}%
{#1}\small\normalsize} \spacingset{1}

\date{}
\if0\blind
{\title{\bf On the Impact of Serial Dependence on Penalized Regression Methods}
  \author{\small Simone Tonini\thanks{Corresponding author: simone.tonini@santannapisa.it}\hspace{.2cm}\\
Institute of Economics \& EMbeDS, Sant'Anna School\\
 of Advanced Studies, Pisa, Italy\\
    and \\
    Francesca Chiaromonte\hspace{.2cm}\\
    Institute of Economics \& EMbeDS,\\
    Sant'Anna School of Advanced Studies, Pisa, Italy.\\
Department of Statistics, \\
Penn State University, USA.\\
  and  \\
    Alessandro Giovannelli\hspace{.2cm}\\
    University of L'Aquila Department of Information Engineering,\\
 Computer Science and Mathematics, L'Aquila, Italy}
  \maketitle
} \fi

\if1\blind
{ \bigskip
  \bigskip
  \bigskip
  \begin{center}
    {\LARGE\bf On the Impact of Serial Dependence on Penalized Regression Methods}
\end{center}
  \medskip
} \fi

\bigskip
\begin{abstract}
This paper characterizes the impact of covariate serial dependence on the non-asymptotic estimation error bound of penalized regressions (PRs). Focusing on the direct relationship between the degree of cross-correlation between covariates and the estimation error bound of PRs, we show that orthogonal or weakly cross-correlated stationary AR processes can exhibit high spurious correlations caused by serial dependence. We provide analytical results on the distribution of the sample cross-correlation in the case of two orthogonal Gaussian AR(1) processes, and extend and validate them through an extensive simulation study. Furthermore, we introduce a new procedure to mitigate spurious correlations in a time series setting, applying PRs to pre-whitened (ARMA filtered) time series. We show that under mild assumptions our procedure allows both to reduce the estimation error and to develop an effective forecasting strategy. The estimation accuracy of our proposal is validated through additional simulations, as well as an empirical application to a large set of monthly macroeconomic time series relative to the Euro Area.
\end{abstract}

\noindent%
\vfill


\blfootnote{{\it Acknowledgments:} The authors wish to thank Marco Lippi for suggesting to develop the theoretical part of this work and for the precious technical support. We are also grateful to Sebastiano Michele Zema and Luca Insolia for their helpful comments and for the stimulating dialogues.}

\newpage
\spacingset{1.9} 

\section{Introduction}\label{Intro}

Much contemporary statistical literature is devoted to the problem of extracting information from large datasets, which are ubiquitous in many fields of science \citep{Fan13}. In the context of high-dimensional regression problems, where the number of variables is comparable to or larger than the sample size, coefficient estimates produced by ordinary least squares (OLS)  
can be hindered by massive variance inflation. The use of 
penalized procedures that introduce a shrinkage in the OLS 
estimator is a widely accepted solution to this problem. In this paper we 
consider the most commonly used penalized regressions (PRs);
namely, those based on the $\ell_1$-penalty \citep{tibshirani96,ada}, the $\ell_2$-penalty \citep{Hoerl70} and their combinations \citep{Zou05}. 
Depending on the form of the penalty, PRs can produce {\it dense solutions}, where coefficients may have small yet non-zero estimates, or {\it sparse solutions}, where less relevant predictors have coefficient estimates equal to zero.


From a theoretical standpoint, we utilize the work of \cite{zhaoyu2006, Bickel2009, Lounici2009, Negahban2009,  Zou2009, Negahban2012, Hastie2015, Xin2017}, extending 
some of their results to the time series
setting. 
Specifically, studying the estimation properties of PRs, these authors 
showed that their non-asymptotic estimation error bound depends critically on the degree of cross-correlation between covariates. In summary, PRs are most effective when covariates are orthogonal or weakly cross-correlated, since the bound is inversely proportional to the minimum eigenvalue of the sample cross-correlation matrix of the covariates themselves. In this respect, two different
situations may occur. 
In the first, the covariates are truly multicollinear; that is, 
cross-correlations exist at the population level. In the second, 
correlations are spurious; the covariates can be orthogonal or nearly orthogonal at the population level, but 
other mechanisms 
generate cross-correlations in the sample. This 
matter becomes more prominent and consequential in 
high dimension, leading to false scientific discoveries and wrong statistical inferences \citep{Fan13, Fan16, Fan18}.


While issues related to
multicollinear time series 
have been extensively studied by \cite{Lippi2000, Lippi2001, SeW2002a, demol2008, Medeiros2012, Fan2020}, the 
effects of spurious correlations on PRs  
in the context of time series data has not been fully addressed to date. 


The main 
objective of this paper is to 
shed light on the role that 
covariates' serial dependence can play for
PRs through spurious cross-correlations. 
We show that, in addition to 
autocorrelation in the residuals
\citep{Bartlett35, Granger73, Granger2001}, also the serial dependence
in the covariates is critical for the estimation of regression coefficients. 
Our work introduces two elements of novelty. 
First, we formalize the impact of covariates' serial dependence on the non-asymptotic estimation error bound of PRs, showing that orthogonal or weakly cross-correlated stationary AR processes can exhibit high spurious correlations caused by serial dependence. Specifically, we demonstrate 
that the probability of spurious correlation between stationary processes depends, in addition to
the sample size, 
on 
their degree of serial dependence. To
prove this we derive the density of the sample cross-correlation between orthogonal stationary Gaussian AR(1) processes. 
Notably, the result can be generalized to high-dimensional correlation matrices
due to the fact that the minimum eigenvalue of a sample correlation matrix is bounded from above by the maximum absolute value of its off diagonal entries.
Second, to improve the estimation performance of PRs in a time series context, we introduce a new procedure based on applying PRs to the pre-whitened variables. We show that, under mild assumptions, our proposal produces more accurate estimates of regression coefficients, as well as 
improved selection of 
relevant 
covariates in sparse regimes. To validate our procedure we conduct an extensive simulation study.
Furthermore, to illustrate the validity of our results in a high-dimension context, we report an empirical exercise on forecasting the consumer price index 
by means of a 
set of 309 macroeconomic time series 
relative to the Euro Area.

Our main results can be summarized as follows. ($i$) Through our 
theoretical 
analysis we show that, whenever 
the autocorrelation coefficients of 
AR(1) processes 
have the same sign, an
increase in the degree of serial dependence 
induces an increase in the probability of large spurious cross-correlations.
($ii$) Through
simulations
we 
show that the association between serial dependence 
and the probability of large spurious
correlations holds 
much more generally, e.g.,~in 
cases where the processes are non-Gaussian, weakly cross-correlated or generated by models 
other than than AR(1). 
This highlights that a small minimum eigenvalue is more likely in finite realizations of serially dependent weakly cross-correlated (or orthogonal) processes, compared to the case of independent samples -- 
and thus that 
serially dependent covariates 
can cause major problems for the estimation accuracy of PRs
-- something that is numerically corroborated when we apply our 
proposed procedure to simulated data. ($iii$) 
Our empirical exercise shows that LASSO \citep{tibshirani96} applied to ARMA residuals generates more parsimonious models 
more accurate forecasts compared to
LASSO applied directly to the time series under consideration.


%

The 
reminder of the paper is organized as follows. In Section~2 we 
describe the problem setup and our contribution. In Section~3 we 
present our theoretical result on the impact of covariates serial dependence on sample cross-correlation. In Section~4 we provide simulation studies to corroborate and extend the theoretical results of Section~3. In Section~5 we introduce and evaluate our proposal for mitigating the adverse effects of serial dependence using simulated and real data. In
Section~6 
we provide some final remarks. 


The following notations 
will be used throughout 
the paper. For any dimension $p$, bold
letters denote vectors and the corresponding regular letters 
their elements, for example $\pmb{a}=(a_1,a_2,\dots,a_p)'$.
Supp$(\pmb{a})$ denotes the support of a vector, that is, 
$\{j\in\{1,2,\dots,p\}:a_j\neq0\}$, and $|\text{Supp}(\pmb{a})|$ the support cardinality. 
The $\ell_q$ norm of 
a vector is $||\pmb{a}||_q\coloneqq\left(\sum_{j=1}^p|a_j|^q\right)^{1/q}$ for $0<q<\infty$, with $||\pmb{a}||_q^k\coloneqq\left(\sum_{j=1}^p|a_j|^q\right)^{k/q}$, and with the usual extension $||\pmb{a}||_0\coloneqq|\text{Supp}(\pmb{a})|$.
Bold capital 
letters denote matrices, for example $\pmb{A}$, where $\left(\pmb{A}\right)_{ij}=a_{ij}$ denotes its $i$-row $j$-column element. Moreover, $\pmb{0}_{p}$ denotes a $p$-length vector of zeros, while $\pmb{I}_{p}$ denotes a $p\times p$ identity matrix. Finally, $Sign(r)$ indicates the sign of a real number $r$.

\section{Problem Setup and Our Contribution}

\subsection{State-of-the-Art}\label{StateArt}

Let $\mathbf{X} = \left\{\mathbf{x}_{t}\right\}_{t=1}^{T}$ denote an $n \times T$ rectangular array of observations 
on $n$ covariates, and 
$\mathbf{y}=\left\{y_{t}\right\}_{t=1}^{T}$ a $1\times T$ response vector. Assume that $\mathbf{y}$  and $\mathbf{X}$ are realizations of strictly Gaussian stationary and absolutely regular processes $\left\{\left( y_{t}, \mathbf{x}_t' \right) \in \mathbb{R}^{1+n}, t \in\mathbb{Z}, n \in\mathbb{N}\right\}$ defined on the probability space $\left( \Omega, \mathcal{F},P\right)$.
Let $\mathbf{C}_x=E[\mathbf{x}_t\mathbf{x}_{t}']$ and
$\widehat{\mathbf{C}}_x=\frac{1}{T-1}\mathbf{X}\mathbf{X}'$ 
denote the cross-covariance matrix and its estimate, with generic element $\widehat{c}_{ij}^x$ and eigenvalues $\widehat{\psi}_{max}^x\geq\ldots \geq\widehat{\psi}_{min}^x$. Finally, assume that each $\mathbf{x}_i$ 
is standardized so that $\frac{1}{T-1}\sum_{t}x_{it}=0$ and $\frac{1}{T-1}\sum_{t}x_{it}^2=1$. We consider the following data generating process (DGP) for the response 
$$
\mathbf{y}=\mathbf{X}'\pmb{\alpha}^*+\pmb{\varepsilon},
$$
where $\pmb{\alpha}^*$ is the $n\times1$ unknown $s$-sparse vector of regression coefficients, i.e.~$||\pmb{\alpha}^*||_0=s<n$, and $\pmb{\varepsilon}\in\mathbb{R}^T$ is a random noise vector. 
If $n \geq T$,
$\pmb{\alpha}^*$ is 
estimated solving 
a convex optimization 
that combines a quadratic loss
and a regularization penalty:
\begin{equation}\label{PRs}
\widehat{\pmb{\alpha}}= \aggregate{argmin}{\pmb{\alpha}\in \mathbb{R}^n}\:\left \{\:\frac{1}{2T}\ \ ||\mathbf{y}-\mathbf{X}'\pmb{\alpha}||_2^2 + \lambda\ell(\pmb{\alpha}) \right \}.
\end{equation}
Here $\lambda>0$ 
represents the weight of the penalty, and $\ell:\mathbb{R}^n\rightarrow\mathbb{R}^+$ is a norm. 
The following definitions will be used.
%
%
\begin{definition}\label{def:uno}(Strong Convexity):
Given a differentiable function $\mathcal{L}:\mathbb{R}^n\rightarrow\mathbb{R}$ and the vector differential operator $\bigtriangledown$, we say that $\mathcal{L}$ is strongly convex with parameter $\gamma>0$ at $\pmb{a}\in\mathbb{R}^n$ if, for all $\pmb{b}\in\mathbb{R}^n$, 
$\mathcal{L}(\pmb{b})-\mathcal{L}(\pmb{a})\geq\bigtriangledown \mathcal{L}(\pmb{a})' (\pmb{b}-\pmb{a})+\frac{\gamma}{2}\left |\left |\pmb{b}-\pmb{a}\right |\right |_2^2$. 
\end{definition}
%
%
\noindent Strong Convexity guarantees a small coefficient estimation error (see \citealt{Negahban2009, Negahban2012}). In particular, when $\mathcal{L}$ (the loss function) is \enquote{sharply curved} around its optimum $\widehat{\pmb{\alpha}}$, a small
$|\mathcal{L}(\widehat{\pmb{\alpha}})-\mathcal{L}(\pmb{\alpha}^*)|$
guarantees that 
$||\widehat{\pmb{\alpha}}-\pmb{\alpha}^*||_2$ is also small.
The parameter $\gamma$ 
governs the strength of convexity;
when 
$\mathcal{L}$ is twice differentiable,
strong convexity requires the minimum eigenvalue of the Hessian $\bigtriangledown^2\mathcal{L}(\pmb{\alpha})$
to be at least $\gamma$ for all $\pmb{\alpha}$
in a neighborhood of 
$\pmb{\alpha}^*$. Thus, since its Hessian is $\bigtriangledown^2\mathcal{L}(\pmb{\alpha})=\widehat{\mathbf{C}}_x$ for all $\pmb{\alpha}\in\mathbb{R}^n$, 
the quadratic loss $\mathcal{L}(\pmb{\alpha})=\frac{1}{2T}||\mathbf{y}-\mathbf{X}'\pmb{\alpha}||_2^2$
is strongly convex with parameter $\gamma$ if and only if $\widehat{\psi}_{min}\geq\gamma$ (see \citealt[p.~293]{Hastie2015}). 
Consequently, in this case $||\widehat{\pmb{\alpha}}-\pmb{\alpha}^*||_2$ depends on $\widehat{\psi}_{min}$. It
is 
important to note that when $n>T$ the
quadratic loss cannot be strongly convex since $\widehat{\mathbf{C}}_x$ is 
singular and thus $\widehat{\psi}_{min}=0$. In this
case \cite{Bickel2009} proposed a
{\it Restricted Eigenvalue Condition}, which is essentially a restriction on the eigenvalues of $\widehat{\mathbf{C}}_x$ as a function of the degree of sparsity, $s$. The Restricted Eigenvalue Condition allows for 
strong convexity (Definition \ref{def:uno}) to hold in the singular case, and we refer to this 
as {\it Restricted Strong Convexity} (see \citealt{Negahban2012}; we provide more details in
Supplement \ref{ResEin}). 
%
%
%
\begin{definition}\label{def:due}(Dual Norm and Subspace Compatibility Constant): Given a norm $\ell$ and the inner product $\langle\cdot,\cdot\rangle$, we define the dual norm of $\ell$ 
as
$\ell^*(\pmb{v})\coloneqq\sup_{\mathbf{u}\in\mathbb{R}^n\backslash\{0\}}\frac{\langle\mathbf{u},\mathbf{v}\rangle}{\ell(\mathbf{u})}$.
%
%
For any subspace $\mathcal{A}$ of $\mathbb{R}^n$ that captures the constraints 
underlying \eqref{PRs}, we define the subspace compatibility constant with respect to the pair ($\ell,||\cdot||_2$)
as
$\Psi(\mathcal{A})\coloneqq \sup_{\mathbf{u}\in\mathcal{A}:\mathbf{u}\neq0}\frac{\ell(\mathbf{u})}{||\mathbf{u}||_2}$.
\end{definition}
%
%
%
\noindent The following Proposition is derived 
from Corollary~1 
of \cite{Negahban2012} and provides the non-asymptotic coefficient estimation error bound for PRs.

\begin{prop}\label{prop:uno} 
Consider the convex optimization problem in \eqref{PRs}. 
Suppose that the penalty parameter $\lambda$ is strictly positive and $\geq2\ell^*\left(\frac{1}{T}\mathbf{X}\pmb{\varepsilon}\right)$, 
and that strong convexity 
holds with parameter $\gamma>0$. Then, any optimal solution $\pmb{\widehat{\alpha}}$
satisfies the bound
$||\widehat{\pmb{\alpha}}-\pmb{\alpha}^*||_2\leq3\frac{\lambda}{\gamma}\Psi(\mathcal{A})$.
\end{prop}

\noindent\textbf{Proof:} See  Corollary 1 in \cite{Negahban2012}. \hfill $\blacksquare$

\vspace{0.3cm}

\noindent 
The coefficient estimation error bound in Proposition~\ref{prop:uno}
increases with the 
penalty parameter 
$\lambda$, which must be strictly positive and 
$\geq 2\ell^*\left(\frac{1}{T}\mathbf{X}\pmb{\varepsilon}\right)$; 
increases with the subspace compatibility constant $\Psi(\mathcal{A})$, which in turn increases with the size of the model subspace $\mathcal{A}$; 
and decreases with the convexity parameter $\gamma$. \cite{Negahban2009,Negahban2012} derive the 
bound 
for PRs in the case of 
independent observations 
(no serial dependence). To this end, the authors compute the probability that $\lambda\geq2\ell^*\left(\frac{1}{T}\mathbf{X}\pmb{\varepsilon}\right)$ 
when the entries of $\mathbf{X}$ and $\pmb{\varepsilon}$ are sub-Gaussian, and assume that 
strong convexity (or 
restricted strong convexity) holds with parameter $\gamma$, i.e.~that $\widehat{\psi}_{min}^x\geq\gamma$ (see Corollary~2 in \citealt{Negahban2012} 
and Corollary~6 in \citealt{Negahban2009} for examples of sparse and
dense PRs, respectively).


This analysis shows the role of covariates cross-correlation 
in determining the estimation accuracy of PRs. In particular, Proposition~\ref{prop:uno} shows that PRs 
perform better if covariates are orthogonal or weakly correlated in the sample, since 
strong sample cross-correlations 
correspond to small $\widehat{\psi}_{min}^x$. 
As mentioned in the Introduction, 
strong sample cross-correlations may 
reflect true multicollinearities at the population level. 
In this case \cite{Fan2020} focus on the fact that time series 
multicollinearities can be 
captured 
with Factor Models, and
propose to apply PRs on the estimated idiosyncratic components obtained by filtering the observed time series through estimated factors.
However, 
strong sample cross-correlation may also be spurious; 
this is the case we wish to tackle 
in the particular context of time series.

\subsection{Our Contribution}\label{MainCont}

We argue that spurious correlations are one of the causes that potentially limits the use of PRs in time series. In particular, we focus on the implications of covariates serial dependence on $\widehat{\psi}_{min}^x$, which determines the \enquote{strength} of strong convexity (see Definition \ref{def:uno})
and is one of the main components of the PRs error bound 
in Proposition \ref{prop:uno}.
In this respect, we relax the 
assumption that strong convexity (or restricted strong convexity) holds with parameter $\gamma>0$, showing that the probability of getting $\widehat{\psi}_{min}^x\leq\gamma$ increases with the covariates' serial dependence. Note that 
in order to focus on $\widehat{\psi}_{min}^x$ 
we are assuming that $\widehat{\mathbf{C}}_x$ is positive definite. 
If $n>T$ the matrix $\widehat{\mathbf{C}}_x$ is singular,
but our results are still valid considering the probability of a restricted eigenvalue 
$\leq \gamma$ (see \citealt{Bickel2009}).


The 
structure 
of our theoretical contribution is
as follows. Given $\gamma=1-\tau$, $\tau\in[0,1)$, and the upper bound $\widehat{\psi}_{min}^x\leq1-\max_{i\neq j}|\widehat{c}_{ij}^x|$, we emphasize the role of a generic off-diagonal element of 
$\widehat{\mathbf{C}}_x$ in determining the probability that $\widehat{\psi}_{min}^x\leq\gamma$
through the 
inequalities
\begin{equation}\small\label{Eigen_Bound}\Pr\left \{\widehat{\psi}_{min}^x\leq1-\tau\right \}\geq \Pr\left \{1-\max_{i\neq j}|\widehat{c}_{ij}^x|\leq1-\tau\right \} \geq \Pr\left \{1-|\widehat{c}_{i\neq j}^x|\leq1-\tau\right \}=\Pr\left \{|\widehat{c}_{i\neq j}^x|\geq\tau\right \}\footnotesize.\end{equation} 
%
%
%
\noindent Thus, 
$\Pr\{|\widehat{c}_{i\neq j}^x|\geq\tau\}$ plays a role in determining the probability of dealing with a small $\widehat{\psi}_{min}^x$
and consequently, through strong convexity, on
the PRs
estimation error bound presented in Proposition~\ref{prop:uno}.
It follows that any impact of the degree of serial dependence on $\Pr\{|\widehat{c}_{i\neq j}^x|\geq\tau \}$ results in an impact 
on such bound.

\begin{figure}[H]
\begin{center}
\graphicspath{{images/}}
\includegraphics[scale=0.4]{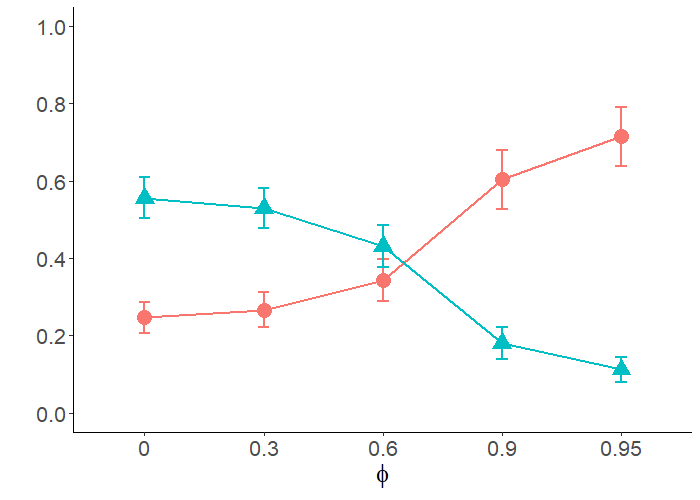}
\caption{\footnotesize Results from a numerical toy example. The orange circles and bars
represent 
means and standard deviations of $\max_{i\neq j}|\widehat{c}_{ij}^x|$ for various values of the autocorrelation $\phi$, as obtained from 5000 Monte Carlo replications. Similarly, 
blue triangles and bars 
represent means and standard deviations of $\widehat{\psi}_{min}^x$. 
}\label{Fig:toy}
\end{center}
\end{figure}

To better 
illustrate our reasoning, we introduce a toy example where we show numerically the impact of serial dependence on $\max_{i\neq j}|\widehat{c}_{ij}^x|$ and $\widehat{\psi}_{min}^x$. 
We generate 10 processes from the model 
$\mathbf{x}_t= \mathbf{D}_{\phi}\mathbf{x}_{t-1}+\mathbf{u}_t$, $t=1,\dots,100$, where $\mathbf{D}_{\phi}$ is a $10\times10$ diagonal matrix with the same autocorrelation coefficient $\phi$ in all positions along the main diagonal, and $\mathbf{u}_t\sim N\left(\pmb{0}_{10},\pmb{I}_{10}\right)$. Note that for these AR(1) processes the degree of serial dependence is determined by $|\phi|$ and, 
since the processes are orthogonal, 
the minimum eigenvalue 
of the population cross-correlation matrix $\mathbf{C}_x$ is $\psi_{min}^x=1$.
%
We consider five values for 
$\phi$, namely $0.0,0.3,0.6,0.9,0.95$, and 
for each 
we calculate the average and the standard deviation of both $\max_{i\neq j}|\widehat{c}_{ij}^x|$ and $\widehat{\psi}_{min}^x$ on 5000 Monte Carlo replications. Results are reported in Figure \ref{Fig:toy}. 
We observe that the stronger 
the persistence of the process ($\phi$ closer to 1) the higher is the probability of a large spurious sample correlation (orange circle), which in turn leads to a small minimum eigenvalue of the sample cross-correlation matrix (blue triangle). 

%
In 
light of these results, 
our next task is to derive the finite sample density of $\widehat{c}_{ij}^x$ for the purpose of formalizing the impact of serial 
dependence on $\Pr\{|\widehat{c}_{i\neq j}^x|\geq\tau\}$. 
%
%
It is noteworthy that, when the covariates have a factor-based structure, 
strong spurious correlations in the sample may affect the idiosyncratic components if they are serially dependent, 
reducing the accuracy of the procedure proposed by \cite{Fan2020} (see Supplement \ref{PopCrosCorr}) .

\section{\setlength{\baselineskip}{0.5\baselineskip}Density of the Sample Correlation between two Orthogonal AR(1) Gaussian Processes}\label{Density_c}

In this section we present 
our main theoretical contribution 
concerning the impact of covariates' serial dependence on the non-asymptotic estimation error bound of PRs,
showing that the probability of incurring in strong spurious correlation increases with serial dependence. 

Consider a first order bivariate autoregressive process $\mathbf{x}_t=\mathbf{D}_{\phi}\mathbf{x}_{t-1}+\mathbf{u}_t$, $t=1,\dots,T$, where $\mathbf{D}_{\phi}$ is a $2\times2$ diagonal matrix with main diagonal elements
$\phi_1, \phi_2 < 1$. We make the following assumption about the bivariate vector of autoregressive residuals:
\begin{ass}\label{ass:uno}
$\mathbf{u}_t\sim N\left(\pmb{0}_2, \pmb{I}_2\right)$.
\end{ass}
\noindent Therefore $\pmb{x}_t \sim$ 
$N\left(\pmb{0}_2,\mathbf{C}_x\right)$
with $\left(\mathbf{C}_x\right)_{ii}=\frac{1}{1-\phi_i^2}$, $i=1,2$, and $\left(\mathbf{C}_x\right)_{12}={c}_{12}^x=0$.
%
%
In this 
setting we focus on the density of the sample correlation coefficient defined as
\begin{equation}\label{SampleCorr}
\widehat{c}_{12}^x=\frac{a_{12}}{\sqrt{a_{11}}\sqrt{a_{22}}},
\end{equation}
where $a_{i,j}=\sum_{t=1}^T(x_{it}-\overline{x}_i)(x_{jt}-\overline{x}_j)=\sum_{t=1}^Tx_{it}x_{jt}$, since $\overline{x}_i=0$,  $i,j=1,2$. 
In particular, when $c_{12}^u=0$, $b=a_{21}/a_{11}$ and $v=a_{22}-a_{21}^2/a_{11}$, then
\begin{equation}\label{AndersonEq}\cfrac{\sqrt{a_{11}}\ b}{\sqrt{v/(T-2)}}=\sqrt{T-2}\cfrac{a_{12}/\sqrt{a_{11}a_{22}}}{\sqrt{1-a_{12}^2/(a_{11}a_{22})}}=\sqrt{T-2}\cfrac{\widehat{c}_{12}^x}{\sqrt{1-(\widehat{c}_{12}^x)^2}}.\end{equation}
%
%
\noindent Note that $b$ is the least squares regression coefficient of $x_{2t}$ on $x_{1t}$, and $v$ is the sum of the square of residuals of such regression. Thus, to derive the finite sample density of $\widehat{c}_{12}^x$ we need the sample densities of $b$ and $v$.
%
%
\begin{rem}\label{rem:uno}
In contrast to asymptotic statements, 
our theoretical analysis 
is intended to derive distributions and densities of estimators that hold for $T<\infty$. 
Hence we will not employ the usual concepts of convergence in probability
and in distribution; 
rather 
we will use a notion of approximation, whose \enquote{precision} needs to be evaluated. The precision of 
our approximations will be extensively tested under several finite $T$ scenarios in both the simulation study provided in Section \ref{MonteCarlo} and in the Supplement.
\end{rem}
%
%
\paragraph{\large Sample Distribution of $b$.}
We start by deriving the sample distribution of $b$, 
the OLS regression coefficient
for $x_2$ on $x_1$. The same holds if we regress $x_1$ on $x_2$.
%
%
%
\begin{prop}\label{prop:due} Under Assumption \ref{ass:uno} the sample distribution of $b$ is approximately \\
$N\left(0, \frac{(1-\phi_1^2\phi_2^2)(1-\phi_1^2)}{(T-1)(1-\phi_2^2)(1-\phi_1\phi_2)^2}\right)$.
\end{prop}

\noindent\textbf{Proof:} See Supplement \ref{Proof_Prop2}

\smallskip

%
\noindent Proposition~\ref{prop:due} shows that the OLS estimate $b$ 
is normally distributed with a variance that strongly depends on the degree of covariates serial dependence. In this
context, it is common
to adjust the standard error of the OLS 
to achieve consistency 
in the presence of heteroskedasticity and/or serial dependence;
this leads, for instance, to the Heteroskedasticity and Autocorrelation Consistent (HAC) estimator of
\cite{NeweyWest87} (NW). However, NW estimates can be highly sub-optimal (or inefficient) in the presence of strong serial dependence \citep{Kapetanios2022}. In Supplement \ref{Dist_b} we provide a simulation study to corroborate the result in Proposition \ref{prop:due}.

\paragraph{\large Sample Distribution of $v$.}
Here we derive the sample distribution of the sum of the square of residuals obtained by regressing $x_2$ on $x_1$. 

%
%
%
%
\begin{prop}\label{prop:tre} 
Under Assumption \ref{ass:uno} the sample distribution of $v$ is approximately \\
$\Gamma\left(\frac{T-2}{2}, \frac{2}{1-\phi_2^2}\right)$.
\end{prop}

\noindent\textbf{Proof:} See Supplement \ref{Proof_Prop3}

\paragraph{\large Sample Density of $\widehat{c}_{12}^x$.}
Note that $b$ and $v$ are independent. Using 
Propositions~\ref{prop:due} and \ref{prop:tre} and Equation~\eqref{AndersonEq} we can now derive the
density 
of the sample distribution of $\widehat{c}_{12}^x$.
%
%
\begin{theorem}\label{theorem:uno} 
Let $\mathbf{x}_t$ be a stationary bivariate Gaussian AR(1) process with autoregressive residuals distributed according to $N\left(\pmb{0}_2,\pmb{I}_2\right)$. Further, let $\phi_{12}=\phi_1\phi_2$ where $\phi_i$, $i=1,2$, are the autoregressive coefficients. Then, the sample density of $\widehat{c}_{12}^x$ is approximated by
\begin{equation}\label{r_dist}
\small\mathcal{D}_{\widehat{c}_{12}^x} = \frac{\Gamma\left(\frac{T-1}{2}\right)(1-\phi_{12})}{\Gamma\left(\frac{T-2}{2}\right)\sqrt{\pi}}(1-(\widehat{c}_{12}^x)^2)^{\frac{T-4}{2}}\left(1-\phi_{12}^2\right)^{\frac{T-2}{2}}\left(\frac{1}{1-\phi_{12}^2+2(\widehat{c}_{12}^x)^2\phi_{12}(\phi_{12}-1)}\right)^{\frac{T-1}{2}}.
\end{equation}
\end{theorem}

\noindent
\textbf{Proof:} See Supplement \ref{Proof_Theorem}

\vspace{0.3in}

\begin{rem}\label{rem:tre} 
$\mathcal{D}_{\widehat{c}_{12}^x}$ is the density of the sample correlation coefficient \eqref{SampleCorr} based on a finite $T$,
with serial dependence expressed by $\phi_1$ and $\phi_2$, and under the assumption of orthogonal Gaussian AR(1) processes.
\end{rem}

\begin{rem}\label{rem:quattro} 
From \eqref{r_dist} we see that $\phi_{12}$ determines the density of $\widehat{c}_{12}^x$ trough both its magnitude and sign. More precisely, when $Sign(\phi_1)=Sign(\phi_2)$, the probability 
in the tails increases as $|\phi_{12}|$ grows. On the other hand, when $Sign(\phi_1)\neq Sign(\phi_2)$, an
increase in $|\phi_{12}|$ leads to a density more concentrated 
around the origin. This 
peculiarity on the effect of $Sign(\phi_{12})$ will be numerically
explored and validated in Section~\ref{MonteCarlo}. 
\end{rem}
%
\noindent Theorem \ref{theorem:uno} shows that,
in a finite $T$ context, the probability of observing
sizeable spurious correlation between orthogonal Gaussian Autoregressive processes crucially depends on the degree of serial dependence. This 
has important consequences on the non-asymptotic performance of PRs for the reasons that we pointed out at the beginning of Section \ref{MainCont}, related to the 
role of $\Pr\left\{|\widehat{c}_{12}^x|\geq\tau\right\}$ (see inequality \eqref{Eigen_Bound}, Definition \ref{def:uno} and Proposition \ref{prop:uno}). The implication of Theorem \ref{theorem:uno} 
for such probability can be summarized in the following remark.

\begin{rem}\label{rem:cinque}
Because of Theorem \ref{theorem:uno} 
$
\ \Pr\left\{|\widehat{c}_{12}^x|\geq\tau\right\}\approx\int_{-1}^{-\tau}\mathcal{D}_{\widehat{c}_{12}^x}d\widehat{c}_{12}^x+\int_{\tau}^1\mathcal{D}_{\widehat{c}_{12}^x}d\widehat{c}_{12}^x,\ 
$
which depends on the degrees of serial dependence of the processes.
\end{rem}

\section{Monte Carlo Experiments}\label{MonteCarlo}
In this Section we conduct Monte Carlo experiments to assess numerically the approximation of the density of $\widehat{c}_{12}^x$ 
described in Section \ref{Density_c}. Then, we expand the theoretical results in more generic contexts, relaxing the assumption
that 
the covariates are orthogonal Gaussian AR(1) processes. 
We 
indicate the density of $\widehat{c}_{12}^x$ obtained by simulations as $d_s(\widehat{c}_{12}^x)$.

\subsection{Numerical Approximation of $d_s(\widehat{c}_{12}^x)$ to $\mathcal{D}_{\widehat{c}_{12}^x}$}

\begin{figure}[H]
\graphicspath{{images/}}
\centering
\subfloat[$d_s(\widehat{c}_{12}^x)$, $T=50$]{\includegraphics[width=4.5cm]{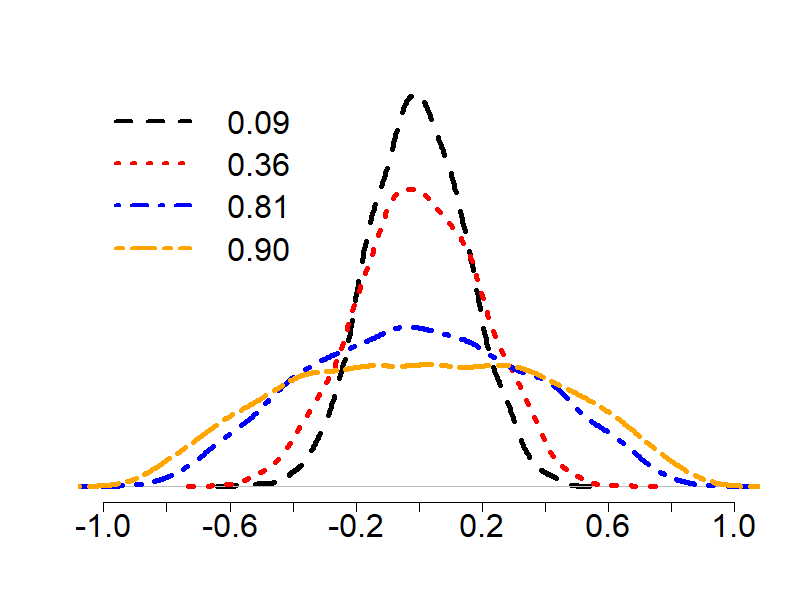}}\hfil
\subfloat[$d_s(\widehat{c}_{12}^x)$, $T=100$]{\includegraphics[width=4.5cm]{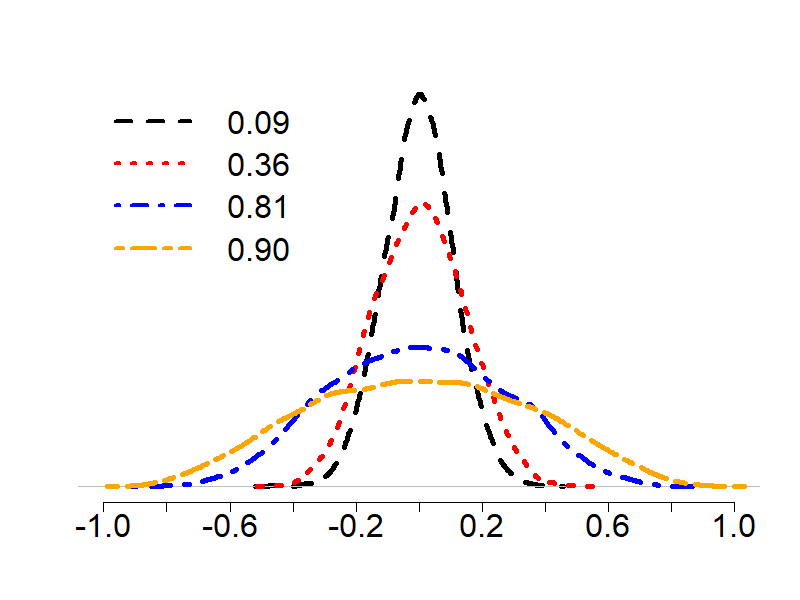}}\hfil 
\subfloat[$d_s(\widehat{c}_{12}^x)$, $T=250$]{\includegraphics[width=4.5cm]{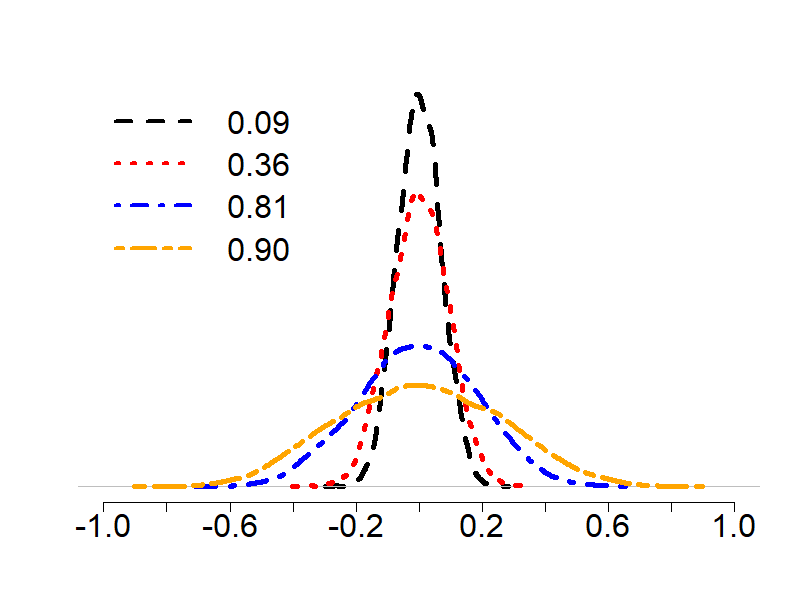}} 

\subfloat[$\mathcal{D}_{\widehat{c}_{12}^x}$, $T=50$]{\includegraphics[width=4.5cm]{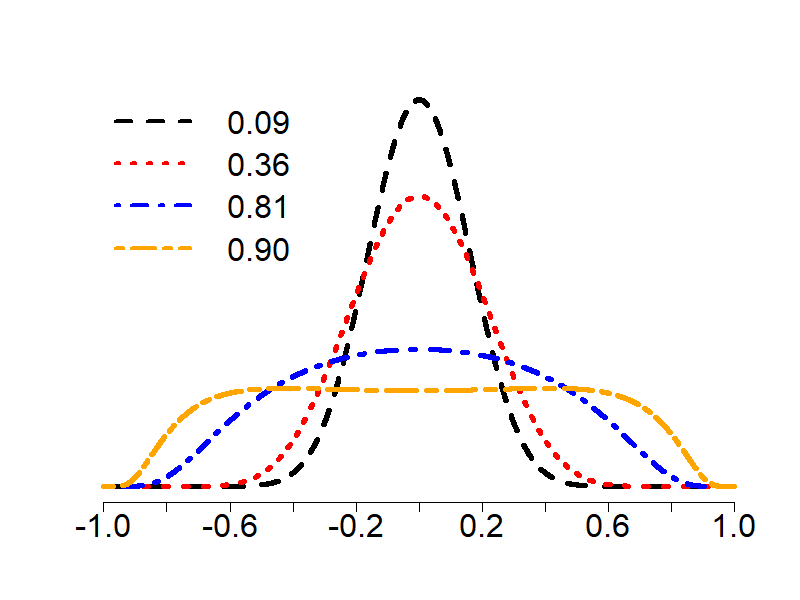}}\hfil   
\subfloat[$\mathcal{D}_{\widehat{c}_{12}^x}$, $T=100$]{\includegraphics[width=4.5cm]{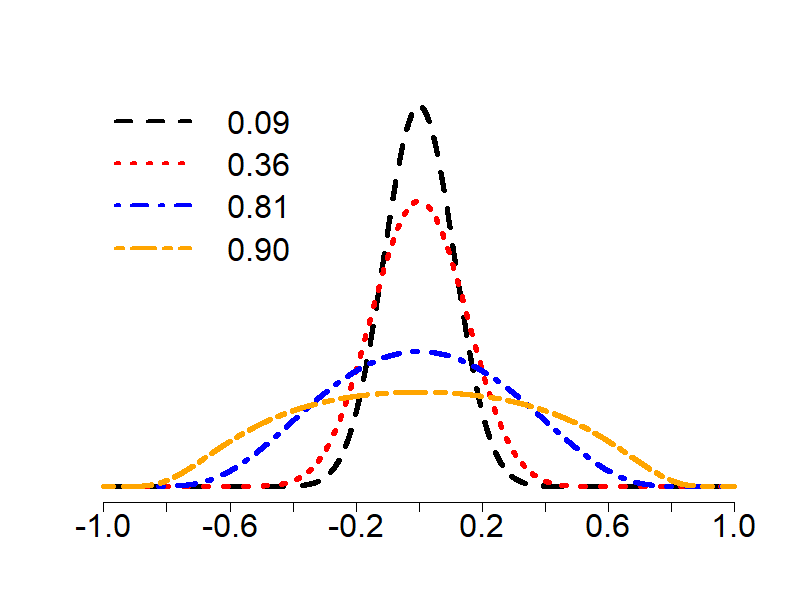}}\hfil
\subfloat[$\mathcal{D}_{\widehat{c}_{12}^x}$, $T=250$]{\includegraphics[width=4.5cm]{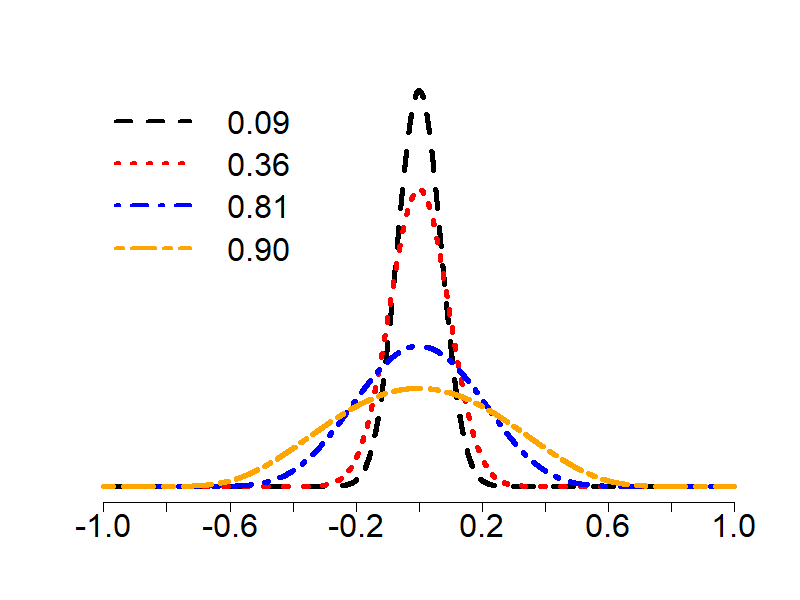}}
\caption{\footnotesize Monte Carlo 
densities for $\widehat{c}_{12}^x$ 
(top) and corresponding $\mathcal{D}_{\widehat{c}_{12}^x}$ 
(bottom) for various 
$T$ and $\phi_{12}$.}\label{DistributionsMC}
\end{figure}

We generate data from the bivariate process $\mathbf{x}_{t}=\mathbf{D}_{\phi}\mathbf{x}_{t-1}+\mathbf{u}_{t}$ for $t=1,\dots,T$, where $\mathbf{D}_{\phi}$ is a $2\times2$ diagonal matrix with the same autocorrelation coefficient $\phi$ in both 
positions along the diagonal, and $u_{t}\sim N(\pmb{0}_2, \pmb{I}_2)$. We consider $T=50, 100, 250$ and $\phi=0.3, 0.6, 0.9, 0.95$ -- thus, the parameter $\phi_{12}$ in $\mathcal{D}_{\widehat{c}_{12}^x}$, here equal to $\phi^2$, takes values 0.09, 0.36, 0.81, 0.90. The first row of Figure \ref{DistributionsMC} (Plots (a), (b), (c)) shows, for various values of $T$ and $\phi_{12}$, the density $d_s(\widehat{c}_{12}^x)$ generated through 5000 Monte Carlo replications. The second row of Figure \ref{DistributionsMC} (Plots (d), (e), (f)) shows the corresponding $\mathcal{D}_{\widehat{c}_{12}^x}$. These were plotted using 5000 values of the argument starting at -1 and increasing by steps of size 0.0004 until 1. As expected, we observe that the degree of approximation of $d_s(\widehat{c}_{12}^x)$ to $\mathcal{D}_{\widehat{c}_{12}^x}$ improves as $T$ increases and/or $\phi_{12}$ decreases. In particular, Plots (a) and (d) in Figure \ref{DistributionsMC},
where $T=50$, show that $\mathcal{D}_{\widehat{c}_{12}^x}$ approximates 
$d_s(\widehat{c}_{12}^x)$ well for a low-to-intermediate degree of serial dependence ($\phi_{12}\leq0.36$, i.e.~$\phi\leq0.6$). In contrast, in cases with high degree of serial dependence ($\phi_{12}\geq0.81$, i.e.~$\phi\geq0.9$), $\mathcal{D}_{\widehat{c}_{12}^x}$ has larger tails compared to $d_s(\widehat{c}_{12}^x)$; that is, 
the latter over-estimates the probability of large spurious correlations. However, it is noteworthy that the difference between the two densities is negligible for $T\geq100$ 
(Figure~\ref{DistributionsMC}, Plots (b) and (e) for $T=100$, 
and Plots (c) and (f) for $T=250$), also with high degree of serial dependence  ($\phi_{12}\approx0.90$, i.e.~$\phi=0.95$). These numerical experiments corroborate the fact that the sample cross-correlation between orthogonal Gaussian AR(1) processes is affected by the degree of serial dependence in a way that is well approximated by $\mathcal{D}_{\widehat{c}_{12}^x}$. In fact, for a sufficiently large finite $T$, we observe that $\Pr\left\{|\widehat{c}_{12}^x|\geq\tau\right\}$, $\tau>0$, increases with $\phi_{12}$ in a similar way for $d_s(\widehat{c}_{12}^x)$ and $\mathcal{D}_{\widehat{c}_{12}^x}$.

\noindent
{\bf The Impact of $Sign(\phi_{12})$}
\normalsize

\noindent 
In Remark \ref{rem:quattro} we pointed out that the impact of $\phi_{12}$ on 
$\mathcal{D}_{\widehat{c}_{12}^x}$ 
depends on $Sign(\phi_{12})$. In particular, when $-1<\phi_{12}<0$, an increment on $|\phi_{12}|$ makes the density of $\widehat{c}_{12}^x$ more concentrated around $0$. 
In order to 
validate this result, 
we 
run simulations with
$T=100$ and 
different values for 
the second element of the diagonal of $\mathbf{D}_{\phi}$; namely,
$-0.3,-0.6,-0.9,-0.95$. 
Results are shown in Plots (a) and (b) of Figure \ref{negSignPhi}. 
In this case, we see that when $Sign(\phi_1)\neq Sign(\phi_2)$ and $|\phi_{12}|$ increases,
$d_s(\widehat{c}_{12}^x)$ 
increases its concentration around $0$ in a way that is, again, well approximated by $\mathcal{D}_{\widehat{c}_{12}^x}$.

\subsection{General Case}
\begin{figure}[H]
\graphicspath{{images/}}
\centering
  \captionsetup[subfigure]{oneside,margin={0.5cm,0cm}}
\subfloat[\scriptsize $d_s(\widehat{c}_{12}^x)$, $Sign(\phi_1)\neq Sign(\phi_2)$]{\includegraphics[width=5cm]{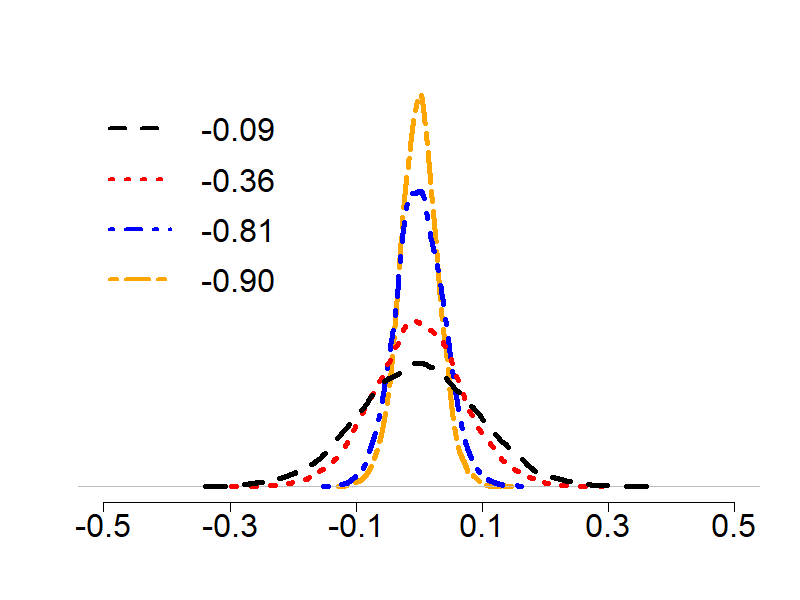}}\hfil
\subfloat[\scriptsize $\mathcal{D}_{\widehat{c}_{12}^x}$, $Sign(\phi_1)\neq Sign(\phi_2)$]{\includegraphics[width=5cm]{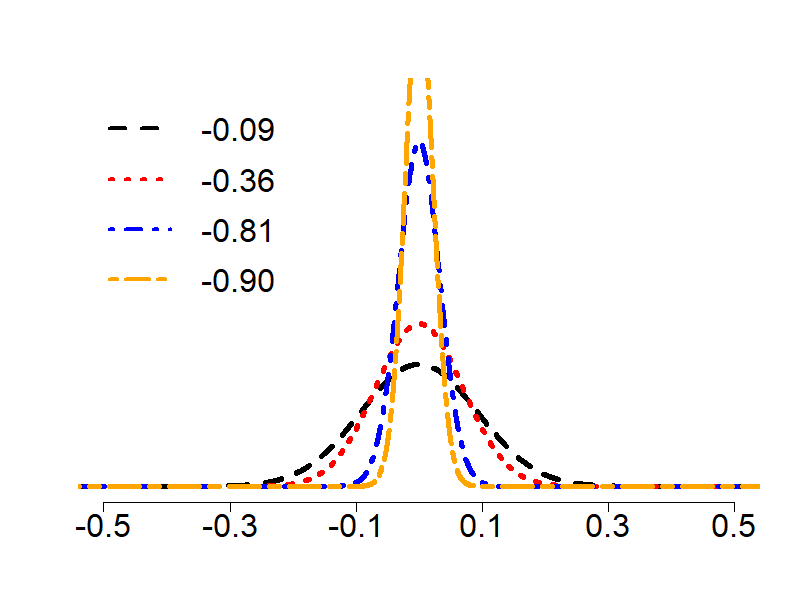}}\hfil
\subfloat[\scriptsize $d_s(\widehat{c}_{12}^x)$, AR and ARMA]{\includegraphics[width=5cm]{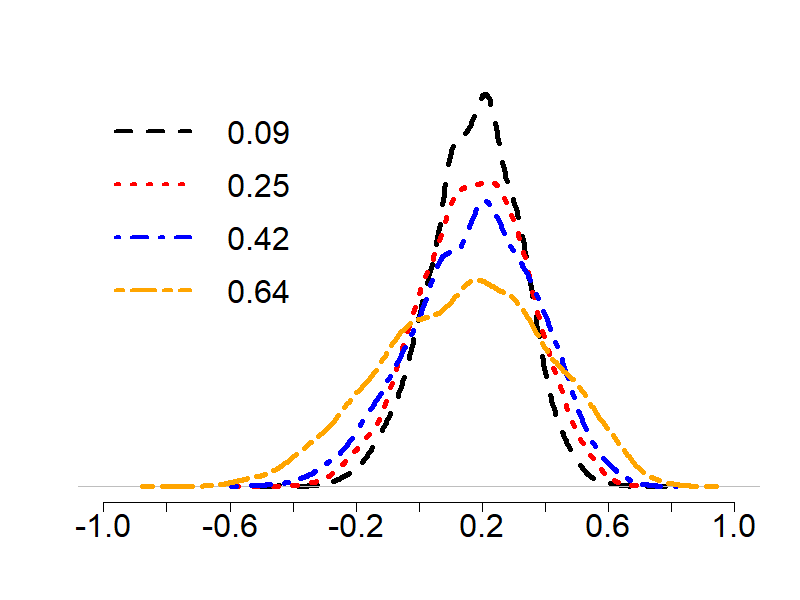}}

\caption{\footnotesize 
Monte Carlo 
densities for
$\widehat{c}_{12}^x$ (a)
and corresponding $\mathcal{D}_{\widehat{c}_{12}^x}$ (b), 
for $T$=100 and various (negative) 
$\phi_{12}$'s. Monte Carlo
densities for $\widehat{c}_{12}^x$
in the case of Laplace weakly correlated AR(3) and ARMA(2,1) processes, for $T=100$ and various 
$\phi$'s (c).}
\label{negSignPhi}
\end{figure}

To generalize 
our findings to the case of non-Gaussian weakly correlated AR and ARMA processes,
we generate covariates according to the following DGPs:
${x}_{1t}=(\phi+0.1){x}_{1t-1}+(\phi+0.1){x}_{1t-2}-0.2{x}_{1t-3}+u_{1t}$, and ${x}_{2t}=\phi x_{2t-1}+\phi x_{2t-2}+{u}_{2t}+0.8{u}_{2t-1}$, where $t=1,
\ldots,100$ and $\phi=0.15, 0.3, 0.45, 0.475$. Moreover, we generate $u_{1t}$ and $u_{2t}$ 
from a bivariate Laplace distribution with means $0$, 
variances $1$, and $c_{12}^u=0.2$.
In these more general cases, we do not know an approximate 
theoretical density 
for $\widehat{c}_{12}^u$.
Therefore, we rely entirely on 
simulations to show the effect of serial dependence on $\Pr\left\{|\widehat{c}_{12}^x|\geq\tau\right\}$. Plot (c) of Figure \ref{negSignPhi} shows $d_s(\widehat{c}_{12}^x)$ obtained 
from 5000 Monte Carlo 
replications 
for the different values of $\phi$. In short, also in 
the more general 
cases where covariates are non-Gaussian, weakly correlated AR(3) and ARMA(2,1) processes, the probability of getting large sample cross-correlations 
depends on the degree of serial dependence. 
More simulation results are provided in Supplement \ref{MoreCases}.

\section{\setlength{\baselineskip}{0.5\baselineskip}A Remedy for Serial Dependence-Induced Spurious Correlation}

In this section we propose 
an approach to mitigate the issues caused by serial dependence-induced spurious correlations 
for the 
performance of PRs. Our proposal consists 
of a two-step procedure. In the first step, we estimate a univariate model on each covariate time series (for example, an ARMA model); in the second step, we
run PRs 
using the residuals of the models fitted
in the first step instead of the original covariates. 
In more detail, let $x_{it}$ 
(the $i$-th
time series 
at 
time $t$)
be generated by the model
\begin{equation}\label{filter}
x_{it}=\sum_{l=1}^{p_i}\phi_{il}x_{it-l}+\sum_{k=1}^{q_i}\theta_{ik}u_{it-k}+u_{it}  
\end{equation}
where $i=1,\dots,n, t=1\dots,T $. This describes an ARMA($p_i$,$q_i$) process where $p_i$ is the order of autocorrelation, which
determines the order of the weighted moving average over past values of the covariate, and $q_i$ is the order of 
the weighted moving average over past errors. 
Note that the AR (i.e.~$p_i\geq1$, $q_i=0$) and MA (i.e.~$p_i=0$, $q_i\geq1$) models 
are special cases of \eqref{filter}. For 
notational simplicity let $x_{it|t-1}=\sum_{l=1}^{p_i}\phi_{il}x_{it-l}+\sum_{k=1}^{q_i}u_{it-k}$
and let $\widehat{x}_{it|t-1}$ be an estimate of $x_{it|t-1}$.
We propose to 
run PRs 
using the estimated residuals $\widehat{u}_{it}=x_{it}-\widehat{x}_{it|t-1}$.

\subsection{The Working Model on ARMA residuals}
Assume that 
response variable and 
covariates are generated by the following DGPs:
\begin{align}
 & y_t=\sum_{i=1}^n\alpha_i^*x_{it-1}+\varepsilon_{t}\ ,
 \label{Mod_xY} \\
 & x_{it}=\phi_ix_{it-1}+u_{it}\ ,  \label{X_AR1} \\
 & \varepsilon_{t}=\phi_{\varepsilon}\varepsilon_{t-1}+\omega_{t} \ ,\label{noise_AR1}
\end{align}
where $i=1,\dots,n$, $t=1,\dots,T$,
$|\phi_i|<1$ 
$|\phi_{\varepsilon}|<1$, and $u_{it}$ and $\omega_{t}$ are the $i.i.d.$ random 
errors of the processes. 
Consider a regression model where one lag of $y_t$ is included in the set of potential predictors
\begin{equation}\label{Mod_xxY}
y_t= \sum_{i=1}^n\alpha_i^*x_{it-1} + \phi_yy_{t-1} + \omega_t \ .
\end{equation}
This strategy is usually adopted to eliminate any residual serial correlation \citep{keele2006}. 
The following two assumptions are crucial for our proposal.
\begin{ass}\label{ass:due} $u_{it}\perp u_{jt-l}$ for any $i$, $j$, $t$ and $l\neq0$.\end{ass} 
\begin{ass}\label{ass:tre} $u_{it-l}\perp\omega_t$ for any $i$, $t$ and $l$.\end{ass}
\noindent Moreover, we assume that $n$ is comparable to or larger than $T$ and,
temporarily and for the sake of the argument, that
the $u_{it-1}$'s
are 
observable --
so that we do not need to estimate them through $x_{it}-\widehat{x}_{it|t-1}$. Our proposal utilizes the
{\it working model}
%
\begin{equation}\label{Mod_uY}
y_t= \sum_{i=1}^n\alpha_i^*u_{it-1} + \phi_yy_{t-1} + \omega_t 
\end{equation}
in place of model~\eqref{Mod_xxY}. In the following, we refer to PRs when estimating the coefficients of
\eqref{Mod_xxY} and to 
$u$-PRs when estimating the coefficients of
\eqref{Mod_uY}. 
Taking a step back, 
we consider the asymptotic behavior of the unpenalized OLS least squares estimates 
for the two models. The reason for this digression is that, in the presence of serial dependence 
in the covariates and 
in the error,
OLS estimates 
for \eqref{Mod_xxY} are asymptotically biased \citep{Achen2000,keele2006}. 
We show that those 
for \eqref{Mod_uY} are not.
The following  Proposition contrasts the asymptotic behavior of the two OLS estimates in terms of convergence in probability.

\begin{prop}\label{prop:quattro}
Let 
the data be generated by \eqref{Mod_xY}, \eqref{X_AR1} and \eqref{noise_AR1},
with $\phi_i=\phi$, $i=1,\dots,n$. Moreover, let $\sigma_y^2$ be the variance of $y$ and $R^2={\pmb{\alpha}^*}'\mathbf{C}_x\pmb{\alpha}^*/\sigma_y^2$ be the asymptotic coefficient of determination, 
where $\mathbf{X}\mathbf{X}'/T\ \overset{p}{\to}\ \mathbf{C}_x$. 
Under Assumptions \ref{ass:due} and \ref{ass:tre}, as $T \to \infty$,
OLS estimates 
for model~\eqref{Mod_xxY} converge as follows
$$\widehat{\phi}_{y}\ \overset{p}{\to}\ \cfrac{\phi_{\varepsilon}(1-R^2)}{1-\phi^2R^2}\ \ \ ,\ \ \ \widehat{\pmb{\alpha}}\ \overset{p}{\to}\ \pmb{\alpha}^*\left(1-\cfrac{\phi\phi_{\varepsilon}(1-R^2)}{1-\phi^2R^2}\right)$$
while 
OLS estimates 
for model~\eqref{Mod_uY} converge as follows
$$\widehat{\phi}_{y}\ \overset{p}{\to}\ \phi R^2+\phi_{\varepsilon}(1-R^2)\ \ \ ,\ \ \  \widehat{\pmb{\alpha}}\ \overset{p}{\to}\ \pmb{\alpha}^*.$$
\end{prop}

\noindent\textbf{Proof:} See Supplement \ref{Proof_Prop4}

\vspace{0.3in}

\noindent Proposition \ref{prop:quattro} shows that applying OLS to the covariates induces an asymptotic bias 
in the estimation of $\pmb{\alpha}^*$ which does not occur when applying OLS to the residuals. The same problem will be inherited by penalized versions of the least squares. Thus, 
we can articulate the impact of serial dependence on the coefficient estimation error of model \eqref{Mod_xxY} 
in relation to $T$. When $T$ is fixed, 
serial dependence leads both to biased coefficients for the relevant variables and to the proliferation of false positives. The latter may be the consequence of serial dependence-induced spurious correlations between relevant and irrelevant variables \citep{Fan08,Fan13}.
When $T$ grows and the issues caused by serial dependence-induced spurious correlations fade, 
we still have a negative effect of serial dependence on the estimated coefficients, as shown in Proposition \ref{prop:quattro} (in Supplement \ref{uOLS} we compare OLS applied on ARMA residuals with some of the best-known OLS estimators used to address serial dependence in regression). Note that, regardless of whether $T<<\infty$ or $T\rightarrow\infty$,
the estimation error of $\widehat{u}$-PRs is smaller than that of PRs.  


We 
point out that the working model \eqref{Mod_xxY} never corresponds to the ``true model'' \eqref{Mod_xY} (i.e.~the data generating process for the response), and its estimates are downwardly biased as both $\phi$ and $\phi_{\varepsilon}$ increase. Thus, 
we proceed by comparing \eqref{Mod_xY} and 
\eqref{Mod_uY}.
The latter allows us to estimate the true vector $\pmb{\alpha}^*$ through the $u_{it-1}$'s, regardless of the possible issues in estimating the serial dependence of $y_t$. This is possible because $u_{it}\perp x_{it|t-1}$ for any specification of $x_{it|t-1}$, which is a consequence of Assumptions \ref{ass:due} and \ref{ass:tre}. For this reason,
omitting autoregressive and/or
moving average component(s) from 
\eqref{Mod_uY} 
does not induce an omitted-variables bias 
for the estimated coefficients. However, 
it does 
reduce 
the explained variance of $y_t$ -- which is mitigated by including its lags among the regressors of the working model. 
How well the ``true model'' \eqref{Mod_xY} and the working model \eqref{Mod_uY} match depends on the DGPs of the covariates and the
noise. We show this in Examples~\ref{ex:uno}
and \ref{ex:due}. 


\vspace{0.3cm}

\begin{example}\label{ex:uno}
(Equal degrees of serial dependence). 
Suppose $\phi_i=\phi_{\varepsilon}=\phi$, $i=1,\dots,n$. Then, model \eqref{Mod_xY} can be rewritten as 
$$y_t  = \sum_{i=1}^n\alpha_i^*(\phi x_{it-2}+u_{it-1})+\phi\varepsilon_{t-1}+\omega_{it} = \sum_{i=1}^n\alpha_i^*u_{1t-1}+\phi y_{t-1} + \omega_t\ .$$
Thus, in an \enquote{ideal regime} in terms of degree of serial dependence (also known as \enquote{common factor restriction}), the working model \eqref{Mod_uY} is equivalent to  the true model \eqref{Mod_xY} because of the decomposition of the AR(1) processes $x_{1t-1}$, $x_{2t-1}$ and $\varepsilon_{t}$. Note that by Proposition \ref{prop:quattro}, if the common factor restriction holds, $\widehat{\phi}_{y}\ \overset{p}{\to}\ \phi$.
\end{example}


\begin{example}\label{ex:due}
(Different degrees of serial dependence). 
Suppose $\phi_1\neq\dots\neq\phi_n\neq\phi_{\varepsilon}$. Then, with some simple steps, model \eqref{Mod_xY} can be rewritten as 
$$y_t  = \sum_{i=1}^n\alpha_i^*(\phi_i x_{it-2}+u_{it-1})+\phi_{\varepsilon}\varepsilon_{t-1}+\omega_{it} =   \sum_{i=1}^n\alpha_i^*u_{it-1}+\sum_{i=1}^n\alpha_i^*(\phi_i x_{it-2})+ \phi_{\varepsilon}\varepsilon_{t-1} + \omega_t\ .$$
Thus, in this perhaps more realistic regime, the working model \eqref{Mod_uY} is {\em not} equivalent to the true model \eqref{Mod_xY} since the predictors and the error do not have the same degree of serial dependence, and therefore the use of $y_{t-1}$ does not allow us to summarize the serial dependence of $y_t$. 
\end{example}
 Two more examples (equal degrees of serial dependence and different models, either for the predictors or for the error) are provided in Supplement \ref{Examples}.
We note that, even when 
true and 
working models do not match, as in 
Example~\ref{ex:due}, 
$u$-PRs 
moves us from estimating coefficients in a context characterized by high spurious correlation, to one characterized by very weak (or absent) spurious correlations. 
Of course the parameters we estimate about 
the past of $y_t$ 
change, but we can still formulate an effective estimation strategy, e.g.,
if $\varepsilon_{t}=\sum_{j=1}^{p_{\varepsilon}}\phi_{\varepsilon j}\varepsilon_{t-j}+\omega_{t}$, 
even in cases such as 
Example~\ref{ex:due}, 
the variability due to
misspecification of the 
serial dependence of $y_t$ 
is less than that introduced by estimating the model directly on the $\mathbf{x}_i$'s, $i=1,\dots,n$. 

Since the error term is not correlated with the regressors included in the model, its serial dependence does not violate the assumption of exogeneity and the OLS estimator remains unbiased and consistent. To prevent serial dependence of the error term, a possible solution could be to increase the number of lags of $y_t$ considered in the working model.


Of course, in practice, the $\pmb{u}_{i}$'s, $i=1,\dots,n$,
are not observable and need to be replaced by estimated residuals of 
ARMA, 
AR 
or MA processes. When fitting the working models with such residuals we refer to our proposal as $\widehat{u}$-PRs.
%
In the following Sections we demonstrate the estimation and forecasting performance of $\widehat{u}$-LASSO 
(LASSO applied on ARMA residuals) 
through both simulations and an empirical application. 

\subsection{$\widehat{u}$-LASSO}

\subsubsection{Coefficient Estimation Error Bound}
\begin{figure}[H]
\graphicspath{{images/}}
\centering
{\includegraphics[width=5.5cm]{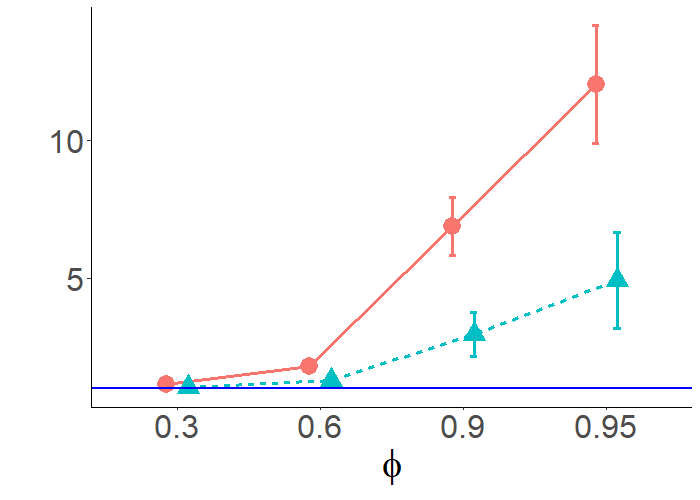}}\hfil
{\includegraphics[width=5.5cm]{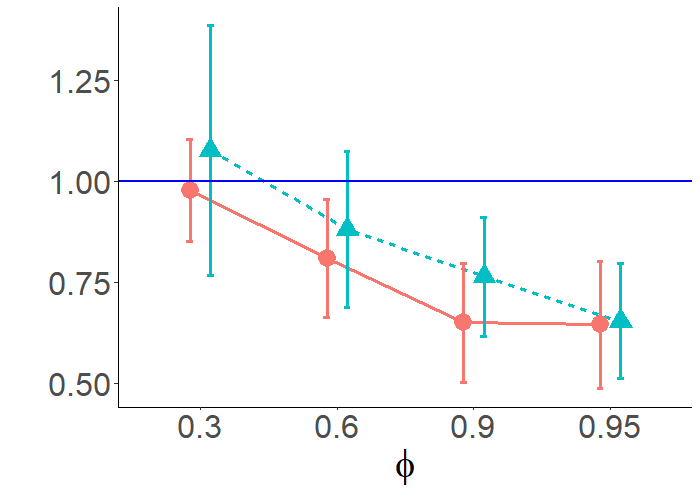}}
\caption{\footnotesize Average and standard deviation of $\widehat{\psi}_{min}^{\widehat{u}}/\widehat{\psi}_{min}^{x}$ (left)
and of $\scriptsize||\widehat{\pmb{\alpha}}_{\widehat{u}}-\pmb{\alpha}^*||_2/||\widehat{\pmb{\alpha}}_x-\pmb{\alpha}^*||_2$ (right)
across 1000 Monte Carlo replications, for various values of $\phi$. 
Orange circles and bars
represent 
means and standard deviations for $n/T=0.5$, 
blue triangles and bars 
represent means and standard deviations for $n/T=1.5$. 
In each panel, 
the horizontal blue line marks a ratio of 1.}\label{Fig:CoefEst}
\end{figure}

Here, we present 
Monte Carlo experiments to assess 
the effectiveness of $\widehat{u}$-LASSO in reducing the coefficient estimation error. We generate the response as
$y_t=\sum_{i=1}^n\alpha_i^*x_{it-1}+\varepsilon_t,$ 
where $\varepsilon_t=\phi\varepsilon_{t-1}+\omega_t$, and $\omega_t\sim N(0,\sigma_{\omega}^2)$. The coefficient vector $\pmb{\alpha}^*=(\alpha_1^*,\dots,\alpha_n^*)'$ is sparse with $||\pmb{\alpha}^*||_0=10$.
The active covariates 
are the first 10, followed by 
$n-10$ inactive ones, and $\alpha_1=\dots=\alpha_{10}=1$. 
We generate the $n$ covariates
as $x_{it}=\phi x_{it-1}+u_{it}$, $i=1,\dots,n$, $t=1,\dots,T$ with $T=100$, where $u_{it}\sim N(0,1)$ and $\left(\mathbf{C}_u\right)_{ij}=c_{ij}^u=0.3^{|i-j|}$. We consider $n=50,
150$ and $\phi=0.3, 0.6, 0.9, 0.95$. 
Left panel of Figure~\ref{Fig:CoefEst} 
displays 
mean and 
standard deviation of the ratio between $\widehat{\psi}_{min}^{\widehat{u}}$ (the minimum eigenvalue of $\widehat{\mathbf{C}}_{\widehat{u}}$) and $\widehat{\psi}_{min}^x$ obtained 
from 1000 Monte Carlo simulations run 
with $n=50$ (orange circles) and $150$ (blu triangles);
for $n=150$ we consider the minimum eigenvalues of the correlation matrices
restricted to the 10 relevant variables. As expected, the correlation matrix of the $\widehat{\mathbf{u}}_i$'s, $i=1,\dots,n$ does not suffer from spurious correlation induced by serial dependence, and this leads to an increment of $\widehat{\psi}_{min}^{\widehat{u}}/\widehat{\psi}_{min}^{x}$ as $\phi$ increases. To observe how this result translates into coefficients estimation accuracy, we compare the estimation error of LASSO $\left(||\widehat{\pmb{\alpha}}_x-\pmb{\alpha}^*||_2\right)$ with that of $\widehat{u}$-LASSO $\left(||\widehat{\pmb{\alpha}}_{\widehat{u}}-\pmb{\alpha}^*||_2\right)$, where the tuning parameter $\lambda$ is selected 
by BIC. Right panel of Figure~\ref{Fig:CoefEst}
shows how the mean and standard deviation of $||\widehat{\pmb{\alpha}}_{\widehat{u}}-\pmb{\alpha}^*||_2/||\widehat{\pmb{\alpha}}_x-\pmb{\alpha}^*||_2$
vary as a function of $\phi$.
Also here, as expected, the application of LASSO on serially uncorrelated data reduces the coefficient estimation error, 
with a gain in estimation accuracy that increases with $\phi$. 
Of course the gains in accuracy shown in right panel of Figure \ref{Fig:CoefEst} may in part be due to a reduction in bias, as the LASSO inherits the OLS bias illustrated 
in Proposition \ref{prop:quattro}. However, bias is not the whole story here; the gains in accuracy are also linked to overcoming 
the spurious correlation induced by serial dependence. 
In this regard, Table \ref{Tab:TPFP} reports 
average 
percentages of true and false positives (\%TP, \%FP) 
with LASSO and $\widehat{u}$-LASSO; 
the latter clearly improves variable selection.
\cite{Fan08} and \cite{Fan13} pointed out the role of the spurious correlations between relevant and irrelevant variables in the proliferation of false positives;
the improved variable selection performance of $\widehat{u}$-LASSO compared to LASSO can be interpreted as partial evidence of the fact that, concurrently with the bias illustrated in Proposition \ref{prop:quattro}, 
spurious correlation induced by serial dependence strongly contributes to false positives.

In summary, results shown in Figure \ref{Fig:CoefEst} and Table \ref{Tab:TPFP} corroborate the theoretical analysis according to which an increase in the degree of serial dependence leads to an increase in the probability of large spurious correlations, which in turn increases the probability of a small minimum eigenvalue for the sample correlation matrix. This negatively affects the estimation accuracy of PRs (see Proposition \ref{prop:uno}).
%
\begin{table}[H]
\centering
\caption{\footnotesize Average 
percentages of true positives (\%TP) and false positives (\%FP) by LASSO and $\widehat{u}$-LASSO across 1000 Monte Carlo replications, for various values of
$n/T$ and $\phi$.
}\label{Tab:TPFP}
\scalebox{0.65}{
\begin{tabular}{c@{\hspace{0.8cm}}cc@{\hspace{0.8cm}}c@{\hspace{0.8cm}}cc@{\hspace{0.8cm}}ccc@{\hspace{0.8cm}}c@{\hspace{0.8cm}}cc@{\hspace{0.8cm}}c}\hline
$\phi$  && \multicolumn{5}{c}{$n/T=0.5$} & & \multicolumn{5}{c}{$n/T=1.5$}  \\\cline{3-7} \cline{9-13}
&  & \multicolumn{2}{c}{LASSO} && \multicolumn{2}{c}{$\widehat{u}$-LASSO} && \multicolumn{2}{c}{LASSO} && \multicolumn{2}{c}{$\widehat{u}$-LASSO}  \\\cline{3-4} \cline{6-7} \cline{9-10} \cline{12-13}
&&\% TP &\% FP &&\% TP &\% FP &&\% TP &\% FP &&\% TP &\% FP \\\cline{3-3}\cline{4-4}\cline{6-6}\cline{7-7}\cline{9-9}\cline{10-10}\cline{12-12}\cline{13-13}
0.3&& 100.0 &13.9 &&100.0 &13.0 &&100.0 &40.4 &&100.0 &48.7\\
0.6&& 100.0 &25.7 &&100.0 &12.8 &&100.0 &63.5 &&100.0 &28.7\\
0.9&& 99.9  &45.2 &&100.0 &13.7 &&99.8 &55.7 &&100.0 &20.3\\
0.95&& 99.9 &47.5 &&100.0 &13.9 &&99.2 &43.8 &&100.0 &19.9\\
\hline
\end{tabular}
}
\begin{tablenotes}
\begin{spacing}{1}
      \footnotesize
      \item 
      Note: \%TP=$\frac{1}{1000}\sum_{k=1}^{1000}\text{TP}_k/10$
      and \%FP=$\frac{1}{1000}\sum_{k=1}^{1000}\text{FP}_k/(n-10)$ where, for each Monte Carlo replication $k$, 
      TP$_k=\#\left \{\widehat{\alpha}_j\neq0 : j\in S \right \}_k$, 
FP$_k=\#\left \{\widehat{\alpha}_j\neq0 : j\in S^c \right \}_k$,
$S$=Supp$(\pmb{\alpha}^*)$, 
$S^c=$ 
complement.

\normalsize
 \end{spacing}
   \end{tablenotes}
\end{table}

\subsubsection{Empirical Application}\label{EmpApp}
We consider 
Euro Area (EA) data composed by 309 monthly macroeconomic time series spanning the period between January 1997 and December 2018. 
The series are listed in Supplement \ref{Stat_EA}, grouped according to their measurement domain: Industry \& Construction Survey (ICS), Consumer Confidence Indicators (CCI), Money \& Interest Rates (M\&IR), Industrial Production (IP), Harmonized Consumer Price Index (HCPI), Producer Price Index (PPI), Turnover \& Retail Sale (TO), Harmonized Unemployment Rate (HUR), and Service Surveys (SI). Supplement \ref{Stat_EA} also reports
transformations applied to 
the series to achieve
stationarity
(we did not attempt to identify or remove
outliers).
The target variable is the Overall EA Consumer Price 
Index (CPI), which is transformed as I(2)
(i.e.~integration of order 2)
following \cite{SeW2002b}:
$$\small y_{t+h}=(1200/h)\text{log}(CPI_{t+h}/CPI_t)-1200\text{log}(CPI_t/CPI_{t-1}) \ ,$$
where $y_t=1200\text{log}(CPI_t/CPI_{t-1})-1200\text{log}(CPI_{t-1}/CPI_{t-2})$, and $h$ is the forecasting horizon. We compute forecasts of $y_{t+h}$ at horizons $h=12$
and $24$ using a rolling $\omega$\--year window $[t-\omega, t+1]$;
the models are re-estimated at each $t$, adding one observation on the right of the window and removing  one observation on the left.
The last forecast is 
December 2018.
The
methods employed for 
our empirical exercise
are:
\begin{itemize}
\item {\it Univariate AR($p$):} the autoregressive forecasting model 
based on $p$ lagged 
values of the target variable, i.e.~$\widehat{y}_{t+h}=\widehat{\alpha}_0+\sum_{i=1}^p\widehat{\phi}_iy_{t-i+1}$,
which serves as a benchmark.
\item {\it LASSO}
\citep{tibshirani96}:
forecasts are obtained from the equation
$\widehat{y}_{t+h}^{x}=\widehat{\alpha}_0+\widehat{\pmb{\beta}}_x'\mathbf{v}_{t},$
where $\widehat{\pmb{\beta}}_x=(\widehat{\phi}_{1},\dots, \widehat{\phi}_{12}, \widehat{\pmb{\alpha}}_x')'$ is the sparse vector of penalized regression coefficients estimated by the LASSO on the original time series, and $\mathbf{v}_{t}=(y_{t},\dots, y_{t-11},\mathbf{x}'_{t})'$.

\item {\it ${\widehat{u}}$-LASSO}: our proposal, where LASSO is applied to the estimated ARMA residuals.~
Forecasts are obtained from the equation
$\widehat{y}_{t+h}^{x}=\widehat{\alpha}_0+\widehat{\pmb{\beta}}_{\widehat{u}}'\mathbf{w}_{t},$
where $\widehat{\pmb{\beta}}_{\widehat{u}}=(\widehat{\phi}_{1},\dots, \widehat{\phi}_{12}, \widehat{\pmb{\alpha}}_{\widehat{u}}')'$ is the sparse vector of penalized regression coefficients estimated by the LASSO on the estimated ARMA residuals, and $\mathbf{w}_{t}=(y_{t},\dots, y_{t-11}, \widehat{\mathbf{u}}'_{t})'$.
\end{itemize}
For the AR($p$) benchmark the lag order $p$ is selected 
by BIC within $0\leq p\leq 12$.
For the $\widehat{u}$-LASSO, 
estimated residuals are obtained 
filtering each time series with an ARMA($p_i,q_i$), where $p_i$ and $q_i$ are selected 
by BIC
within $0 \leq p_i,q_i\leq 12$, $i=1,\dots,n$.
The shrinkage parameter $\lambda$ of LASSO and $\widehat{u}$-LASSO is selected with BIC and 10-folds cross-validation (CV). Forecasting accuracy 
for 
all three methods is evaluated using the root mean square forecast error (RMSFE), 
defined as
\begin{equation*}\small 
 RMSFE = \sqrt{\frac{1}{T_1 - T_0}\sum_{\tau=T_0}^{T_1}\Bigl( \widehat{y}_{\tau} - y_{\tau}\Bigr)^2}
 \end{equation*}
where $T_0$ and $T_1$ are the first and last 
time points used for the out of sample evaluation.
For LASSO and $\widehat{u}$-LASSO we also consider the number of selected variables. 

Table \ref{Tab:RMSFE} reports ratios of 
RMSFEs between pairs of methods,
as well as 
significance of the corresponding Diebold-Mariano test \citep{DiebMar1995}. We also report the ratio between the average number of selected variables with $\widehat{u}$-LASSO 
and with LASSO.
Notably, $\widehat{u}$-LASSO 
produces significantly better forecasts 
than both the classical LASSO and the AR($p$),
and 
provides a more parsimonious model than the LASSO; 
the ratio between the average number of selected variables is much smaller than $<$1.
This 
is, in principle, consistent with the theoretical analysis we provided earlier. 
The sparser $\widehat{u}$-LASSO output may be due 
to fewer false positives, as compared to the LASSO -- since the latter suffers from the effects of spurious correlations 
induced by serial dependence.
However, since in this real data application we do not know the true DGP, any comments
regarding 
accuracy in variable selection 
is necessarily speculative.

\begin{figure}[H]
\graphicspath{{images/}}
\centering
\subfloat[LASSO, $h$=12]{\includegraphics[width=6.5cm]{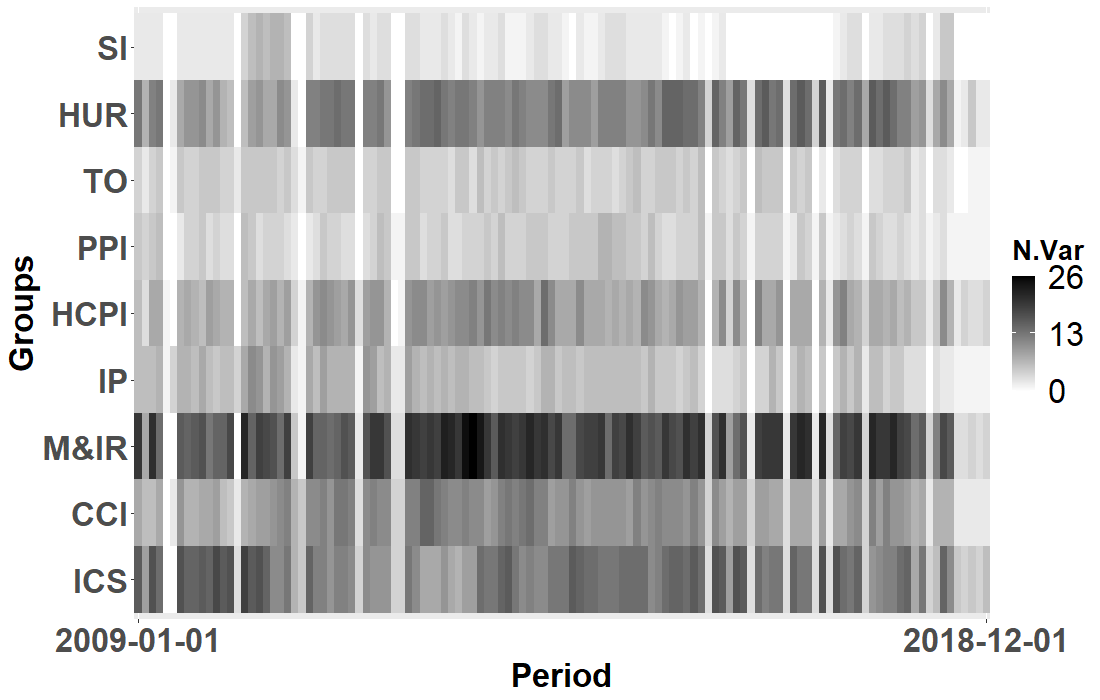}}\hfil
\subfloat[$\widehat{u}$-LASSO, $h$=12]{\includegraphics[width=6.5cm]{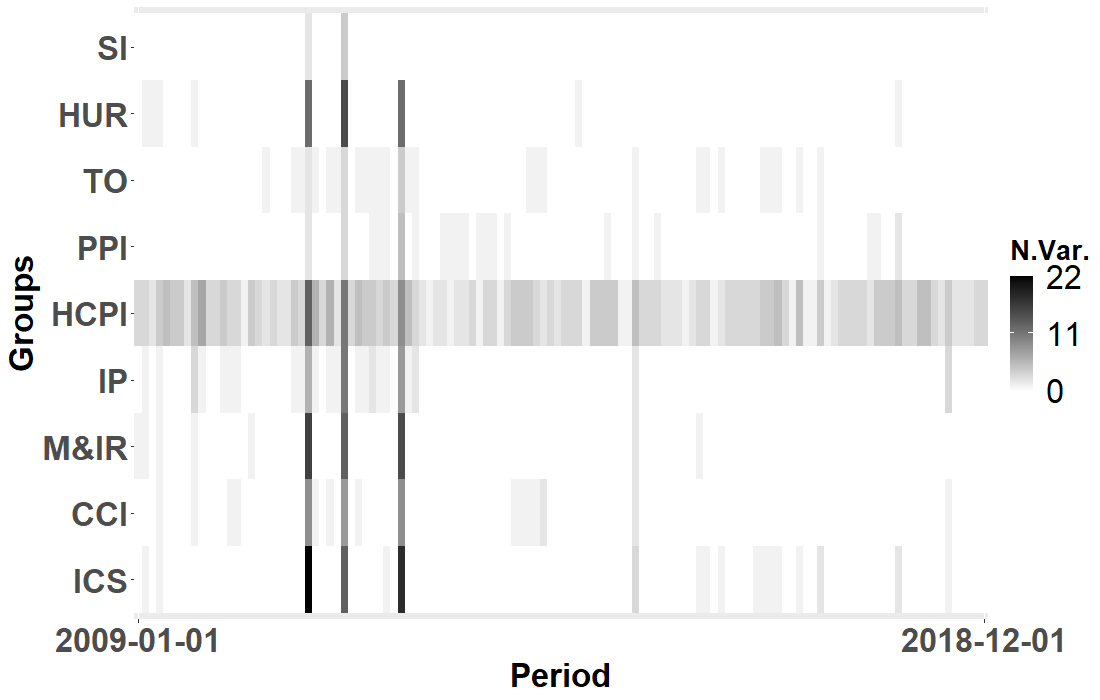}}

\subfloat[LASSO, $h$=24]{\includegraphics[width=6.5cm]{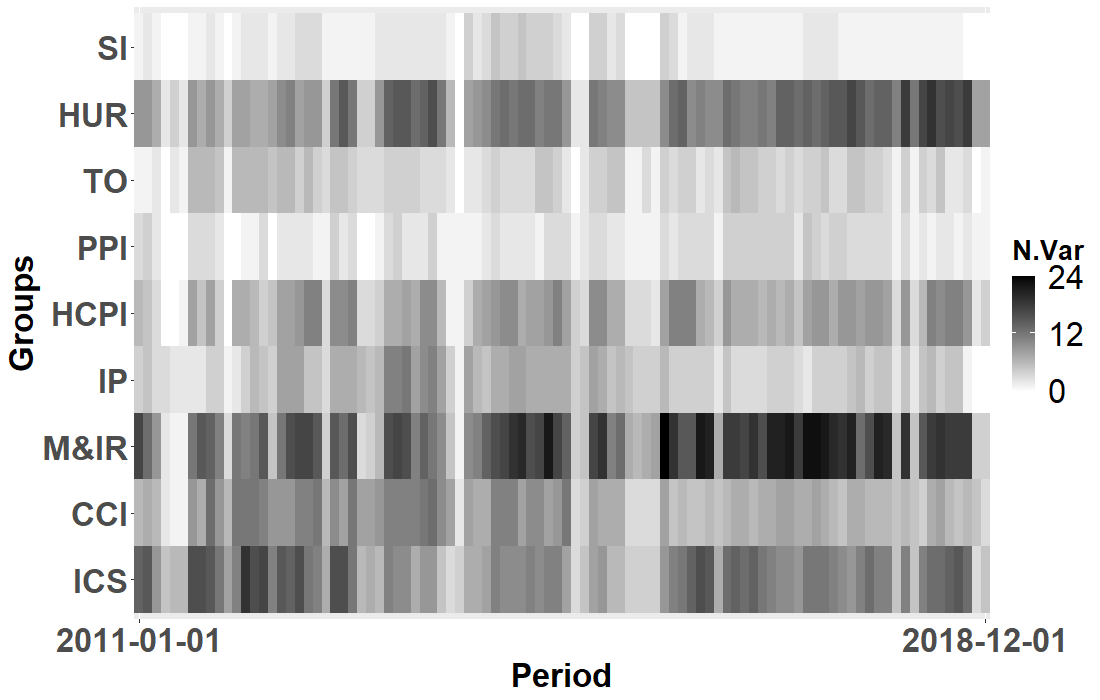}}\hfil
\subfloat[$\widehat{u}$-LASSO, $h$=24]{\includegraphics[width=6.5cm]{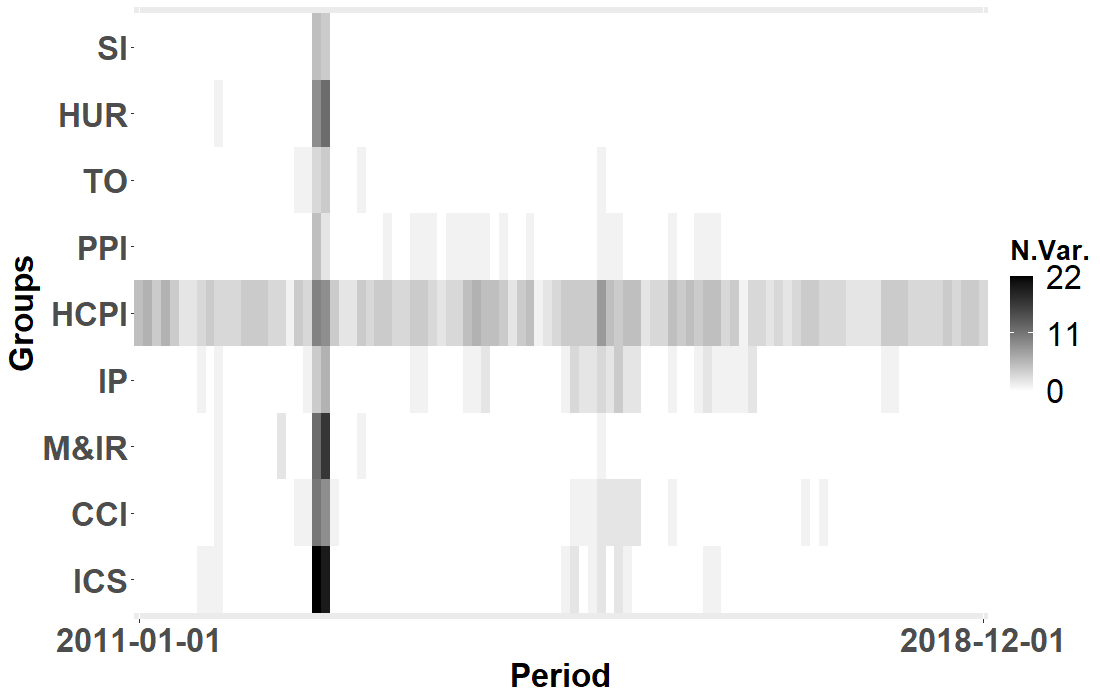}}

\caption{\footnotesize Heatmaps representing the number of 
variables selected by LASSO (left)
and $\widehat{u}$-LASSO (right)
in the nine main domains. The tuning procedure is BIC.
}\label{Heatmaps}
\end{figure}

\begin{table}[H]
\centering
\caption{\footnotesize 
Left: ratios of RMSFE contrasting pairs of employed methods; for each ratio 
we perform a Diebold-Mariano test 
(alternative: 
the second method is less accurate in forecasting) and report p-values as
0 '***'  0.001 '**' 0.01 '*' 0.05 '$^{\sbullet}$' 0.1''. 
Right: average of the number of variables selected 
by $\hat{u}$-LASSO (left of the vertical bar) and LASSO (right of the vertical bar).
}\label{Tab:RMSFE}
\scalebox{0.65}{
\begin{tabular}{l@{\hspace{0.8cm}}ll@{\hspace{0.8cm}}ll@{\hspace{0.8cm}}l@{\hspace{0.8cm}}ll@{\hspace{0.8cm}}lll@{\hspace{0.8cm}}l@{\hspace{0.8cm}}ll@{\hspace{0.8cm}}l}\hline
Method 1 & & Method 2 & & \multicolumn{5}{c}{RMSFE (ratio)} & & \multicolumn{5}{c}{Average of Selected Variables}  \\\cline{1-1} \cline{3-3} \cline{5-9} \cline{11-15}
& &  & & \multicolumn{2}{c}{$h$=12} && \multicolumn{2}{c}{$h=24$} && \multicolumn{2}{c}{$h$=12} && \multicolumn{2}{c}{$h=24$}  \\\cline{5-6} \cline{8-9} \cline{11-12} \cline{14-15}
&&&&BIC&CV&&BIC&CV&&BIC&CV&&BIC&CV\\\cline{5-5}\cline{6-6}\cline{8-8}\cline{9-9}\cline{11-11}\cline{12-12}\cline{14-14}\cline{15-15}
$\widehat{u}$-LASSO && LASSO   & & 0.69***  &  0.66*** &&  0.82*$ \ \ $    & 0.83**$ \ $ &&6.0$|$67.9 &14.8$|$56.8 &&6.2$|$60.9 &14.7$|$57.9\\
$\widehat{u}$-LASSO && AR($p$) & & 0.94$ \ \ \ $   & 0.91*$ \ \ $ &&   0.89*$ \ \ $    & 0.88**    &&--&--&&--&--\\
LASSO && AR($p$)     &    & 1.36$ \ \ \ $  &1.38$ \ \ \ $  && 1.08$ \ \ \ $     &  1.07$ \ \ \ $ &&--&--&&--&--  \\
\hline
\end{tabular}
}
\begin{tablenotes}
\begin{spacing}{1}
      \footnotesize
      \item Note: For 
      AR($p$) 
      coefficients are estimated using the R package {\it lm}. For 
      $\widehat{u}$-LASSO
      estimated residuals are obtained by means of an ARMA($p_i,q_i$) filter. 
      The penalty parameter $\lambda$ 
      is selected with BIC using the R package {\it HDeconometrics}, 
      and with 10-folds cross-validation (CV) 
      using the R package {\it glmnet}. \normalsize
 \end{spacing}
   \end{tablenotes}
\end{table}
\vspace{-0.3in}
In terms of selected variables, Figure \ref{Heatmaps} summarizes patterns over time obtained with BIC tuning (results obtained with CV are reported in Supplement \ref{HeatmapsCV}). The heatmaps represent the number of selected variables categorized according to the nine main domains (see above). LASSO selects variables largely, though not exclusively, from the domains ICS, M\&IR and HUR.
$\widehat{u}$-LASSO is more targeted, selecting variables in the HCPI domain.
The top 5 variables in terms of selection frequency across forecasting samples are listed in Table \ref{Tab:VarSel}. Regardless of tuning (BIC, CV) and forecasting horizon $h$, the top predictor for $\widehat{u}$-LASSO is the Goods Index. The other top predictors, also in the HCPI domain, include EA measurements (e.g.,~Services Index), or are specific to France and Germany (e.g.,~All-Items).
In summary,
$\widehat{u}$-LASSO exploits cross-sectional information mainly focusing on prices, and accrues a forecasting advantage -- as LASSO uses many more variables to produce significantly worse forecasts.
\begin{table}[H]
\caption{\footnotesize 
Five most frequently selected variables.
Selection percentages are ratios between the number of times a variable appears in a forecast and total number of forecasts (120 for $h$=12 and 96 for $h=24$). 
}\label{Tab:VarSel}
\scalebox{0.55}{
\begin{tabular}{l@{\hspace{0.8cm}}ll@{\hspace{0.8cm}}l@{\hspace{0.8cm}}ll@{\hspace{0.8cm}}l}\hline
Rank & & \multicolumn{5}{c}{Selected Variables} \\\cline{1-1} \cline{3-7}
& & \multicolumn{2}{c}{$h$=12} && \multicolumn{2}{c}{$h=24$} \\\cline{3-4} \cline{6-7} 
&&BIC&CV&&BIC&CV\\\cline{3-3}\cline{4-4}\cline{6-6}\cline{7-7}
I°  & &Goods, Index &Goods, Index & &Goods, Index &Goods, Index\\
& &85.8\% &80\% & &85.4\% &85.4\% \\
II° & &Industrial Goods, Index &Services, Index & &Services, Index &Services, Index\\
& &47.5\% &77.5\% & &43.8\% &82.3\% \\
III°& &Services, Index &Industrial Goods, Index & &All-Items (De) &All-Items (Fr)\\
& &40.8\% &56.7\% & &35.4\% &62.5\% \\
IV° & &All-Items Excluding Tobacco, Index &All-Items (Fr) & &All-Items Excluding Tobacco, Index &All-Items (De)\\
& &32.5\% &40\% & &32.3\% &42.7\% \\
V°  & &All-Items (Fr) &All-Items Excluding Tobacco, Index & &Industrial Goods, Index &Industrial Goods, Index\\
& &24.2\% &39.2\% & &30.2\% &35.4\% \\
\hline
\end{tabular}
}
\end{table}

\section{Concluding Remarks}
In this paper we demonstrated that the probability of spurious correlations between stationary orthogonal or weakly cross-correlated processes depends not only on the sample size, but also on the degree of covariates serial dependence. Through this result, we pointed out that serial dependence negatively affects the behavior of the sample cross-correlation matrix, leading to a large probability of getting a small minimum eigenvalue. Considering the role of the minimum eigenvalue 
in the non-asymptotic 
estimation error bounds of PRs, 
our findings 
highlight the limitations of these methods in a time series 
context. In order to improve the estimation performance of PRs in 
such context, we propose an approach based on applying PRs to pre-whitened (ARMA filter) time series. This proposal allows us to solve the problem of large spurious correlation as well as the problem of biased estimates due to the inclusion of response lags \citep{Achen2000}.

We assessed the performance of our proposal through Monte Carlo simulations and an empirical application to Euro Area macroeconomic time series.
Through simulations we observed that $\widehat{u}$-LASSO, i.e.~the LASSO applied on ARMA residuals, reduces the probability of large spurious correlation, performing better than 
LASSO applied on the original covariates in terms of coefficients estimation. 
Through the empirical application
we observed that $\widehat{u}$-LASSO improves the forecasting performance of LASSO, and produces 
more parsimonious 
models. These findings encourage us to further investigate the
potential of 
$\widehat{u}$-PRs 
-- and especially 
sparse 
$\widehat{u}$-PRs.

\bibliographystyle{chicago}
\addcontentsline{toc}{chapter}{Bibliography}
\bibliography{biblio}


\newtheorem{innercustomgeneric}{\customgenericname}
\providecommand{\customgenericname}{}
\newcommand{\newcustomtheorem}[2]{%
  \newenvironment{#1}[1]
  {%
   \renewcommand\customgenericname{#2}%
   \renewcommand\theinnercustomgeneric{##1}%
   \innercustomgeneric
  }
  {\endinnercustomgeneric}
}

\newcustomtheorem{customprop}{Proposition}
\newcustomtheorem{customrem}{Remark}
\newcustomtheorem{customex}{Example}

\newpage

\appendix
\huge{\textbf{Supplement}}
\normalsize
\begin{appendices}

\section{Restricted Eigenvalue}\label{ResEin}

In the specific case of $n>T$ the loss function $\mathcal{L}(\pmb{\alpha})$ cannot be strongly convex since $\mathbf{X}\mathbf{X}'/T$ is not positive definite. In this specific case \cite{Bickel2009} proposed a solution based on a kind of strong convexity for some subset $\mathcal{C}\subset\mathbb{R}^n$ of possible perturbation vectors $\Delta=|\widehat{\pmb{\alpha}}-\pmb{\alpha}^*|\in\mathbb{R}^n$, named {\it Restricted Eigenvalue Condition}. In particular, for any subset $S\subseteq\left\{1,2,\dots,n\right\}$ with cardinality $s$, let $\Delta_S\in\mathbb{R}^{S}$ and $\Delta_{S^c}\in\mathbb{R}^{S^c}$. The restricted eigenvalue condition requires that there exists a positive number $\nu$ such that
$$\min_{\Delta\in\mathbb{R}^n:\Delta\neq0}\frac{\left |\left |\mathbf{X}'\Delta\right |\right |_2}{\sqrt{T}\left |\left |\Delta_S\right |\right |_2}\geq\nu.$$

\noindent Such condition is essentially a restriction on the eigenvalues of $\mathbf{X}\mathbf{X}'/T$ as a function of sparsity, which allows for the strong convexity to hold whit parameter $\gamma=\nu$, which characterizes how strong the covariates depend on each other. According to \cite{Bickel2009}, the restricted eigenvalue condition restricts the LASSO error to a set of the form:
$$\mathcal{C}\left (S\right )\coloneqq\left \{\widehat{\Delta}\in\mathbb{R}^n : ||\widehat{\Delta}_{S^c}||_1\leq 3||\widehat{\Delta}_S||_1 \right \}.$$

\section{\setlength{\baselineskip}{0.5\baselineskip}On the Population Cross-Correlation in Time Series}\label{PopCrosCorr}
Let the response variable be generated according to the following DGP: $y_t=\sum_{i=1}^n\alpha_i^*x_{it-1}+\varepsilon_t$. We consider the case where the covariates are generated as follows
\begin{equation}\label{FactorDGP}
x_{it}=\lambda_iF_t+u_{it},\end{equation}
with $i=1,\dots,n$, $t=1,\dots,T$, where $n$ is comparable to or larger than $T$ and therefore PRs are used in order to estimate $\alpha$'s. $F_t$ represents a common factor that introduce population cross-correlation between covariates, $\lambda_i$ is the factor loading relative to $x_i$, and $u_{it}$ is the idiosyncratic component relative to $x_i$ at time $t$.

In this case, \cite{Fan2020} propose a method to reduce the cross-correlation between covariates in order to improve the estimation accuracy of sparse PRs. It consists in using the principal component analysis to obtain $\widehat{\lambda}_i,\ \widehat{F}_t$ and $\widehat{u}_{it}=x_{it}-\widehat{\lambda}_i\widehat{F}_t$, i.e. estimates of $\lambda_i$, $F_t$ and $u_t$. 
Hence, when covariates are generated by model \eqref{FactorDGP}, the procedure proposed by \cite{Fan2020} allows us to deal with the problem from PRs estimation with highly cross-correlated covariates $\mathbf{x}_1,\dots,\mathbf{x}_n$ to PRs estimation with weakly or orthogonal covariates $\widehat{\mathbf{u}}_1,\dots,\widehat{\mathbf{u}}_n$.

\vspace{0.3cm}

However, we stress that if the idiosyncratic components are orthogonal or weakly cross-correlated AR processes, the methodology proposed by \cite{Fan2020} would not solve the problem of the high spurious cross-correlation caused by serial dependence.

\section{Proofs}

To achieve the theoretical results we generalize the approach of \cite{Anderson03} to our time series context. In particular the equation
$$\cfrac{\sqrt{a_{11}}\ b}{\sqrt{v/(T-2)}}=\sqrt{T-2}\cfrac{a_{12}/\sqrt{a_{11}a_{22}}}{\sqrt{1-a_{12}^2/(a_{11}a_{22})}}=\sqrt{T-2}\cfrac{\widehat{c}_{12}^x}{\sqrt{1-(\widehat{c}_{12}^x)^2}},$$
where $c_{12}^u=0$, $b=a_{21}/a_{11}$ and $v=a_{22}-a_{21}^2/a_{11}$, is obtained by \citet[p.~119]{Anderson03}.

\subsection{Proposition 2}\label{Proof_Prop2}
We first focus on the distribution of the sample covariance between $x_1$ and $x_2$, which is
\begin{multline}
\widehat{Cov}(x_1,x_2)= \frac{a_{12}}{(T-1)}=\\
=\left[\sum_{l=1}^{T-1}\phi_1^l\widehat{Cov}(u_{1[-l]},u_2)+
\sum_{l=1}^{T-1}\phi_2^l\widehat{Cov}(u_{2[-l]},u_1)+\widehat{Cov}(u_1,u_2)\right](1-\phi_1\phi_2)^{-1},
\end{multline}
\noindent where $\widehat{Cov}(u_{i[-l]},u_j)=\sum_{t=l+1}^{T}(u_{it-l}-\overline{u}_i)(u_{jt}-\overline{u}_j)/(T-l-1)$, for $i\neq j=1,2$,
$\overline{u}_i=\frac{1}{T-l-1}\sum_{t=l+1}^{T}u_{it-l}$ and $\overline{u}_j=\frac{1}{T-l-1}\sum_{t=l+1}^{T}u_{jt}$.
\noindent Since $u_1$ and $u_2$ are standard 
Normal, we have (see \citealt{Glen2004} and Supplement \ref{Dist_Cu})
$$\widehat{Cov}(u_1,u_2)\approx  N\left (0\:,\:\cfrac{1}{(T-1)}\right ).$$
\noindent Moreover, the quantity
$$\eta_{12}=\sum_{l=1}^{T-1}\phi_1^l\widehat{Cov}(u_{1[-l]},u_2)+\sum_{l=1}^{T-1}\phi_2^l\widehat{Cov}(u_{2[-l]},u_1),$$
\noindent is a linear combination of the sample covariances between the  residual of a time series at time $t$ and the  lagged residuals of the other time series. Note that $\eta_{12}$ is a linear combination of $N\left(0, \frac{\phi_i^{2l}}{T-l-1}\right)$, $i=1,2$. However, because $|\phi_i|<1$, we can approximate $\eta_{12}$ as a linear combination of
centered Normals with variance $\frac{1}{T-1}$, so that
$$\eta_{12}\approx  N\left(0, \frac{\phi_1^2+\phi_2^2-2\phi_1^2\phi_2^2}{(T-1)(1-\phi_1^2)(1-\phi_2^2)}\right)$$
and
$$\widehat{Cov}(x_1,x_2)\approx  N\left (0\:,\:\frac{1-\phi_1^2\phi_2^2}{(T-1)(1-\phi_1^2)(1-\phi_2^2)(1-\phi_1\phi_2)^2}\right ).$$
%
%
\noindent Since $a_{11}$ is $T-1$ times the sample variance of  $x_1$, 
$a_{11}\approx \frac{T-1}{1-\phi_1^2}$. Therefore, $b=\frac{a_{21}}{a_{11}}$ is Normally distributed and, 
based on the approximation of mean and variance of a ratio (see \citealt{Stuart1998}), we have $E(b)=0$ and
%
\begin{align*}
\begin{split}
Var(b)
=(T-1)^2\frac{Var(a_{12})}{E(a_{11})^2}
& \approx  (T-1)^2Var\left(\widehat{Cov}(x_1,x_2)\right)\left(\frac{1-\phi_1^2}{T-1}\right)^2 \\
 & =\frac{(1-\phi_1^2\phi_2^2)(1-\phi_1^2)}{(T-1)(1-\phi_2^2)(1-\phi_1\phi_2)^2} \ \ .
\end{split}
\end{align*}
 \hfill $\blacksquare$

\begin{customrem}{C.1}
For $T<\infty$ the quantity $\eta_{12}$ has a variance that increases with the degree of serial dependence. This quantity strongly affects the impact of the degree of serial dependence on the variance of $a_{12}$ and, as a consequence, on the variance of both $\widehat{Cov}(x_1,x_2)$ and $b$.
\end{customrem}

\begin{customrem}{C.2}
It is important to note that if $x_1$ and $x_2$ are generated by independent MA($q$) processes then serial dependence increases the variance of $b$ as the order $q$ increases; see \cite{Granger2001}. This is due to the fact that any MA($\infty$) can be represented as an AR(1). Thus, increasing $q$, we are faced with the same spurious component $\eta_{12}$ that impacts the sample covariance between orthogonal AR(1) processes. 
\end{customrem}

\subsection{Proposition 3}\label{Proof_Prop3}
To obtain the sample distribution of $v$ we 
adapt Theorem~3.3.1 
in \citet[p.~75]{Anderson03} to the case of AR(1) processes.

\vspace{0.3cm}

\noindent Consider a $(T-1)\times(T-1)$ orthogonal matrix $D=\left(d_{ht}\right)$ with first row $\mathbf{x}_1'/\sqrt{a_{11}}$
and let $s_h=\sum_{t=1}^{T-1}d_{ht}x_{2t}$, $h=1,\dots,(T-1)$, $t=1,\dots,(T-1)$. We have
$$b=\frac{\sum_{t=1}^{T-1}x_{1t}x_{2t}}{\sum_{t=1}^{T-1}x_{1t}^2}=\frac{\sum_{t=1}^{T-1}d_{1t}x_{2t}}{\sqrt{a_{11}}}=\frac{s_1}{\sqrt{a_{11}}}\ \ .$$
Then, from Lemma 3.3.1 in \citet[p.~76]{Anderson03}, we have
$$v=\sum_{t=1}^{T-1}x_{2t}^2-b^2\sum_{t=1}^{T-1}x_{1t}^2=\sum_{t=1}^{T-1}s_t^2-s_1^2=\sum_{t=2}^{T-1}s_t^2\ \ .$$
Thus, $v$ approximates the sum of $T-2$ Normal variables with variance $1/(1-\phi_2^2)$. 
Now, let $z_t$ be the variable obtained by standardizing $x_{2t}$.
We have
\begin{equation}\label{Gamma}
v=\sum_{t=2}^{T-1}s_t^2\approx \sum_{t=2}^{T-1}\frac{z_t^2}{1-\phi_2^2}\ \ .
\end{equation}
The right side of \eqref{Gamma} is a Gamma distribution with shape parameter $\frac{T-2}{2}$ and rate parameter $\frac{2}{1-\phi_2^2}$.
 \hfill $\blacksquare$

\subsection{Theorem 1}\label{Proof_Theorem}

Because of Proposition 2, $\sqrt{a_{11}}b$ is approximately
$N\left(0,\frac{1-\phi_1^2\phi_2^2}{(1-\phi_2^2)(1-\phi_1\phi_2)^2}\right)$. 
Let $\delta^2=\frac{1-\phi_1^2\phi_2^2}{(1-\phi_2^2)(1-\phi_1\phi_2)^2}$, $\theta^2=\frac{1}{1-\phi_2^2}$ and $t=\frac{\sqrt{a_{11}}\ b}{\sqrt{v/(T-2)}}$. 
In the reminder of the proof, we consider the distributions of $b$ and $v$
in Propositions 2 and 3 as 
exact, not approximate.
Thus, we have the 
densities
\begin{align}
g\left(\sqrt{a_{11}}b\right) & = \frac{1}{\delta\sqrt{2\pi}}\text{exp}\left(-\frac{a_{11}b^2}{2\delta^2}\right) \ ,
\label{gb} \\
h(v) & = \frac{1}{(2\theta^2)^{\frac{T-2}{2}}\Gamma\left(\frac{T-2}{2}\right)}v^{\frac{T-2}{2}-1}\text{exp}\left(-\frac{v}{2\theta^2}\right)\ .
\label{hv}
\end{align}
We focus on
\begin{align}
\begin{split}
f(t) &  =  \int\sqrt{\frac{v}{T-2}}g\left(\sqrt{\frac{v}{T-2}}t\right)h(v)dv
\nonumber \\
& =\int_{0}^{\infty}\sqrt{\frac{v}{T-2}}\frac{1}{\delta\sqrt{2\pi}}\text{exp}\left(-\frac{vt^2}{(T-2)2\delta^2}\right)\frac{v^{\frac{T-2}{2}-1}\text{exp}\left(-\frac{v}{2\theta^2}\right)}{(2\theta^2)^{\frac{T-2}{2}}\Gamma\left(\frac{T-2}{2}\right)}dv
\\
& =\frac{1}{\sqrt{2\pi(T-2)}\delta(2\theta^2)^{\frac{T-2}{2}}\Gamma\left(\frac{T-2}{2}\right)}\int_0^{\infty}v^{\frac{1}{2}}v^{\frac{T-2}{2}-1}\text{exp}\left(-\frac{vt^2}{(T-2)2\delta^2}\right)\text{exp}\left(-\frac{v}{2\theta^2}\right)dv
\\
& =\frac{1}{\sqrt{2\pi(T-2)}\delta(2\theta^2)^{\frac{T-2}{2}}\Gamma\left(\frac{T-2}{2}\right)}\int_0^{\infty}v^{\frac{T-3}{2}}\text{exp}\left(-\left(\frac{1}{\theta^2}\frac{t^2}{(T-2)\delta^2}\right)\frac{v}{2}\right)dv \ .
\end{split}
\end{align}
Now define $\Upsilon=\frac{1}{\sqrt{2\pi(T-2)}\delta(2\theta^2)^{\frac{T-2}{2}}\Gamma\left(\frac{T-2}{2}\right)}$ and $x=\left(\frac{1}{\theta^2}+\frac{t^2}{(T-2)\delta^2}\right)\frac{v}{2}$. Then
\begin{align}
\begin{split}
f(t) &  = \Upsilon \int_0^{\infty}\left[2x\left(\frac{1}{\theta^2}+\frac{t^2}{(T-2)\delta^2}\right)^{-1}\right]^{\frac{T-3}{2}}\text{exp}\left(-x\right)dx
\nonumber \\
& = \Upsilon \ 2^{\frac{T-1}{2}}\left(\frac{1}{\theta^2}+\frac{t^2}{(T-2)\delta^2}\right)^{-\frac{T-1}{2}}
\int_0^{\infty}x^{\frac{T-1}{2}-1}\text{exp}\left(-x\right)dx\ .
\end{split}
\end{align}
The integral on the right hand side can be represented by using the gamma function
$$\Gamma(\alpha)=\int_0^{\infty}x^{\alpha-1}\text{exp}\left(-x\right)dx\ .$$ 
Thus we obtain
\begin{align}
\begin{split}
f(t) & = \Upsilon\ 2^{\frac{T-1}{2}}\left[\frac{(T-2)\delta^2+t^2\theta^2}{\theta^2(T-2)\delta^2}\right]^{-\frac{T-1}{2}}\Gamma\left(\frac{T-1}{2}\right)
\nonumber \\
& = \frac{\Gamma\left(\frac{T-1}{2}\right)2^{\frac{T-1}{2}}}{\sqrt{2\pi(T-2)}\delta(2\theta^2)^{\frac{T-2}{2}}\Gamma\left(\frac{T-2}{2}\right)}\left[\frac{(T-2)\delta^2+t^2\theta^2}{\theta^2(T-2)\delta^2}\right]^{-\frac{T-1}{2}}
\\
& = \frac{\Gamma\left(\frac{T-1}{2}\right)\theta}{\sqrt{\pi(T-2)}\delta \Gamma\left(\frac{T-2}{2}\right)}\left[\frac{(T-2)\delta^2+t^2\theta^2}{(T-2)\delta^2}\right]^{-\frac{T-1}{2}}.
\end{split}
\end{align}
Substituting $\delta^2$ with $\frac{1-\phi_1^2\phi_2^2}{(1-\phi_2^2)(1-\phi_1\phi_2)^2}$ and $\theta^2$ with $\frac{1}{1-\phi_2^2}$, 
we obtain the density 
\begin{align}
\begin{split}
f(t) & = \frac{\Gamma\left(\frac{T-1}{2}\right)(1-\phi_1\phi_2)\sqrt{(1-\phi_2^2})}{\Gamma\left(\frac{T-2}{2}\right)\sqrt{\pi(T-2)(1-\phi_1^2\phi_2^2)(1-\phi_2^2)}}\left[1+\frac{t^2(1-\phi_1\phi_2)^2(1-\phi_2^2)}{(T-2)(1-\phi_1^2\phi_2^2)(1-\phi_2^2)}\right]^{-\frac{T-1}{2}}
\nonumber\\
& = \frac{\Gamma\left(\frac{T-1}{2}\right)(1-\phi_1\phi_2)}{\Gamma\left(\frac{T-2}{2}\right)\sqrt{\pi(T-2)(1-\phi_1^2\phi_2^2)}}\left[1+\frac{t^2(1-\phi_1\phi_2)^2}{(T-2)(1-\phi_1^2\phi_2^2)}\right]^{-\frac{T-1}{2}}
\ .
\end{split}
\end{align}
The density of $w=\widehat{c}_{12}^x\left[1-(\widehat{c}_{12}^x)^2\right]^{-\frac{1}{2}}$ is thus
\begin{align}
\begin{split}
f(w)& = \frac{\Gamma\left(\frac{T-1}{2}\right)(1-\phi_1\phi_2)}{\Gamma\left(\frac{T-2}{2}\right)\sqrt{\pi(1-\phi_1^2\phi_2^2)}}\left[1+\frac{w^2(1-\phi_1\phi_2)^2}{(1-\phi_1^2\phi_2^2)}\right]^{-\frac{T-1}{2}} \ . \nonumber
\end{split}
\end{align}
Next,
define $\kappa(\widehat{c}_{12}^x)=w=\widehat{c}_{12}^x\left[1-(\widehat{c}_{12}^x)^2\right]^{-\frac{1}{2}}$, from which $\kappa'(\widehat{c}_{12}^x)=\left[1-(\widehat{c}_{12}^x)^2\right]^{-\frac{3}{2}}$,
$\phi_{12}=\phi_1\phi_2$ and $\Theta=\left[\Gamma\left(\frac{T-1}{2}\right)(1-\phi_{12})\right]/\left[\Gamma\left(\frac{T-2}{2}\right)\sqrt{\pi(1-\phi_{12}^2)}\right]$.
We can use these quantities to write
\begin{align}
\begin{split}
\mathcal{D}_{\widehat{c}_{12}^x}  & = f_w(\kappa(\widehat{c}_{12}^x))\kappa'(\widehat{c}_{12}^x) = \Theta\left[1+\left(\widehat{c}_{12}^x(1-(\widehat{c}_{12}^x)^2)^{-\frac{1}{2}}\right)^2\frac{(1-\phi_{12})^2}{(1-\phi_{12}^2)}\right]^{-\frac{T-1}{2}}\left[1-(\widehat{c}_{12}^x)^2\right]^{-\frac{3}{2}}
\nonumber \\
& = \Theta\left[1+\frac{(\widehat{c}_{12}^x)^2(1-\phi_{12})^2}{(1-(\widehat{c}_{12}^x)^2)(1-\phi_{12}^2)}\right]^{-\frac{T-1}{2}}\left[1-(\widehat{c}_{12}^x)^2\right]^{-\frac{3}{2}}
\\
& = \Theta\left[\frac{(1-(\widehat{c}_{12}^x)^2)(1-\phi_{12}^2)+(\widehat{c}_{12}^x)^2(1-\phi_{12})^2}{(1-(\widehat{c}_{12}^x)^2)(1-\phi_{12}^2)}\right]^{-\frac{T-1}{2}}\left[1-(\widehat{c}_{12}^x)^2\right]^{-\frac{3}{2}}
\\
& = \Theta\left[\frac{1-\phi_{12}^2+2(\widehat{c}_{12}^x)^2\phi_{12}(\phi_{12}-1)}{(1-(\widehat{c}_{12}^x)^2)(1-\phi_{12}^2)}\right]^{-\frac{T-1}{2}}\left[1-(\widehat{c}_{12}^x)^2\right]^{-\frac{3}{2}} 
\\
& =\Theta\left[1-(\widehat{c}_{12}^x)^2\right]^{\frac{T-4}{2}}\left[\frac{(1-\phi_{12}^2)}{1-\phi_{12}^2+2(\widehat{c}_{12}^x)^2\phi_{12}(\phi_{12}-1)}\right]^{\frac{T-1}{2}} \ .
\end{split}
\end{align}
%
%
Thus, the 
(finite) sample density of $\widehat{c}_{12}^x$, 
taking the densities in \eqref{gb} and \eqref{hv} 
as exact, is
$$\mathcal{D}_{\widehat{c}_{12}^x} = \frac{\Gamma\left(\frac{T-1}{2}\right)(1-\phi_{12})}{\Gamma\left(\frac{T-2}{2}\right)\sqrt{\pi}}\left[1-(\widehat{c}_{12}^x)^2\right]^{\frac{T-4}{2}}\left(1-\phi_{12}^2\right)^{\frac{T-2}{2}}\left[\frac{1}{1-\phi_{12}^2+2(\widehat{c}_{12}^x)^2\phi_{12}(\phi_{12}-1)}\right]^{\frac{T-1}{2}} \ .$$
 \hfill $\blacksquare$

\subsection{Proof of Proposition 4}\label{Proof_Prop4}

In this Section we show that the probability limits of the OLS coefficients are biased as a consequence of serial dependence. To simplify notation, we report results in matrix form, where the subscript $-j$ denotes the corresponding matrix or vector of $j$ period lagged values. Then, we consider the following equations\begin{align}
 & \mathbf{y}=\phi_y\mathbf{y}_{-1}+\mathbf{X}_{-h}'\pmb{\alpha}^*+\pmb{\varepsilon}, \\
 & \mathbf{X}=\phi\mathbf{X}_{-1}+\mathbf{U},\\
 & \pmb{\varepsilon}=\phi_{\varepsilon}\pmb{\varepsilon}_{-1}+\pmb{\omega}
\end{align}
\noindent We provide the convergence in probability of the OLS estimates of $\phi_y$ and $\pmb{\alpha}^*$, namely $\widehat{\phi}_y$ and $\widehat{\pmb{\alpha}}$, when $\phi_y=0$, i.e. when we incorrectly include $y_{t-1}$ in our model. This results are based on the contibution of \cite{Achen2000}. Note that the following results hold for any lag period $h$ in which $\mathbf{X}_{-h}$ loads on $\mathbf{y}$.  

\vspace{0.3cm}

\noindent As a consequence of Assumptions 2 and 3 in the text of the main paper:
$$\mathbf{X}_{-h}\pmb{\varepsilon}/T\ \overset{p}{\to}\ 0; \ \ \ \ \  \phi\mathbf{X}_{-1}\mathbf{U}'/T\ \overset{p}{\to}\ 0; \ \ \ \ \  \phi_{\varepsilon}\pmb{\varepsilon}_{-1}'\pmb{\omega}/T\ \overset{p}{\to}\ 0; \ \ \ \ \  \mathbf{U}_{-h}\pmb{\omega}/T\ \overset{p}{\to}\ 0.$$
\noindent Let:
$$\mathbf{X}_{-h}\mathbf{X}_{-h}'/T\ \overset{p}{\to}\ \mathbf{C}_x; \ \ \ \ \   \pmb{\varepsilon}'\pmb{\varepsilon}/T\ \overset{p}{\to}\ c_{\varepsilon}.$$
\noindent Preliminary results:
$$\mathbf{X}_{-h}\mathbf{y}/T\ \overset{p}{\to}\ \mathbf{C}_x\pmb{\alpha}^*; \ \ \ \ \  \mathbf{X}_{-h}\mathbf{y}_{-1}/T\ \overset{p}{\to}\ \phi\mathbf{C}_x\pmb{\alpha}^*;$$  
$$\pmb{\varepsilon}'\mathbf{y}_{-1}\pmb{\omega}/T\ \overset{p}{\to}\ \phi_{\varepsilon}c_{\varepsilon}; \ \ \ \ \  \mathbf{y}_{-1}'\mathbf{y}/T\ \overset{p}{\to}\ \phi{\pmb{\alpha}^*}'\mathbf{C}_x\pmb{\alpha}^*+\phi_{\varepsilon}c_{\varepsilon}.$$
\noindent Moreover, $var(y)=\sigma_y^2={\pmb{\alpha}^*}'\mathbf{C}_x\pmb{\alpha}^*+c_{\varepsilon}$.
$$\begin{bmatrix}
\widehat{\phi}_y\\
\widehat{\pmb{\alpha}}
\end{bmatrix}
= 
\begin{bmatrix}
\mathbf{y}_{-1}'\mathbf{y}_{-1}/T & \mathbf{y}_{-1}'\mathbf{X}_{-h}'/T \\
\mathbf{X}_{-h}\mathbf{y}_{-1}/T & \mathbf{X}_{-h}\mathbf{X}_{-h}'/T 
\end{bmatrix}^{-1}
\begin{bmatrix}
\mathbf{y}_{-1}'\mathbf{y}/T\\
\mathbf{X}_{-h}\mathbf{y}/T
\end{bmatrix}.$$
\noindent Applying the results above gives:
$$\begin{bmatrix}
\widehat{\phi}_y\\
\widehat{\pmb{\alpha}}
\end{bmatrix}
\ \overset{p}{\to}\  
\begin{bmatrix}
\sigma_y^2 & \phi{\pmb{\alpha}^*}'\mathbf{C}_x \\
\phi\mathbf{C}_x\pmb{\alpha}^* & \mathbf{C}_x 
\end{bmatrix}^{-1}
\begin{bmatrix}
\phi{\pmb{\alpha}^*}'\mathbf{C}_x\pmb{\alpha}^*+\phi_{\varepsilon}c_{\varepsilon}\\
\mathbf{C}_x\pmb{\alpha}^*
\end{bmatrix}.$$
\indent Setting $s=\sigma_y^2-\phi^2{\pmb{\alpha}^*}'\mathbf{C}_x\pmb{\alpha}^*$, we have
$$\begin{bmatrix}
\widehat{\phi}_y\\
\widehat{\pmb{\alpha}}
\end{bmatrix}
\ \overset{p}{\to}\ 
\cfrac{1}{s} 
\begin{bmatrix}
1 & -\phi{\pmb{\alpha}^*}' \\
-\phi\pmb{\alpha}^* & s\mathbf{C}_x^{-1}+\phi^2\pmb{\alpha}^*{\pmb{\alpha}^*}'
\end{bmatrix}
\begin{bmatrix}
\phi{\pmb{\alpha}^*}'\mathbf{C}_x\pmb{\alpha}^*+\phi_{\varepsilon}c_{\varepsilon}\\
\mathbf{C}_x\pmb{\alpha}^*
\end{bmatrix}.$$
After some algebra we get:
$$\widehat{\phi}_{y}\ \overset{p}{\to}\ \cfrac{\phi_{\varepsilon}c_{\varepsilon}}{s}\ , \ \ \ \ \ \ \ \ \ \widehat{\pmb{\alpha}}\ \overset{p}{\to}\ \pmb{\alpha}^*\left(1-\cfrac{\phi\phi_{\varepsilon}c_{\varepsilon}}{s}\right).$$
Considering that $c_{\varepsilon}=\sigma_y^2-{\pmb{\alpha}^*}'\mathbf{C}_x\pmb{\alpha}^*$, and ${\pmb{\alpha}^*}'\mathbf{C}_x\pmb{\alpha}^*=R^2\sigma_y^2$, then
$$\widehat{\phi}_{y}\ \overset{p}{\to}\ \cfrac{\phi_{\varepsilon}(1-R^2)}{1-\phi^2R^2}\ ,\ \ \ \ \ \ \ \ \ \widehat{\pmb{\alpha}}\ \overset{p}{\to}\ \pmb{\alpha}^*\left(1-\cfrac{\phi\phi_{\varepsilon}(1-R^2)}{1-\phi^2R^2}\right).$$
\noindent If we replace $\mathbf{X}_{-h}$ with $\mathbf{U}_{-h}$ in the estimated model, then we have the following preliminary results:
$$\mathbf{U}_{-h}\mathbf{U}_{-h}'/T\ \overset{p}{\to}\ \mathbf{C}_u; \ \ \ \ \  \mathbf{U}_{-h}\mathbf{y}/T\ \overset{p}{\to}\ \mathbf{C}_u\pmb{\alpha}^*; \ \ \ \ \ \mathbf{U}_{-h}\mathbf{y}_{-1}/T\ \overset{p}{\to}\ 0.$$
Thus, 
$$\begin{bmatrix}
\widehat{\phi}_y\\
\widehat{\pmb{\alpha}}
\end{bmatrix}
= 
\begin{bmatrix}
\mathbf{y}_{-1}'\mathbf{y}_{-1}/T & \mathbf{y}_{-1}'\mathbf{U}_{-h}'/T \\
\mathbf{U}_{-h}\mathbf{y}_{-1}/T & \mathbf{U}_{-h}\mathbf{U}_{-h}'/T 
\end{bmatrix}^{-1}
\begin{bmatrix}
\mathbf{y}_{-1}'\mathbf{y}/T\\
\mathbf{U}_{-h}\mathbf{y}/T
\end{bmatrix},$$
\noindent and applying the results above gives:
$$\begin{bmatrix}
\widehat{\phi}_y\\
\widehat{\pmb{\alpha}}
\end{bmatrix}
\ \overset{p}{\to}\  
\begin{bmatrix}
\sigma_y^2 & 0 \\
0 & \mathbf{C}_u 
\end{bmatrix}^{-1}
\begin{bmatrix}
\phi{\pmb{\alpha}^*}'\mathbf{C}_x\pmb{\alpha}^*+\phi_{\varepsilon}c_{\varepsilon}\\
\mathbf{C}_u\pmb{\alpha}^* 
\end{bmatrix},$$
from which
$$\begin{bmatrix}
\widehat{\phi}_y\\
\widehat{\pmb{\alpha}}
\end{bmatrix}
\ \overset{p}{\to}\ 
\cfrac{1}{\sigma_y^2} 
\begin{bmatrix}
1 & 0 \\
0 & \cfrac{\sigma_y^2}{\mathbf{C}_u}
\end{bmatrix}
\begin{bmatrix}
\phi{\pmb{\alpha}^*}'\mathbf{C}_x\pmb{\alpha}^*+\phi_{\varepsilon}c_{\varepsilon}\\
\mathbf{C}_u\pmb{\alpha}^* 
\end{bmatrix}.$$
After some algebra we get:
$$\widehat{\phi}_{y}\ \overset{p}{\to}\ \phi R^2+\phi_{\varepsilon}(1-R^2)\ ,\ \ \ \ \ \ \ \ \ \widehat{\pmb{\alpha}}\ \overset{p}{\to}\ \pmb{\alpha}^*.$$ 

\hfill $\blacksquare$

\section{Distribution of $b$}\label{Dist_b}

\begin{figure}[H]
\begin{center}
\begin{minipage}[c]{0.5\textwidth}
\graphicspath{{images/}}
\includegraphics[scale=0.35]{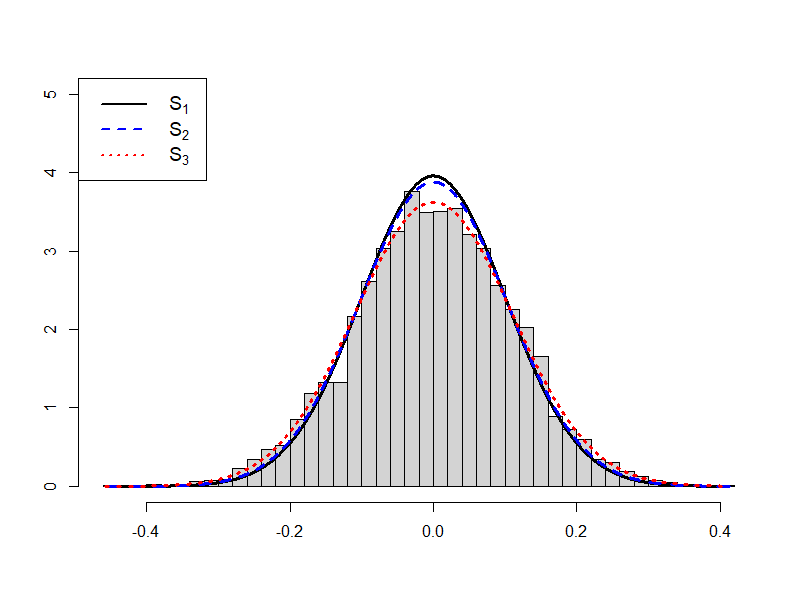}
\subcaption{$\phi=0.3,\ T=100$}
\end{minipage}\hfill{}%
\begin{minipage}[c]{0.5\textwidth}%
\graphicspath{{images/}}
\includegraphics[scale=0.35]{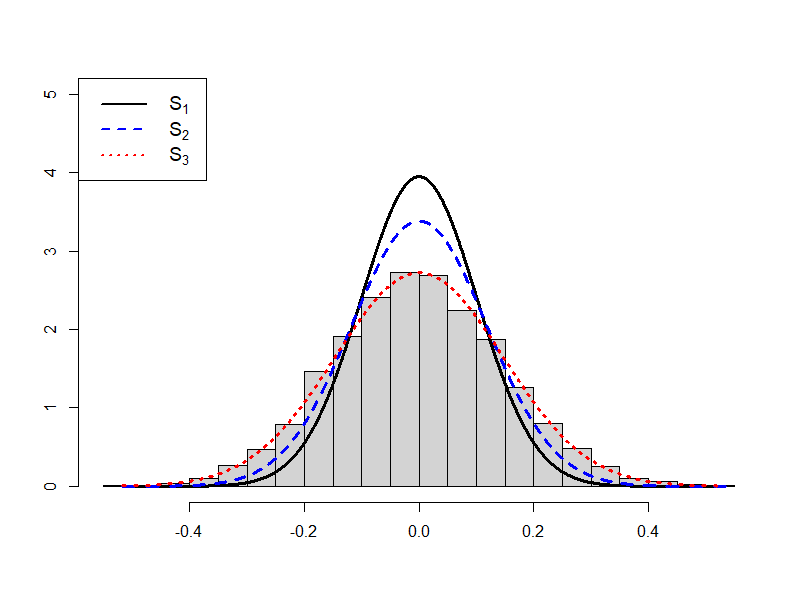}
\subcaption{$\phi=0.6,\ T=100$}
\end{minipage}
\begin{minipage}[c]{0.5\textwidth}%
\graphicspath{{images/}}
\includegraphics[scale=0.35]{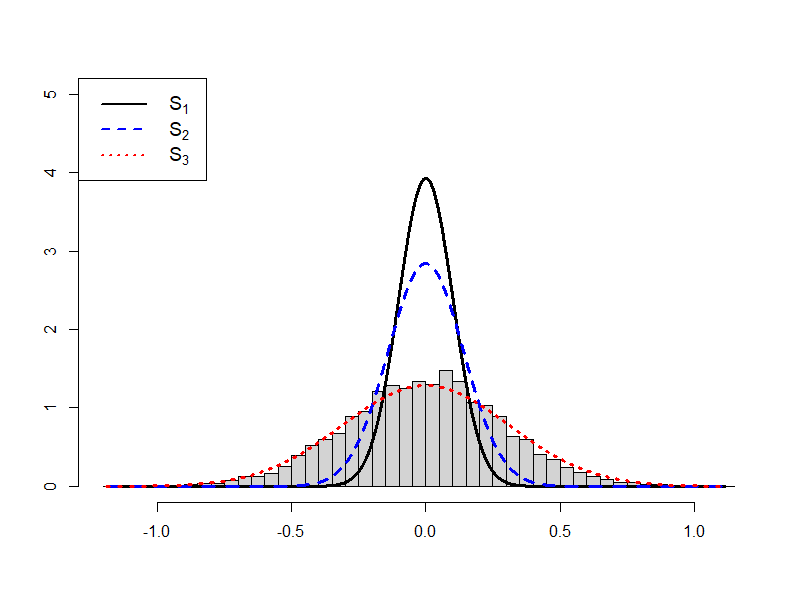}
\subcaption{$\phi=0.9,\ T=100$}
\end{minipage}\hfill{}%
\begin{minipage}[c]{0.5\textwidth}%
\graphicspath{{images/}}
\includegraphics[scale=0.35]{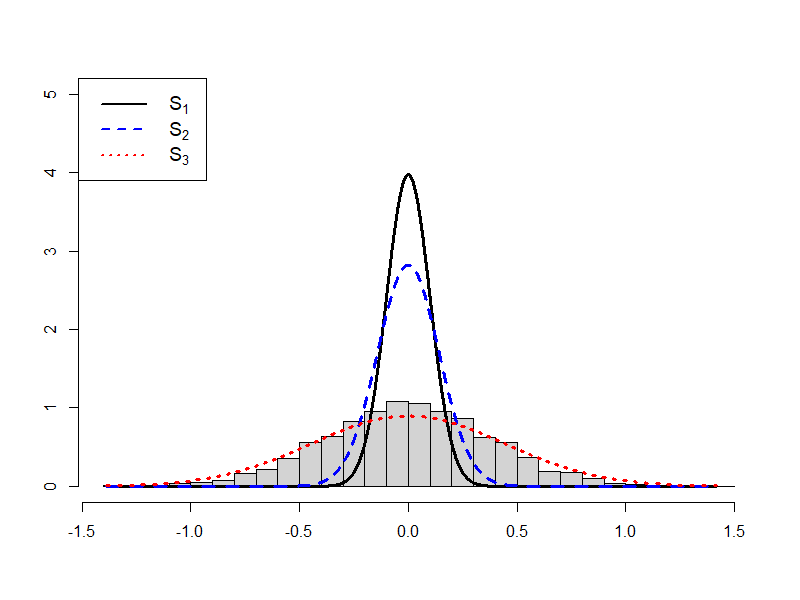}
\subcaption{$\phi=0.95,\ T=100$}
\end{minipage}
\end{center}
\caption{\footnotesize Density of $b$ between uncorrelated AR(1) Gaussian processes. Solid line indicates the approximated density obtained by using the classical OLS estimator, dashed line indicates the approximated density obtained by using the NW estimator, and, finally, dotted line shows the theoretical approximated density obtained in Proposition 2.}\label{b_Dist}
\end{figure}

Consider two orthogonal Gaussian AR(1) processes generated according to the model ${x}_{it}=\phi_i{x}_{it-1}+{u}_{it}$, where $u_{it}\sim N(0,1)$, $i=1,2,\ t=1,\dots,100$ and $\phi_1=\phi_2=\phi$. In this simulation exercise we run the model 
\begin{equation}\label{LinMod_x2x1APP}
x_{2t}=\beta x_{1t}+e_t,
\end{equation}
where $e_t\sim N(0, \sigma_e^2)$, and study the distribution of the OLS estimator $b$ of $\beta$ in the following four cases in terms of degrees of serial dependence: $\phi=0.3,0.6,0.9,0.95$. Figure \ref{b_Dist} reports the density of $b$ across the $\phi$ values obtained on 5000 Monte Carlo replications. We compare this density with that of three zero-mean Gaussian variables where the variances are respectively:
\begin{itemize}
\item $S_1^2=\frac{\widehat{\sigma}^2_{\widehat{e}}}{\sum_{t=1}^T(x_{1t}-\overline{x}_1)^2}$, where $\widehat{\sigma}^2_{\widehat{e}}$ is the sample variance of the estimated residual $\widehat{e}_t=x_{2t}-bx_{1t}$. This is the OLS estimator for the variance of $\beta$.
\item $S_2^2=\frac{1}{T}\frac{\frac{1}{T-2}\sum_{t=1}^T(x_{1t}-\overline{x}_1)^2\widehat{e}_t^2}{\left[\frac{1}{T}\sum_{t=1}^T(x_{1t}-\overline{x}_1)^2\right]^2}\widehat{f}_t$, is the Newey-West (NW) HAC estimator \citep{NeweyWest87,StockWatson2008}, where $\widehat{f}_t=\left(1+2\sum_{j=1}^{m-1}\left(\frac{m-j}{m}\right)\widehat{\rho}_j\right)$ is the correction factor that adjusts for serially correlated errors and involves estimates of $m - 1$ autocorrelation coefficients $\widehat{\rho}_j$, and $\widehat{\rho}_j=\frac{\sum_{t=j+1}^T\widehat{v}_t\widehat{v}_{t-j}}{\sum_{t=1}^T\widehat{v}_t^2}$, with $\widehat{v}_t=(x_{1t}-\overline{x}_1)\widehat{e}_t$. A rule of thumb for choosing $m$ is $m=[0.75T^{1/3}]$.
\item $S_3^2=\frac{(1-\phi_1^2\phi_2^2)(1-\phi_1^2)}{(T-1)(1-\phi_2^2)(1-\phi_1\phi_2)^2}$, is the theoretical variance of $b$ obtained in Proposition 2.
\end{itemize}

\noindent From Figure \ref{b_Dist} we observe that the variance of $b$ increases with the degree of serial dependence ($\phi$) in a way that is well approximated by the distribution derived in Proposition 2 (see dotted line). On the contrary, OLS (solid line) and NW (dashed line), are highly sub-optimal in the presence of strong serial dependence, underestimating the variability of $b$ as the serial dependence increases.

\section{More General Cases}\label{MoreCases}

We study the density of $\widehat{c}_{12}^x$ in three different cases: non-Gaussian processes; weakly and high cross-correlated processes; and ARMA processes with different order. Note that for the first two cases the variables are AR(1) processes with $T=100$ and autocorrelation coefficient $\phi=0.3,0.6,0.9,0.95$. Since we do not have $\mathcal{D}_{\widehat{c}_{12}^x}$ for these cases, we rely on the densities obtained on 5000 Monte Carlo replications, i.e. $d_s(\widehat{c}_{12}^x)$, to show the effect of serial dependence on $\Pr\left\{|\widehat{c}_{12}^x|\geq\tau\right\}$.

\vspace{1cm}

\noindent{\bf The Impact of non-Gaussianity}\normalsize

\begin{figure}[H]
\graphicspath{{images/}}
\centering
  \captionsetup[subfigure]{oneside,margin={0.5cm,0cm}}
\subfloat[\footnotesize Laplace]{\includegraphics[width=4.5cm]{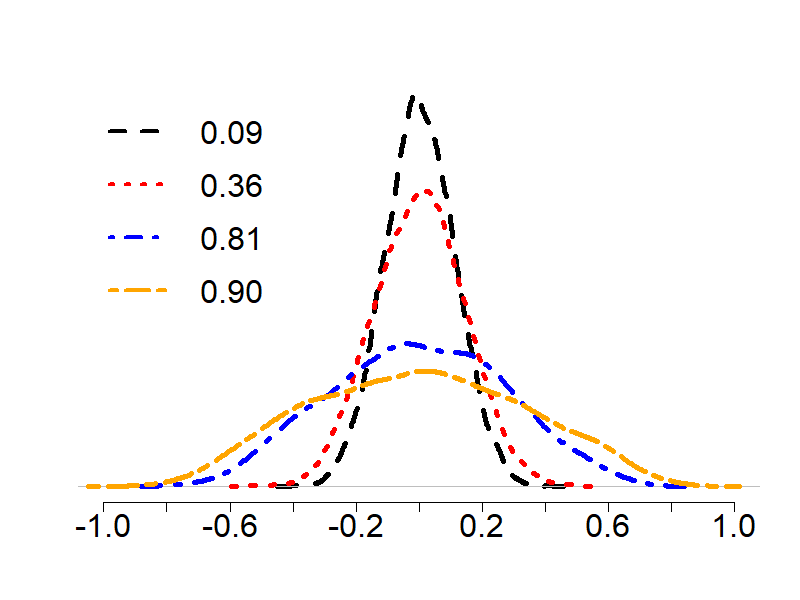}}\hfil
\subfloat[\footnotesize Cauchy]{\includegraphics[width=4.5cm]{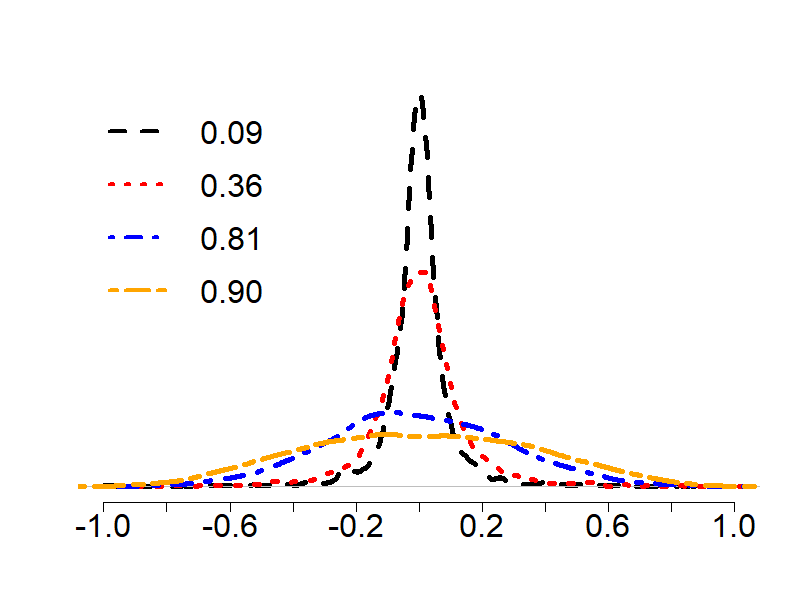}}\hfil 
\subfloat[\footnotesize t-Student]{\includegraphics[width=4.5cm]{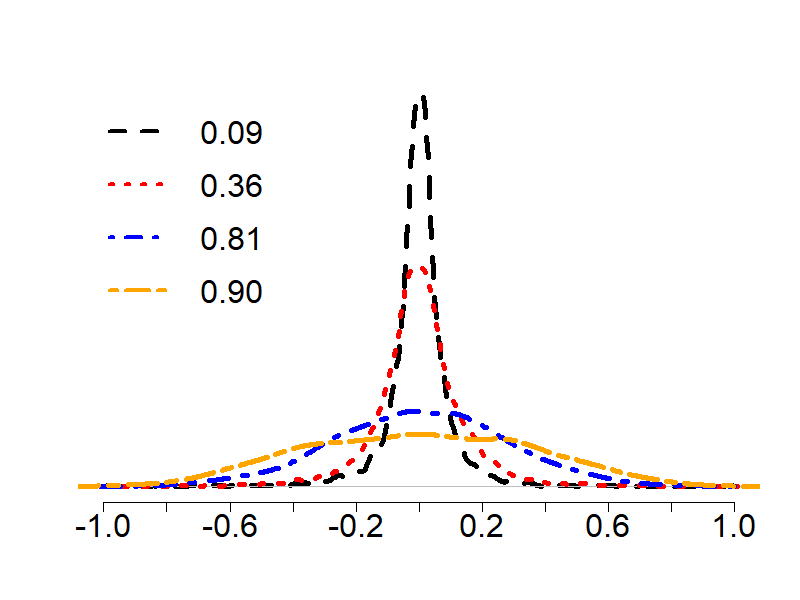}} 
\caption{\footnotesize Simulated density of $\widehat{c}_{12}^x$ in the case of non-Gaussian processes, for $T=100$ and various values of $\phi$.}\label{NoGauss}
\end{figure}

\noindent The theoretical  contribution reported in Section 3 requires the Gaussianity of $u_1$ and $u_2$. With the following simulation experiments we show that the impact of $\phi_{12}$ on the density of $\widehat{c}_{12}^x$ is relevant also when $u_{1t}$ and $u_{2t}$ are non-Gaussian  random variables. To this end, we generate $u_{1t}$ and $u_{2t}$ from the following distributions: Laplace with mean 0 and variance 1 (case (a)); Cauchy with location parameter 0 and scale parameter 1 (case (b)); and from a $t$-student with 1 degree of freedom (case (c)). Figure \ref{NoGauss} reports the results of the simulation experiment. We can state that regardless the distribution of the processes, whenever $Sign(\phi_1)=Sign(\phi_2)$, the probability of large values of $\widehat{c}_{12}^x$ increases with $\phi_{12}$. As a curiosity, this result is more evident for the case of Laplace variables, whereas for Cauchy and $t$-student the effect of $\phi_{12}$ declines. 

\vspace{1cm}

\noindent{\bf The Impact of Population Cross-Correlation}\normalsize

\begin{figure}[H]
\graphicspath{{images/}}
\centering
  \captionsetup[subfigure]{oneside,margin={0.5cm,0cm}}
\subfloat[\footnotesize $c_{12}^u=0.2$]{\includegraphics[width=5cm]{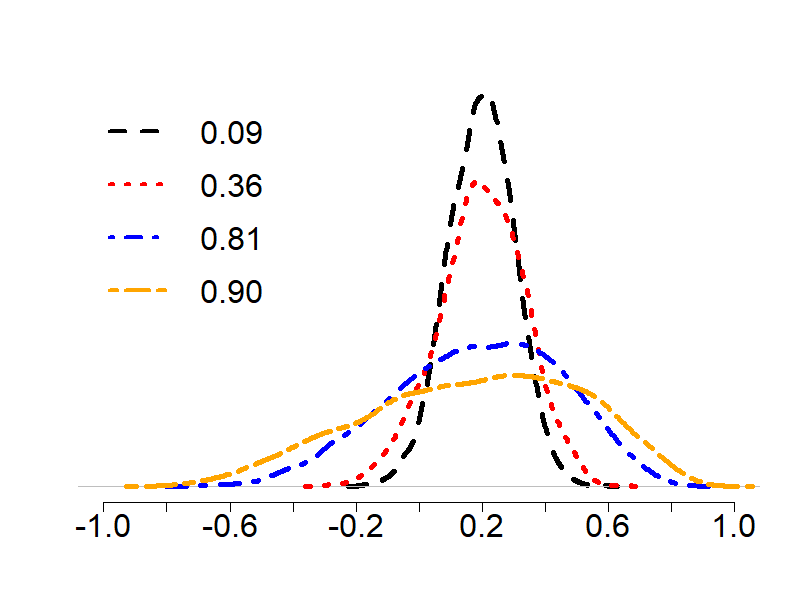}}\hfil
\subfloat[\footnotesize $c_{12}^u=0.8$]{\includegraphics[width=5cm]{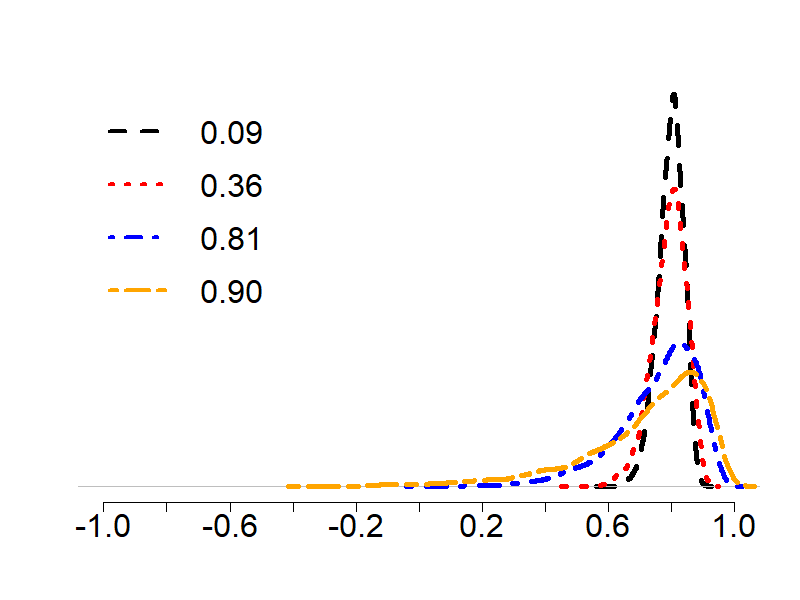}}
\caption{\footnotesize $d_s(\widehat{c}_{12}^x)$ obtained through simulations in the case of $c_{12}^x=0.2$ (a) and $c_{12}^x=0.8$ (b), for $T=100$ and various values of $\phi$.}
\label{fig:rho}
\end{figure}

\noindent Since orthogonality is an unrealistic assumption for most economic applications, here we admit population cross-correlation. In Figure \ref{fig:rho} we report $d_s(\widehat{c}_{12}^x)$ when the processes are weakly cross-correlated with $c_{12}^u=0.2$, and when the processes are multicollinear with $c_{12}^u=0.8$ (usually we refer to multicollinearity when $c_{12}^u\geq0.7$). We observe that the impact of $\phi_{12}$ on $d_s(\widehat{c}_{12}^x)$ depends on the degree of (population) cross-correlation as follows. In the case of weakly correlated processes, an increase in $\phi_{12}$ yields a high probability of observing large sample correlations in absolute value. 
In the case of multicollinear processes, on the other hand, an increase in $\phi_{12}$ leads to a high probability of underestimating the true population cross-correlation.

\vspace{1cm}

\noindent{\bf Density of $\widehat{c}_{12}^x$ in the case of ARMA$(p_i,q_i)$ processes}\normalsize

\begin{figure}[H]
\begin{center}
\graphicspath{{images/}}
\includegraphics[scale=0.35]{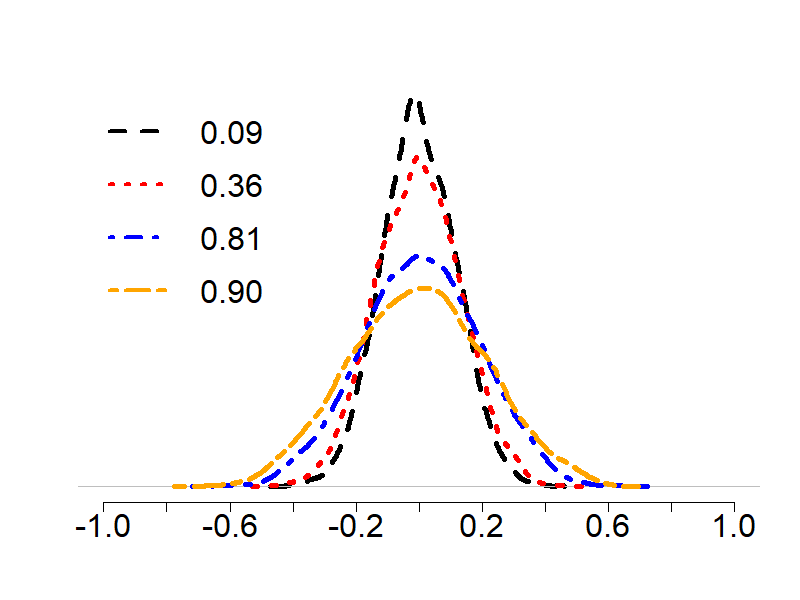}
\end{center}
\caption{\footnotesize Densities of $d_s(\widehat{c}_{12}^x)$ between two uncorrelated ARMA Gaussian processes, for $T=100$ and various values of $\phi$.}\label{fig:ARMA}
\end{figure}

\noindent To show the effect of serial dependence on a more general case, we generate $x_1$ and $x_2$ through the following ARMA processes
\begin{align*}
 & {x}_{1t}={\phi}{x}_{1t-1}+\phi{x}_{1t-2}-\phi{x}_{1t-3}{u}_{1t}+0.5u_{1t-1},\\
 & {x}_{2t}={\phi}{x}_{2t-1}+\phi{x}_{2t-2}+{u}_{2t}+0.7{u}_{2t-1}-0.4u_{3t-2},
\end{align*}
where $t=1,\dots,100$ and $u_i\sim N(0,1)$. In Figure \ref{fig:ARMA} we report the density of $\widehat{c}_{12}^x$ in the case of $T=100$ and $\phi=0.1,0.2,0.3,0.33$. With no loss of generality we can observe that $d_s(\widehat{c}_{12}^x)$ gets larger as $\phi$ increases, that is $\Pr\left\{|\widehat{c}_{12}^x|\geq\tau\right\}$ increases with $|\phi|$.

\section{\setlength{\baselineskip}{0.5\baselineskip}$\widehat{u}$-OLS: Coefficients Estimation and Prediction Accuracy in Low Dimension}\label{uOLS}

We simulate the response through the equation $y_t=\alpha_1^*x_{1t-1}+\alpha_2^*x_{2t-1}+\varepsilon_{t}$, and consider three different scenarios: 

\noindent{\it DGP 1: Equal degrees of serial dependence}, with $x_{it}=\phi_ix_{it-1}+u_{it}$, $\varepsilon_{t}=\phi_{\varepsilon}\varepsilon_{t-1}+\omega_{t}$, for $i=1,2$, $t=1,\dots,T$, where $\phi_1=\phi_2=\phi_{\varepsilon}=0.7$
(as in Example 1). This is the ideal regime in terms of degree of serial dependence, where the working model estimated through $\widehat{u}$-OLS is equivalent to the true model. 

\noindent {\it DGP 2: Different degrees of serial dependence}, with $x_{it}=\phi_ix_{it-1}+u_{it}$, $\varepsilon_{t}=\phi_{\varepsilon}\varepsilon_{t-1}+\omega_{t}$, for $i=1,2$, $t=1,\dots,T$, where $\phi_1=0.75$, $\phi_2=0.6$, and $\phi_{\varepsilon}=0.9$ (as in Example 2). Here the common factor restriction does not hold.

\noindent {\it DGP 3: Different models for predictors and error}. Here we change DGP of covariates. In particular, we consider $x_{1t}=0.6x_{1t-1}+u_{1t}+0.5u_{1t-1}$, $x_{2t}=0.75x_{2t-1}+u_{2t}$, $\varepsilon_{t}=0.6\varepsilon_{t-1}+0.3\varepsilon_{t-2}+\omega_{t}$, $t=1,dots,T$, that is $x_{1t}$, $x_{2t}$ and $\varepsilon_{t}$ are ARMA(1,1), AR(1) and AR(2) processes, respectively.

\vspace{0.3cm}

\noindent For all DGPs  the $i.i.d.$ errors $u_{it}$ and $\omega_t$ are standard Normal random variables for which Assumptions 2 and 3 hold. Here we compare coefficients estimation and forecasting performance of the following methods:
\begin{itemize}
\item \sloppy NW: the Heteroskedasticity and Autocorrelation Consistent (HAC) Newey-West estimator \citep{NeweyWest87}, which accommodates autocorrelation and heteroskedasticity of the error terms. The forecasting equation, in terms of the projection of $y_t$ on the hyperplane spanned by the covariates, is $y_t^{(x)}=\text{Proj}(y_t|y_{t-1},x_{1t-1},x_{2t-1})$.
\item CO: the Cochrane-Orcutt generalized least squares (GLS)
estimator \citep{CocOrc1949}, which adjusts a linear model for serial correlation in the error terms iterating two steps, one to estimate the first order autocorrelation on OLS residuals, and one to transform the variables to eliminate serial dependence in the errors, until a certain criterion is satisfied (e.g., the estimated autocorrelation has converged); transformations are applied from the second observation onward, i.e.~for $t=2,\dots,T$. The forecasting equation is $y_t^{*(x)}=\text{Proj}(y_t^*|x_{1t-1}^*,x_{2t-1}^*)$, where $y_t^*=y_t-\widehat{\phi}_{\varepsilon}^*y_{t-1}$, $x_{it-1}^*=x_{it-1}-\widehat{\phi}_{\varepsilon}^*x_{it-2}$, and $\widehat{\phi}_{\varepsilon}^*$ is the CO estimate of $\phi_{\varepsilon}$. 
\item \sloppy DynReg: the dynamic regression method \citep{Kapetanios2022}, which includes lags of the variables as predictors; if the error is an AR($p$), one adds to the model $p$ lagged values of $y_t$ and $x_{it-1}$, $i=1,\dots,n$. The forecasting equation is $y_t^{(x)}=\text{Proj}(y_t|y_{t-1},x_{1t-1},x_{1t-2},x_{2t-1},x_{2t-2})$ in Scenarios 1 and 2; and $y_t^{(x)}=\text{Proj}(y_t|y_{t-1},y_{t-2},x_{1t-1},x_{1t-2},x_{1t-3},x_{2t-1},x_{2t-2},x_{2t-3})$ in  Scenario 3.
\item $\widehat{u}$-OLS: our proposal, which applies OLS using as predictors $\widehat{u}_{it-1}=x_{it}-\widehat{x}_{it|t-1}$, $i=1,\dots,n$. The forecasting equation is $y_t^{(\widehat{u})}=\text{Proj}(y_t|y_{t-1},\widehat{u}_{1t-1},\widehat{u}_{2t-1})$.
\end{itemize}

\noindent Note that for DGP 1, i.e. under the common factor restriction, also for CO and DynReg the true model and the estimated one coincide. Table \ref{Tab:ExpALL} reports, for each method, the average and standard deviation of the coefficient estimation error $||\widehat{\pmb\alpha}-\pmb{\alpha}^*||_2$ and of the coefficient of determination ($R^2$) over 1000 Monte Carlo replications, considering  $T=100$ (panel (a)) and $T=1000$ (panel (b)). 
Unsurprisingly, NW has the largest coefficient estimation error (it retains OLS estimates and only adjusts standard errors). CO, DynReg and $\widehat{u}$-OLS have smaller and similar coefficient estimation errors. However, in terms of $R^2$, CO is outperformed by DynReg and $\widehat{u}$-OLS, which both include $y_{t-1}$ as predictor in their forecasting equation. We note that, while DynReg and $\widehat{u}$-OLS have similar estimation and prediction performance, DynReg requires the estimation of more parameters. In fact, to express $y_t$ through $n$ covariates and an AR($p$) error, DynReg estimates $p+n+np$ parameters (where $p$ refers to the number of lags of $y_t$). In contrast, $\widehat{u}$-OLS always estimates $p+n$ parameters. This fact highlights the advantage of using our proposal when $n$ is comparable to or larger than $T$. In the next Section we also provide an analysis of the $t$-statistics associated with these methods in the case of spurious regression between uncorrelated autoregressive processes.
\begin{table}[H]
  \centering
\scalebox{0.7}{
\hskip0cm
\begin{tabular}{l@{\hspace{0.8cm}}l@{\hspace{0.8cm}}l@{\hspace{0.8cm}}c@{\hspace{0.8cm}}c@{\hspace{0.8cm}}c@{\hspace{0.8cm}} c c@{\hspace{0.8cm}}c@{\hspace{0.8cm}}c@{\hspace{0.8cm}}c@{\hspace{0.8cm}}c}\hline
DGP & Metric &Stat.& \multicolumn{4}{c}{(a) $T=100$} && \multicolumn{4}{c}{(b) $T=1000$} \\
\cline{4-7} \cline{9-12}
& &  & NW & CO & DynReg & $\widehat{u}$-OLS  && NW & CO & DynReg & $\widehat{u}$-OLS\\\cline{4-7} \cline{9-12}
1 &$||\widehat{\pmb{\alpha}}-\pmb{\alpha}^*||_2$&ave. &0.317
&0.128 &0.129 &0.129 &&0.341 &0.040 &0.040 &0.040\\
&&s.d. &0.132 &0.069 &0.069 &0.068 &&0.046 &0.020 &0.020 &0.020\\
&$R^2$  &ave. &0.747 &0.682 &0.824 &0.817 &&0.747 &0.668 &0.829 &0.829\\
&&s.d. &0.066 &0.056 &0.045 &0.046 &&0.021 &0.017 &0.014 &0.014\\\cline{4-7} \cline{9-12}
2 &$||\widehat{\pmb{\alpha}}-\pmb{\alpha}^*||_2$&ave. &0.351 &0.124 &0.132 &0.134 &&0.379 &0.037 &0.039 &0.040\\
&&s.d. &0.227 &0.066 &0.069 &0.071 &&0.238 &0.020 &0.020 &0.021\\
&$R^2$  &ave. &0.761 &0.704 &0.836 &0.814 &&0.768 &0.695 &0.845 &0.828\\
&&s.d. &0.066 &0.057 &0.059 &0.065 &&0.022 &0.019 &0.046 &0.049\\\cline{4-7} \cline{9-12}
3 &$||\widehat{\pmb{\alpha}}-\pmb{\alpha}^*||_2$&ave. &0.474 &0.126 &0.134 &0.148 &&0.579 &0.038 &0.040 &0.044\\
&&s.d. &0.184 &0.066 &0.070 &0.078 &&0.072 &0.020 &0.021 &0.024\\
&$R^2$  &ave. &0.789 &0.701 &0.888 &0.846 &&0.791 &0.684 &0.900 &0.868\\
&&s.d. &0.072 &0.054 &0.039 &0.051 &&0.034 &0.016 &0.017 &0.020\\
\hline
\end{tabular}
}
\caption{\footnotesize Coefficient estimation error and coefficient of determination ($R^2$) of Newey-West HAC estimator (NW), Cochrane-Orcutt GLS estimator (CO), Dynamic Regression (DynReg) and OLS applied to the estimated ARMA residuals ($\widehat{u}$-OLS) across the three simulation scenarios (DGPs). Panel (a) $T=100$, panel ($b$) $T=1000$. Results are obtained on 1000 Monte Carlo replications.}\label{Tab:ExpALL}
\end{table}

\section{Detecting 
Spurious Regression}\label{SpurReg}
Here we generate data as in Supplement \ref{Dist_b} and compare the t-statistics of the Newey-West-style HAC estimators (NW), Cochrane-Orcutt GLS estimator (CO), Dynamic Regression (DynReg) and the ordinary least squares applied on the estimated AR residuals ($\widehat{u}$-OLS) to evaluate their ability in avoiding spurious regressions. In Table \ref{t_dist} we report the percentage of times that the $t$-statistics are greater than 1.96 in absolute value. Note that according to statistical theory $|t_b|>1.96$ will occur approximately $5\%$ of the time. The main results from this analysis are: ($i$) OLS estimator (Table \ref{t_dist} panel (a)) suffers of spurious regressions for any $\phi>0$, which occurs about $\%50$ when $\phi=0.9$. ($ii$) NW estimator (Table \ref{t_dist} panel (b)) reduces the problem, but for large value of $\phi$ spurious regression occurs frequently. Note that these two results are in line with those in \cite{Granger2001}.  ($iii$) CO, DynReg and $\widehat{u}$-OLS (Table \ref{t_dist} panel (c)-(e)) solve the problem of spurious regression due to serial dependence, by making the variance of $b$ independent of $\phi$. However, $\widehat{u}$-OLS keeps the advantage already mentioned with respect to CO and DynReg, that are a better prediction accuracy and the estimation of less parameters (see Section 5.2 in the text of the main paper).
\begin{table}[H]
\begin{center}
\caption{\footnotesize Percentage of $t$-statistics over 1.96 in absolute value obtained on 1000 Monte Carlo replications.}\label{t_dist}
\scalebox{0.75}{
\begin{tabular}{l@{\hspace{0.8cm}}r@{\hspace{0.8cm}}c@{\hspace{0.8cm}}c@{\hspace{0.8cm}}c@{\hspace{0.8cm}}c@{\hspace{0.8cm}}c}\hline
 & $T$ & $\phi=0.0$ & $\phi=0.3$ & $\phi=0.6$ & $\phi=0.9$ & $\phi=0.95$ \\\cline{2-7}
&50 &5.96 &8.00 &17.16 &47.58 &56.12\\
&100 &5.38 &7.44 &18.00 &50.54 &60.36\\
(a) $|t_{b}^{ols}|>1.96$&250 &6.06 &7.20 &18.26 &51.16 &64.90\\
&1000 &4.94 &7.28 &17.72 &51.82 &66.48\\
&10000 &5.12 &7.08 &19.00 &53.62 &65.76\\\cline{2-7}
&50 &7.32 &9.36 &16.48 &41.82 &50.58\\
&100 &6.96 &7.58 &12.48 &36.36 &47.72\\
(b) $|t_{b}^{nw}|>1.96$&250 &6.96 &6.08 &9.34 &29.00 &43.62\\
&1000 &5.00 &5.24 &7.36 &18.88 &31.72\\
&10000 &5.32 &4.68 &6.08 &9.48 &17.00\\\cline{2-7}
&50 &6.58 &6.76 &7.16 &8.24 &8.42\\
&100 &5.78 &5.60 &5.92 &5.86 &6.28\\
(c) $|t_{b}^{co}|>1.96$&250 &6.16 &4.76 &5.42 &4.52 &5.04\\
&1000 &5.00 &5.06 &4.74 &5.02 &5.20\\
&10000 &5.22 &5.14 &4.80 &4.86 &5.56\\\cline{2-7}
&50 &5.88 &5.94 &6.16 &6.12 &5.52\\
&100 &5.36 &5.16 &5.30 &5.04 &5.52\\
(d) $|t_{b}^{dr}|>1.96$&250 &5.86 &4.62 &5.14 &4.56 &4.86\\
&1000 &4.86 &5.02 &4.84 &4.88 &5.12\\
&10000 &5.22 &5.10 &4.82 &4.86 &5.58\\\cline{2-7}
&50 &6.08 &6.48 &5.58 &6.08 &5.06\\
&100 &5.36 &5.40 &5.34 &4.84 &5.26\\
(e) $|t_{b}^{\widehat{u}-ols}|>1.96$&250 &6.02 &4.66 &5.10 &4.52 &4.70\\
&1000 &4.94 &5.00 &4.78 &4.96 &5.16\\
&10000 &5.12 &5.16 &4.86 &4.84 &5.54\\
\hline
\end{tabular}
}
\end{center}
\end{table}

\noindent As a further analysis, the following Proposition shows that the variability of the limiting distribution of $t_{b}^{ols}$ depends only on the degree of serial dependence of the processes.
\begin{customprop}{G.1}\label{prop:cinque}
Let $S_{ols}^2=\frac{\widehat{\sigma}^2_{\widehat{e}}}{\sum_{t=1}^T(x_{1t}-\overline{x}_1)^2}$, where $\widehat{\sigma}_{\widehat{e}}^2$ is the estimated variance of the residual of model \eqref{LinMod_x2x1APP}. Then 
$$\frac{b}{S_{ols}} \xrightarrow{\text{d}} N\left(0,\frac{1-\phi_1^2\phi_2^2}{(1-\phi_1\phi_2)^2}\right).$$
\end{customprop}

\vspace{0.3cm}

\noindent\textbf{Proof:} From Proposition 2 in main text we know that $b\approx  N\left(0, \frac{(1-\phi_1^2\phi_2^2)(1-\phi_1^2)}{(T-1)(1-\phi_2^2)(1-\phi_1\phi_2)^2}\right)$. Then, considering $S_{ols}^2\approx \frac{1-\phi_1^2}{(T-1)(1-\phi_2^2)}$, we have 
$\frac{b}{S_{ols}}\xrightarrow{\text{d}}  N\left(0,\frac{1-\phi_1^2\phi_2^2}{(1-\phi_1\phi_2)^2}\right).$
 \hfill $\blacksquare$

\vspace{0.3cm}

\noindent Note that the result in Proposition \ref{prop:cinque} has been also derived in \cite{Granger2001}. This result show that the misspecification of $t_b^{ols}$ is only due to the degree of serial dependence. To confirm this, look at the columns of Table \ref{t_dist} and consider that the value of $|t_{b}^{ols}|$ increases with the degree of serial dependence $\phi$, but stay quite constant regardless of the sample size $T$.

\section{More Examples with Different Models}\label{Examples}
For the sake of clarity, models \eqref{Mod_xY} and \eqref{Mod_uY} refer to the true and working model presented in the main text.
\begin{customex}{H.1}\label{ex:tre}
(Equal degrees of serial dependence and different models for the predictors). Consider $x_{it}$ and $x_{jt}$ generated through models $x_{it}=\phi x_{it-1} + \phi x_{it-2} + u_{it}$ and $x_{jt}=\phi x_{jt-1}+\theta u_{jt-1}+u_{jt}$, for $i=1,\dots,q$ and $j=q+1,\dots,n$, where $2|\phi|<1$. Model \eqref{Mod_xY} can be rewritten as 
\begin{align*}
\begin{split}
 y_t & = \sum_{i=1}^q\alpha_i^*(\phi x_{it-2}+\phi x_{it-3}+u_{it-1}) + \sum_{j=1}^{(n-q)}\alpha_j^*(\phi x_{jt-2}+\theta u_{jt-2}+u_{jt-1}) +\phi_{\varepsilon}\varepsilon_{t-1}+\omega_{it} \\
& = \sum_{i=1}^q\alpha_i^*u_{it-1} + \sum_{j=1}^{(n-q)}\alpha_j^*u_{jt-1} + \phi y_{t-1} + \sum_{i=1}^q\phi\alpha_i^*x_{it-3}+\sum_{j=1}^{(n-q)}\alpha_j^*\theta u_{jt-2}+ \omega_t\ .
\end{split}
\end{align*}
Thus, if we have an \enquote{ideal regime} in terms of degree of serial dependence, but different models for the predictors, the working model \eqref{Mod_uY} is not equivalent to the true model \eqref{Mod_xY}. Here, the difference between true and working model is due to the differences between the mechanisms generating $x_{it|t-1}$ and $x_{jt|t-1}$. Again, this makes $y_{t-1}$ not suitable for summarizing the serial dependence of $y_t$.
\end{customex}

\begin{customex}{H.2}\label{ex:quattro}
(Equals degrees of serial dependence and different model for the error). Consider now the case where $\varepsilon_{t}=\phi\varepsilon_{t-1}+\phi\varepsilon_{t-2}+\omega_{t}$ with $2|\phi|<1$. Model \eqref{Mod_xY} becomes
$$ y_t  = \sum_{i=1}^n\alpha_i^*(\phi x_{it-2}+u_{it-1}) +\phi\varepsilon_{t-1}+\phi\varepsilon_{t-2}+\omega_{it} = \sum_{i=1}^n\alpha_i^*u_{1t-1} + \phi y_{t-1} +\phi\varepsilon_{t-2} + \omega_t\ .
$$
Thus, if we have an \enquote{ideal regime} in terms of degree of serial dependence, but a different model for the error, the working model \eqref{Mod_uY} is not equivalent to the true model \eqref{Mod_xY}. Here, the difference between true and working model is due to the differences between the mechanism generating $\varepsilon_{t|t-1}$ and the mechanism generating the predictors. In this case, the residual of the working model would have an autoregressive component.
\end{customex}

\section{Heatmaps CV}\label{HeatmapsCV}
\begin{figure}[H]
\graphicspath{{images/}}
\centering
\subfloat[LASSO, $h$=12]{\includegraphics[width=6.5cm]{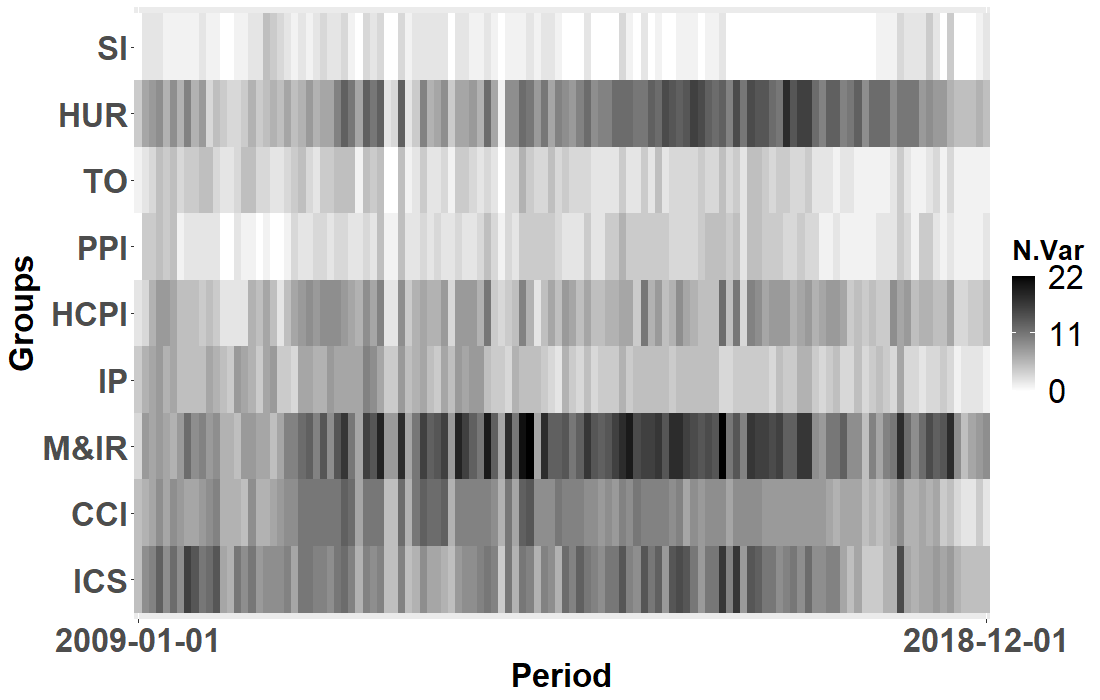}}\hfil
\subfloat[$\widehat{u}$-LASSO, $h$=12]{\includegraphics[width=6.5cm]{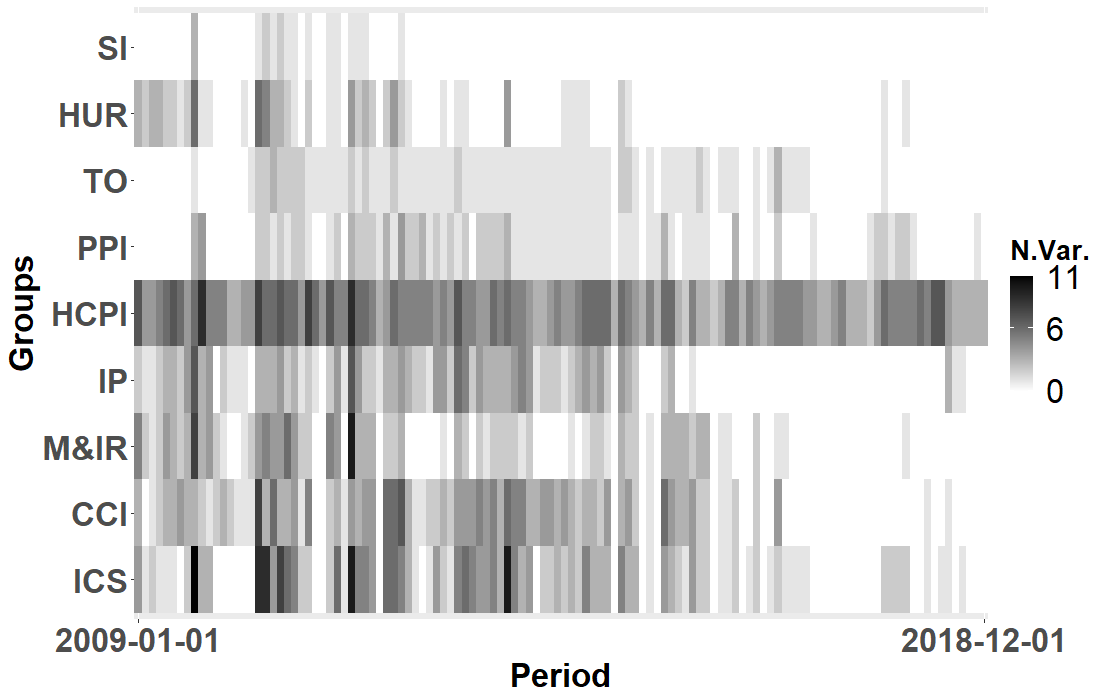}} 

\subfloat[LASSO, $h$=24]{\includegraphics[width=6.5cm]{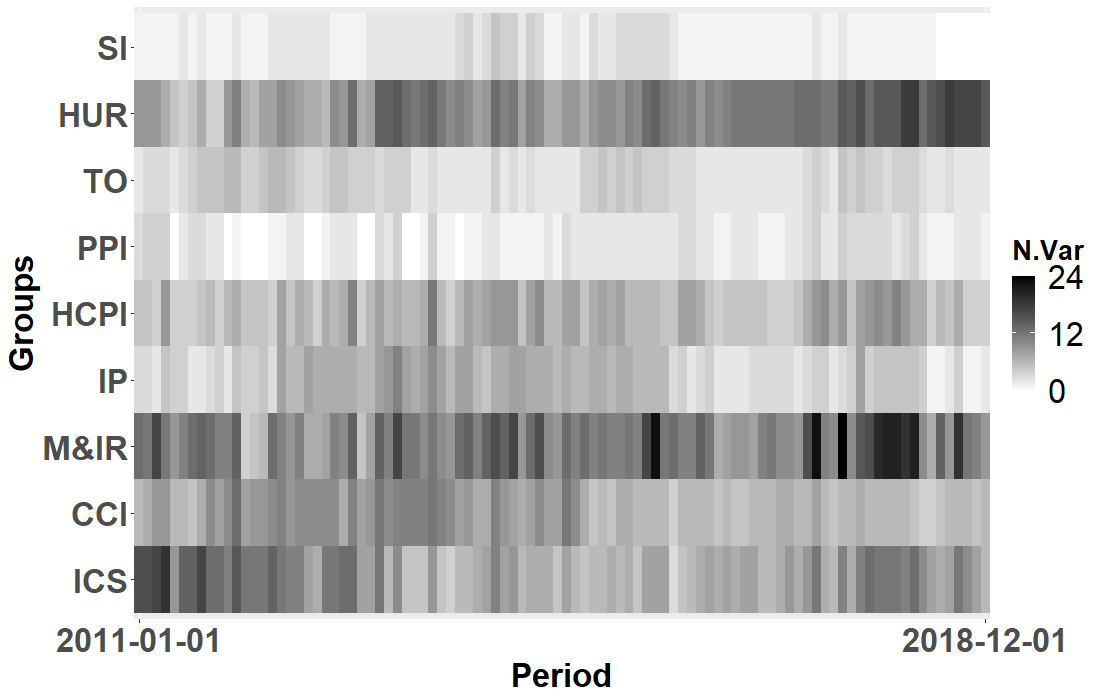}}\hfil
\subfloat[$\widehat{u}$-LASSO, $h$=24]{\includegraphics[width=6.5cm]{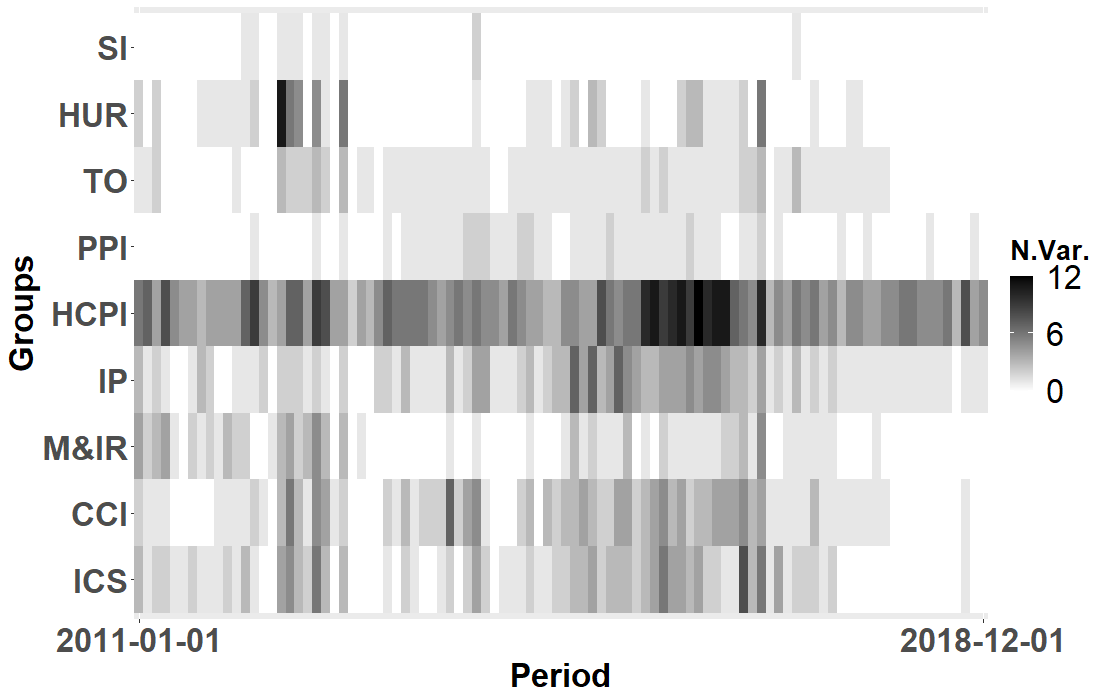}} 
\caption{\footnotesize Heatmaps of the variables selected by LASSO (left column) and $\widehat{u}$-LASSO (right column) when tuning parameter is selected with CV, categorized by groups.}
\end{figure}

\section{List of Time Series 
in the Euro Area Data}\label{Stat_EA}
We report the list of series for the Euro Area dataset adopted in the forecasting exercise (obtained from \cite{Ale2021}). As for the FRED data, the column tcode denotes the data transformation for a given series $x_t$: (1) no transformation; (2) $\Delta x_t$; (3)$\Delta^2 x_t$; (4) $ \text{log}(x_t)$; (5) $\Delta \text{log}(x_t)$; (6) $\Delta^2 \text{log}(x_t)$. (7) $\Delta (x_t/x_{t-} - 1.0)$.

\noindent The acronyms for the sectors refer to:
\begin{itemize}
  \item ICS: Industry \& Construction Survey
  \item CCI: Consumer Confidence Indicators
  \item M\&IR: Money \& Interest Rates
  \item IP: Industrial Production
  \item HCPI: Harm. Consumer Price Index
  \item PPI: Producer Price Index
  \item T\/O: Turnover \& Retail Sale
  \item HUR: Harm. Unemployment rate
  \item SI: Service Svy.
\end{itemize}

\footnotesize
\begin{longtable}{cllcc}
\caption{\footnotesize A sample long table.} \label{tab:long} \\

\hline \multicolumn{1}{c}{\textbf{ID}} &  \multicolumn{1}{c}{\textbf{Description}} & \multicolumn{1}{c}{\textbf{Area}} & \multicolumn{1}{c}{\textbf{Sector}} & \multicolumn{1}{c}{\textbf{Tcode}} \\ \hline
\endfirsthead

\multicolumn{5}{c}%
{{\bfseries \tablename\ \thetable{} -- continued from previous page}} \\
\hline \multicolumn{1}{c}{\textbf{ID}} &  \multicolumn{1}{c}{\textbf{Description}} & \multicolumn{1}{c}{\textbf{Area}} & \multicolumn{1}{c}{\textbf{Sector}} & \multicolumn{1}{c}{\textbf{Tcode}} \\ \hline
\endhead

\hline \multicolumn{5}{r}{{Continued on next page}} \\
\endfoot

\hline \hline
\endlastfoot

    1     & Ind Svy: Employment Expectations  & EA    & ICS   & 1 \\
    2     & Ind Svy: Export Order-Book Levels  & EA    & ICS   & 1 \\
    3     & Ind Svy: Order-Book Levels  & EA    & ICS   & 1 \\
    4     & Ind Svy: Mfg - Selling Price Expectations  & EA    & ICS   & 1 \\
    5     & Ind Svy: Production Expectations  & EA    & ICS   & 1 \\
    6     & Ind Svy: Production Trend  & EA    & ICS   & 1 \\
    7     & Ind Svy: Mfg - Stocks Of Finished Products  & EA    & ICS   & 1 \\
    8     & Constr. Svy: Price Expectations & EA    & ICS   & 1 \\
    9     & Ind Svy: Export Order Book Position  & EA    & ICS   & 1 \\
    10    & Ind Svy: Production Trends In Recent Mth. & EA    & ICS   & 1 \\
    11    & Ind Svy: Selling Prc. Expect. Mth. Ahead  & EA    & ICS   & 1 \\
    12    & Ret. Svy: Employment & EA    & ICS   & 1 \\
    13    & Ret. Svy: Orders Placed With Suppliers & EA    & ICS   & 1 \\
    14    & Constr. Svy: Synthetic Bus. Indicator & FR    & ICS   & 1 \\
    15    & Bus. Svy: Constr. Sector - Capacity Utilisation Rate & FR    & ICS   & 1 \\
    16    & Constr. Svy: Activity Expectations & FR    & ICS   & 1 \\
    17    & Constr. Svy: Price Expectations & FR    & ICS   & 1 \\
    18    & Constr. Svy: Unable To Increase Capacity & FR    & ICS   & 1 \\
    19    & Constr. Svy: Workforce Changes & FR    & ICS   & 1 \\
    20    & Constr. Svy: Workforce Forecast Changes & FR    & ICS   & 1 \\
    21    & Svy: Mfg Output - Order Book \& Demand & FR    & ICS   & 1 \\
    22    & Svy: Mfg Output - Order Book \& Foreign Demand & FR    & ICS   & 1 \\
    23    & Svy: Mfg Output - Personal Outlook & FR    & ICS   & 1 \\
    24    & Svy: Auto Ind - Order Book \& Demand & FR    & ICS   & 1 \\
    25    & Svy: Auto Ind - Personal Outlook & FR    & ICS   & 1 \\
    26    & Svy: Basic \& Fab Met Pdt Ex Mach \& Eq - Personal Outlook & FR    & ICS   & 1 \\
    27    & Svy: Ele \& Elec Eq, Mach Eq - Order Book \& Demand & FR    & ICS   & 1 \\
    28    & Svy: Ele \& Elec Eq, Mach Eq - Order Book \& Foreign Demand & FR    & ICS   & 1 \\
    29    & Svy: Ele \& Elec Eq, Mach Eq - Personal Outlook & FR    & ICS   & 1 \\
    30    & Svy: Mfg Output - Price Outlook & FR    & ICS   & 1 \\
    31    & Svy: Mfg Of Chemicals \& Chemical Pdt - Order Book \& Demand & FR    & ICS   & 1 \\
    32    & Svy: Mfg Of Chemicals \& Chemical Pdt - Personal Outlook & FR    & ICS   & 1 \\
    33    & Svy: Mfg Of Food Pr \& Beverages - Order Book \& Demand & FR    & ICS   & 1 \\
    34    & Svy: Mfg Of Food Pr \& Beverages - Order Book \& Foreign Demand & FR    & ICS   & 1 \\
    35    & Svy: Mfg Of Trsp Eq - Finished Goods Inventories & FR    & ICS   & 1 \\
    36    & Svy: Mfg Of Trsp Eq - Order Book \& Demand & FR    & ICS   & 1 \\
    37    & Svy: Mfg Of Trsp Eq - Order Book \& Foreign Demand & FR    & ICS   & 1 \\
    38    & Svy: Mfg Of Trsp Eq - Personal Outlook & FR    & ICS   & 1 \\
    39    & Svy: Oth Mfg, Mach \& Eq Rpr \& Instal - Ord Book \& Demand & FR    & ICS   & 1 \\
    40    & Svy: Oth Mfg, Mach \& Eq Rpr \& Instal - Ord Book \& Fgn Demand & FR    & ICS   & 1 \\
    41    & Svy: Oth Mfg, Mach \& Eq Rpr \& Instal - Personal Outlook & FR    & ICS   & 1 \\
    42    & Svy: Other Mfg - Order Book \& Demand & FR    & ICS   & 1 \\
    43    & Svy: Rubber, Plastic \& Non Met Pdt - Order Book \& Demand & FR    & ICS   & 1 \\
    44    & Svy: Rubber, Plastic \& Non Met Pdt - Order Book \& Fgn Demand & FR    & ICS   & 1 \\
    45    & Svy: Rubber, Plastic \& Non Met Pdt - Personal Outlook & FR    & ICS   & 1 \\
    46    & Svy: Total Ind - Order Book \& Demand & FR    & ICS   & 1 \\
    47    & Svy: Total Ind - Order Book \& Foreign Demand & FR    & ICS   & 1 \\
    48    & Svy: Total Ind - Personal Outlook & FR    & ICS   & 1 \\
    49    & Svy: Total Ind - Price Outlook & FR    & ICS   & 1 \\
    50    & Svy: Wood \& Paper, Print \& Media - Ord Book \& Fgn Demand & FR    & ICS   & 1 \\
    51    & Trd. \& Ind: Bus Sit & DE    & ICS   & 1 \\
    52    & Trd. \& Ind: Bus Expect In 6Mo & DE    & ICS   & 1 \\
    53    & Trd. \& Ind: Bus Sit & DE    & ICS   & 1 \\
    54    & Trd. \& Ind: Bus Climate & DE    & ICS   & 1 \\
    55    & Cnstr Ind: Bus Climate & DE    & ICS   & 1 \\
    56    & Mfg: Bus Climate & DE    & ICS   & 1 \\
    57    & Mfg: Bus Climate & DE    & ICS   & 1 \\
    58    & Mfg Cons Gds: Bus Climate & DE    & ICS   & 1 \\
    59    & Mfg (Excl Fbt): Bus Climate & DE    & ICS   & 1 \\
    60    & Whsle (Incl Mv): Bus Climate & DE    & ICS   & 1 \\
    61    & Mfg: Bus Sit & DE    & ICS   & 1 \\
    62    & Mfg: Bus Sit & DE    & ICS   & 1 \\
    63    & Mfg (Excl Fbt): Bus Sit & DE    & ICS   & 1 \\
    64    & Mfg (Excl Fbt): Bus Sit & DE    & ICS   & 1 \\
    65    & Cnstr Ind: Bus Expect In 6Mo & DE    & ICS   & 1 \\
    66    & Cnstr Ind: Bus Expect In 6Mo & DE    & ICS   & 1 \\
    67    & Mfg: Bus Expect In 6Mo & DE    & ICS   & 1 \\
    68    & Mfg: Bus Expect In 6Mo & DE    & ICS   & 1 \\
    69    & Mfg Cons Gds: Bus Expect In 6Mo & DE    & ICS   & 1 \\
    70    & Mfg (Excl Fbt): Bus Expect In 6Mo & DE    & ICS   & 1 \\
    71    & Mfg (Excl Fbt): Bus Expect In 6Mo & DE    & ICS   & 1 \\
    72    & Rt (Incl Mv): Bus Expect In 6Mo & DE    & ICS   & 1 \\
    73    & Whsle (Incl Mv): Bus Expect In 6Mo & DE    & ICS   & 1 \\
    74    & Bus. Conf. Indicator & IT    & ICS   & 1 \\
    75    & Order Book Level: Ind & ES    & ICS   & 1 \\
    76    & Order Book Level: Foreign - Ind & ES    & ICS   & 1 \\
    77    & Order Book Level: Investment Goods & ES    & ICS   & 1 \\
    78    & Order Book Level: Int. Goods & ES    & ICS   & 1 \\
    79    & Production Level - Ind & ES    & ICS   & 1 \\
    80    & Cons. Confidence Indicator  & EA    & CCI   & 1 \\
    81    & Cons. Svy: Economic Situation Last 12 Mth. - Emu 11/12 & EA    & CCI   & 1 \\
    82    & Cons. Svy: Possible Savings Opinion & FR    & CCI   & 1 \\
    83    & Cons. Svy: Future Financial Situation & FR    & CCI   & 1 \\
    84    & Svy - Households, Economic Situation Next 12M & FR    & CCI   & 1 \\
    85    & Cons. Confidence Indicator - DE & DE    & CCI   & 1 \\
    86    & Cons. Confidence Index & DE    & CCI   & 5 \\
    87    & Gfk Cons. Climate Svy - Bus. Cycle Expectations & DE    & CCI   & 1 \\
    88    & Cons.S Confidence Index & DE    & CCI   & 5 \\
    89    & Cons. Confidence Climate (Balance) & DE    & CCI   & 1 \\
    90    & Cons. Svy: Economic Climate Index (N.West It) & IT    & CCI   & 5 \\
    91    & Cons. Svy: Economic Climate Index (Southern It) & IT    & CCI   & 5 \\
    92    & Cons. Svy: General Economic Situation (Balance) & IT    & CCI   & 1 \\
    93    & Cons. Svy: Prices In Next 12 Mths. - Lower & IT    & CCI   & 5 \\
    94    & Cons. Svy: Unemployment Expectations (Balance) & IT    & CCI   & 1 \\
    95    & Cons. Svy: Unemployment Expectations - Approx. Same & IT    & CCI   & 5 \\
    96    & Cons. Svy: Unemployment Expectations - Large Increase & IT    & CCI   & 5 \\
    97    & Cons. Svy: Unemployment Expectations - Small Increase & IT    & CCI   & 5 \\
    98    & Cons. Svy: General Economic Situation (Balance) & IT    & CCI   & 1 \\
    99    & Cons. Svy: Household Budget - Deposits To/Withdrawals & ES    & CCI   & 5 \\
    100   & Cons. Svy: Household Economy (Cpy) - Much Worse & FR    & CCI   & 5 \\
    101   & Cons. Svy: Italian Econ.In Next 12 Mths.- Much Worse & FR    & CCI   & 5 \\
    102   & Cons. Svy: Major Purchase Intentions - Balance & FR    & CCI   & 1 \\
    103   & Cons. Svy: Major Purchase Intentions - Much Less & FR    & CCI   & 5 \\
    104   & Cons. Svy: Households Fin Situation - Balance & FR    & CCI   & 1 \\
    105   & Indl. Prod. - Excluding Constr.   & EA    & IP    & 5 \\
    106   & Indl. Prod. - Cap. Goods  & EA    & IP    & 5 \\
    107   & Indl. Prod. - Cons. Non-Durables  & EA    & IP    & 5 \\
    108   & Indl. Prod. - Cons. Durables  & EA    & IP    & 5 \\
    109   & Indl. Prod. - Cons. Goods  & EA    & IP    & 5 \\
    110   & Indl. Prod. & FR    & IP    & 5 \\
    111   & Indl. Prod. - Mfg & FR    & IP    & 5 \\
    112   & Indl. Prod. - Mfg (2010=100) & FR    & IP    & 5 \\
    113   & Indl. Prod. - Manuf. Of Motor Vehicles, Trailers, Semitrailers & FR    & IP    & 5 \\
    114   & Indl. Prod. - Int. Goods & FR    & IP    & 5 \\
    115   & Indl. Prod. - Indl. Prod. - Constr. & FR    & IP    & 5 \\
    116   & Indl. Prod. - Manuf. Of Wood And Paper Products & FR    & IP    & 5 \\
    117   & Indl. Prod. - Manuf. Of Computer, Electronic And Optical Prod & FR    & IP    & 5 \\
    118   & Indl. Prod. - Manuf. Of Electrical Equipment & FR    & IP    & 5 \\
    119   & Indl. Prod. - Manuf. Of Machinery And Equipment & FR    & IP    & 5 \\
    120   & Indl. Prod. - Manuf. Of Transport Equipment & FR    & IP    & 5 \\
    121   & Indl. Prod. - Other Mfg & FR    & IP    & 5 \\
    122   & Indl. Prod. - Manuf. Of Chemicals And Chemical Products & FR    & IP    & 5 \\
    123   & Indl. Prod. - Manuf. Of Rubber And Plastics Products & FR    & IP    & 5 \\
    124   & Indl. Prod. - Investment Goods & IT    & IP    & 5 \\
    125   & Indl. Prod. & IT    & IP    & 5 \\
    126   & Indl. Prod.  & IT    & IP    & 5 \\
    127   & Indl. Prod. - Cons. Goods - Durable & IT    & IP    & 5 \\
    128   & Indl. Prod. - Investment Goods & IT    & IP    & 5 \\
    129   & Indl. Prod. - Int. Goods & IT    & IP    & 5 \\
    130   & Indl. Prod. - Chemical Products \& Synthetic Fibres & IT    & IP    & 5 \\
    131   & Indl. Prod. - Machines \& Mechanical Apparatus & IT    & IP    & 5 \\
    132   & Indl. Prod. - Means Of Transport & IT    & IP    & 5 \\
    133   & Indl. Prod. - Metal \& Metal Products & IT    & IP    & 5 \\
    134   & Indl. Prod. - Rubber Items \& Plastic Materials & IT    & IP    & 5 \\
    135   & Indl. Prod. - Wood \& Wood Products & IT    & IP    & 5 \\
    136   & Indl. Prod. & IT    & IP    & 5 \\
    137   & Indl. Prod. - Computer, Electronic And Optical Products & IT    & IP    & 5 \\
    138   & Indl. Prod. - Basic Pharmaceutical Products & IT    & IP    & 5 \\
    139   & Indl. Prod. - Constr. & DE    & IP    & 5 \\
    140   & Indl. Prod. - Ind Incl Cnstr & DE    & IP    & 5 \\
    141   & Indl. Prod. - Mfg & DE    & IP    & 5 \\
    142   & Indl. Prod. - Rebased To 1975=100 & DE    & IP    & 5 \\
    143   & Indl. Prod. - Chems \& Chem Prds & DE    & IP    & 5 \\
    144   & Indl. Prod. - Ind Excl Cnstr & DE    & IP    & 5 \\
    145   & Indl. Prod. - Ind Excl Energy \& Cnstr & DE    & IP    & 5 \\
    146   & Indl. Prod. - Mining \& Quar & DE    & IP    & 5 \\
    147   & Indl. Prod. - Cmptr, Eleccl \& Opt Prds, Elecl Eqp & DE    & IP    & 5 \\
    148   & Indl. Prod. - Interm Goods & DE    & IP    & 5 \\
    149   & Indl. Prod. - Cap. Goods & DE    & IP    & 5 \\
    150   & Indl. Prod. - Durable Cons Goods & DE    & IP    & 5 \\
    151   & Indl. Prod. - Tex \& Wearing Apparel & DE    & IP    & 5 \\
    152   & Indl. Prod. - Pulp, Paper\&Prds, Pubshg\&Print & DE    & IP    & 5 \\
    153   & Indl. Prod. - Chem Prds & DE    & IP    & 5 \\
    154   & Indl. Prod. - Rub\&Plast Prds & DE    & IP    & 5 \\
    155   & Indl. Prod. - Basic Mtls & DE    & IP    & 5 \\
    156   & Indl. Prod. - Cmptr, Eleccl \& Opt Prds, Elecl Eqp & DE    & IP    & 5 \\
    157   & Indl. Prod. - Motor Vehicles, Trailers\&Semi Trail & DE    & IP    & 5 \\
    158   & Indl. Prod. - Tex \& Wearing Apparel & DE    & IP    & 5 \\
    159   & Indl. Prod. - Paper \& Prds, Print, Reprod Of Recrd Media & DE    & IP    & 5 \\
    160   & Indl. Prod. - Chems \& Chem Prds & DE    & IP    & 5 \\
    161   & Indl. Prod. - Basic Mtls, Fab Mtl Prds, Excl Mach\&Eqp & DE    & IP    & 5 \\
    162   & Indl. Prod. - Repair \& Install Of Mach \& Eqp & DE    & IP    & 5 \\
    163   & Indl. Prod. - Mfg Excl Cnstr \& Fbt & DE    & IP    & 5 \\
    164   & Indl. Prod. - Mining \& Ind Excl Fbt & DE    & IP    & 5 \\
    165   & Indl. Prod. - Ind Excl Fbt & DE    & IP    & 5 \\
    166   & Indl. Prod. - Interm \& Cap. Goods & DE    & IP    & 5 \\
    167   & Indl. Prod. - Fab Mtl Prds Excl Mach \& Eqp & ES    & IP    & 5 \\
    168   & Indl. Prod.  & ES    & IP    & 5 \\
    169   & Indl. Prod. - Cons. Goods  & ES    & IP    & 5 \\
    170   & Indl. Prod. - Cap. Goods  & ES    & IP    & 5 \\
    171   & Indl. Prod. - Int. Goods  & ES    & IP    & 5 \\
    172   & Indl. Prod. - Energy  & ES    & IP    & 5 \\
    173   & Indl. Prod. - Cons. Goods, Non-Durables & ES    & IP    & 5 \\
    174   & Indl. Prod. - Mining  & ES    & IP    & 5 \\
    175   & Indl. Prod. - Mfg Ind  & ES    & IP    & 5 \\
    176   & Indl. Prod. - Other Mining \& Quarrying  & ES    & IP    & 5 \\
    177   & Indl. Prod. - Textile  & ES    & IP    & 5 \\
    178   & Indl. Prod. - Chemicals \& Chemical Products  & ES    & IP    & 5 \\
    179   & Indl. Prod. - Plastic \& Rubber Products  & ES    & IP    & 5 \\
    180   & Indl. Prod. - Other Non-Metal Mineral Products  & ES    & IP    & 5 \\
    181   & Indl. Prod. - Metal Processing Ind  & ES    & IP    & 5 \\
    182   & Indl. Prod. - Metal Products Excl. Machinery  & ES    & IP    & 5 \\
    183   & Indl. Prod. - Electrical Equipment  & ES    & IP    & 5 \\
    184   & Indl. Prod. - Automobile  & ES    & IP    & 5 \\
    185   & Euro Interbank Offered Rate - 3-Month (Mean) & EA    & M\&IR & 5 \\
    186   & Money Supply: Loans To Other Ea Residents Excl. Govt. & EA    & M\&IR & 5 \\
    187   & Money Supply: M3  & EA    & M\&IR & 5 \\
    188   & Euro Short Term Repo Rate & FR    & M\&IR & 5 \\
    189   & Datastream Euro Share Price Index (Mth. Avg.) & FR    & M\&IR & 1 \\
    190   & Euribor: 3-Month (Mth. Avg.) & FR    & M\&IR & 5 \\
    191   & Mfi Loans To Resident Private Sector & FR    & M\&IR & 5 \\
    192   & Money Supply - M1 & FR    & M\&IR & 5 \\
    193   & Money Supply - M3 & FR    & M\&IR & 5 \\
    194   & Share Price Index - Sbf 250 & DE    & M\&IR & 1 \\
    195   & Fibor - 3 Month (Mth.Avg.) & DE    & M\&IR & 5 \\
    196   & Money Supply - M3 & DE    & M\&IR & 5 \\
    197   & Money Supply - M2 & DE    & M\&IR & 5 \\
    198   & Bank Prime Lending Rate / Ecb Marginal Lending Facility & DE    & M\&IR & 5 \\
    199   & Dax Share Price Index, Ep & IT    & M\&IR & 1 \\
    200   & Interbank Deposit Rate-Average On 3-Months Deposits & IT    & M\&IR & 5 \\
    201   & Official Reserve Assets & ES    & M\&IR & 5 \\
    202   & Money Supply: M3 - Spanish  & ES    & M\&IR & 5 \\
    203   & Madrid S.E - General Index & ES    & M\&IR & 5 \\
    204   & Hicp - Overall Index & EA    & HCPI  & 6 \\
    205   & Hicp - All-Items Excluding Energy, Index & EA    & HCPI  & 6 \\
    206   & Hicp - Food Incl. Alcohol And Tobacco, Index & EA    & HCPI  & 6 \\
    207   & Hicp - Processed Food Incl. Alcohol And Tobacco, Index & EA    & HCPI  & 6 \\
    208   & Hicp - Unprocessed Food, Index & EA    & HCPI  & 6 \\
    209   & Hicp - Goods, Index & EA    & HCPI  & 6 \\
    210   & Hicp - Industrial Goods, Index & EA    & HCPI  & 6 \\
    211   & Hicp - Industrial Goods Excluding Energy, Index & EA    & HCPI  & 6 \\
    212   & Hicp - Services, Index & EA    & HCPI  & 6 \\
    213   & Hicp - All-Items Excluding Tobacco, Index & EA    & HCPI  & 6 \\
    214   & Hicp - All-Items Excluding Energy And Food, Index & EA    & HCPI  & 6 \\
    215   & Hicp - All-Items Excluding Energy And Unprocessed Food, Index & EA    & HCPI  & 6 \\
    216   & All-Items Hicp & DE    & HCPI  & 6 \\
    217   & All-Items Hicp & ES    & HCPI  & 6 \\
    218   & All-Items Hicp & FR    & HCPI  & 6 \\
    219   & All-Items Hicp & IT    & HCPI  & 6 \\
    220   & Goods (Overall Index Excluding Services) & DE    & HCPI  & 6 \\
    221   & Goods (Overall Index Excluding Services) & FR    & HCPI  & 6 \\
    222   & Processed Food Including Alcohol And Tobacco & DE    & HCPI  & 6 \\
    223   & Processed Food Including Alcohol And Tobacco & ES    & HCPI  & 6 \\
    224   & Processed Food Including Alcohol And Tobacco & FR    & HCPI  & 6 \\
    225   & Processed Food Including Alcohol And Tobacco & IT    & HCPI  & 6 \\
    226   & Unprocessed Food & DE    & HCPI  & 6 \\
    227   & Unprocessed Food & ES    & HCPI  & 6 \\
    228   & Unprocessed Food & FR    & HCPI  & 6 \\
    229   & Unprocessed Food & IT    & HCPI  & 6 \\
    230   & Non-Energy Industrial Goods & DE    & HCPI  & 6 \\
    231   & Non-Energy Industrial Goods & FR    & HCPI  & 6 \\
    232   & Services (Overall Index Excluding Goods) & DE    & HCPI  & 6 \\
    233   & Services (Overall Index Excluding Goods) & FR    & HCPI  & 6 \\
    234   & Overall Index Excluding Tobacco & DE    & HCPI  & 6 \\
    235   & Overall Index Excluding Tobacco & FR    & HCPI  & 6 \\
    236   & Overall Index Excluding Energy & DE    & HCPI  & 6 \\
    237   & Overall Index Excluding Energy & FR    & HCPI  & 6 \\
    238   & Overall Index Excluding Energy And Unprocessed Food & DE    & HCPI  & 6 \\
    239   & Overall Index Excluding Energy And Unprocessed Food & FR    & HCPI  & 6 \\
    240   & Ppi: Ind Excluding Constr. \& Energy  & EA    & PPI   & 6 \\
    241   & Ppi: Cap. Goods  & EA    & PPI   & 6 \\
    242   & Ppi: Non-Durable Cons. Goods  & EA    & PPI   & 6 \\
    243   & Ppi: Int. Goods  & EA    & PPI   & 6 \\
    244   & Ppi: Non Dom. - Mining, Mfg \& Quarrying  & EA    & PPI   & 6 \\
    245   & Ppi: Non Dom. Mfg  & DE    & PPI   & 6 \\
    246   & Ppi: Int. Goods Excluding Energy & DE    & PPI   & 6 \\
    247   & Ppi: Cap. Goods & DE    & PPI   & 6 \\
    248   & Ppi: Cons. Goods & DE    & PPI   & 6 \\
    249   & Ppi: Fuel & DE    & PPI   & 6 \\
    250   & Ppi: Indl. Products (Excl. Energy) & DE    & PPI   & 6 \\
    251   & Ppi: Machinery & DE    & PPI   & 6 \\
    252   & Deflated T/O: Ret. Sale In Non-Spcld Str With Food, Bev \& Tob & DE    & T/O   & 5 \\
    253   & Deflated T/O: Oth Ret. Sale In Non-Spcld Str & DE    & T/O   & 5 \\
    254   & Deflated T/O: Sale Of Motor Vehicle Pts \& Acces & DE    & T/O   & 5 \\
    255   & Deflated T/O: Wholesale Of Agl Raw Matls \& Live Animals & DE    & T/O   & 5 \\
    256   & Deflated T/O: Wholesale Of Household Goods & IT    & T/O   & 5 \\
    257   & T/O: Ret. Trd, Exc Of Mv , Motorcyles \& Fuel & ES    & T/O   & 5 \\
    258   & T/O: Ret. Sale Of Clth \& Leath Gds In Spcld Str & ES    & T/O   & 5 \\
    259   & T/O: Ret. Sale Of Non-Food Prds (Exc Fuel) & ES    & T/O   & 5 \\
    260   & T/O: Ret. Sale Of Info, Househld \& Rec Eqp In Spcld Str & ES    & T/O   & 5 \\
    261   & Ek Unemployment: All & EA    & HUR   & 5 \\
    262   & Ek Unemployment: Persons Over 25 Years Old  & EA    & HUR   & 5 \\
    263   & Ek Unemployment: Women Under 25 Years Old  & EA    & HUR   & 5 \\
    264   & Ek Unemployment: Women Over 25 Years Old  & EA    & HUR   & 5 \\
    265   & Ek Unemployment: Men Over 25 Years Old  & EA    & HUR   & 5 \\
    266   & Fr Hur All Persons (All Ages)  & FR    & HUR   & 5 \\
    267   & Fr Hur Femmes (Ages 15-24)  & FR    & HUR   & 5 \\
    268   & Fr Hur Femmes (All Ages)  & FR    & HUR   & 5 \\
    269   & Fr Hur Hommes (Ages 15-24)  & FR    & HUR   & 5 \\
    270   & Fr Hur Hommes (All Ages)  & FR    & HUR   & 5 \\
    271   & Fr Hur All Persons (Ages 15-24)  & FR    & HUR   & 5 \\
    272   & Fr Hurall Persons(Ages 25 And Over)  & FR    & HUR   & 5 \\
    273   & Fr Hur Females (Ages 25 And Over)  & FR    & HUR   & 5 \\
    274   & Fr Hur Males (Ages 25 And Over)  & FR    & HUR   & 5 \\
    275   & Bd Hur All Persons (All Ages)  & DE    & HUR   & 5 \\
    276   & Bd Hur Femmes (Ages 15-24)  & DE    & HUR   & 5 \\
    277   & Bd Hur Femmes (All Ages)  & DE    & HUR   & 5 \\
    278   & Bd Hur Hommes (Ages 15-24)  & DE    & HUR   & 5 \\
    279   & Bd Hur Hommes (All Ages)  & DE    & HUR   & 5 \\
    280   & Bd Hur All Persons (Ages 15-24)  & DE    & HUR   & 5 \\
    281   & Bd Hurall Persons(Ages 25 And Over)  & DE    & HUR   & 5 \\
    282   & Bd Hur Females (Ages 25 And Over)  & DE    & HUR   & 5 \\
    283   & Bd Hur Males (Ages 25 And Over)  & DE    & HUR   & 5 \\
    284   & It Hur All Persons (All Ages)  & IT    & HUR   & 5 \\
    285   & It Hur Femmes (All Ages)  & IT    & HUR   & 5 \\
    286   & It Hur Hommes (All Ages)  & IT    & HUR   & 5 \\
    287   & It Hur All Persons (Ages 15-24)  & IT    & HUR   & 5 \\
    288   & It Hurall Persons(Ages 25 And Over)  & IT    & HUR   & 5 \\
    289   & Es Hur All Persons (All Ages)  & ES    & HUR   & 5 \\
    290   & Es Hur Femmes (Ages 16-24)  & ES    & HUR   & 5 \\
    291   & Es Hur Femmes (All Ages)  & ES    & HUR   & 5 \\
    292   & Es Hur Hommes (Ages 16-24)  & ES    & HUR   & 5 \\
    293   & Es Hur Hommes (All Ages)  & ES    & HUR   & 5 \\
    294   & Es Hur All Persons (Ages 16-24)  & ES    & HUR   & 5 \\
    295   & Es Hurall Persons(Ages 25 And Over)  & ES    & HUR   & 5 \\
    296   & Es Hur Females (Ages 25 And Over)  & ES    & HUR   & 5 \\
    297   & Es Hur Males (Ages 25 And Over)  & ES    & HUR   & 5 \\
    298   & De - Service  Confidence Indicator & DE    & SI    & 1 \\
    299   & De Services - Buss. Dev. Past 3 Months & DE    & SI    & 1 \\
    300   & De Services - Evol. Demand Past 3 Months & DE    & SI    & 1 \\
    301   & De Services - Exp. Demand Next 3 Months & DE    & SI    & 1 \\
    302   & De Services - Evol. Employ. Past 3 Months & DE    & SI    & 1 \\
    303   & Fr - Service  Confidence Indicator & FR    & SI    & 1 \\
    304   & Fr Services - Buss. Dev. Past 3 Months & FR    & SI    & 1 \\
    305   & Fr Services - Evol. Demand Past 3 Months & FR    & SI    & 1 \\
    306   & Fr Services - Exp. Demand Next 3 Months & FR    & SI    & 1 \\
    307   & Fr Services - Evol. Employ. Past 3 Months & FR    & SI    & 1 \\
    308   & Fr Services - Exp. Employ. Next 3 Months & FR    & SI    & 1 \\
    309   & Fr Services - Exp. Prices Next 3 Months & FR    & SI    & 1 \\

\end{longtable}

\section{Distribution of $\widehat{Cov}(u_1,u_2)$}\label{Dist_Cu}

\begin{figure}[H]
\begin{center}
\hskip1.2cm
\begin{minipage}[c]{0.4\textwidth}
\graphicspath{{images/}}
\includegraphics[scale=0.28]{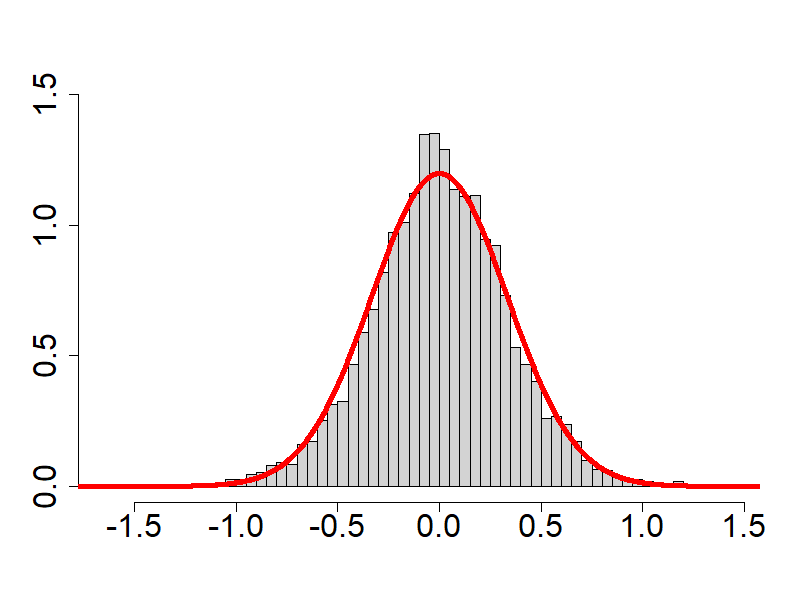}
\subcaption{$T=10$}
\end{minipage}\hfill{}%
\begin{minipage}[c]{0.4\textwidth}%
\graphicspath{{images/}}
\includegraphics[scale=0.28]{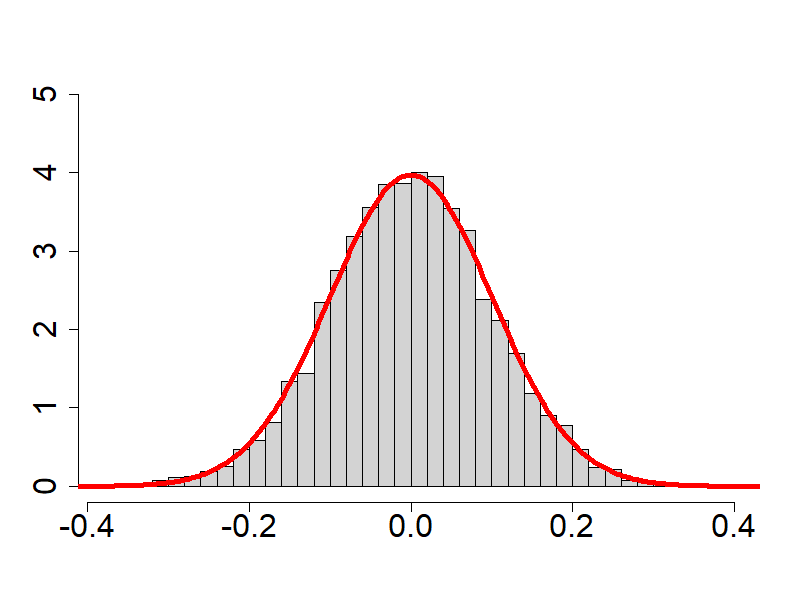}
\subcaption{$T=100$}
\end{minipage}
\end{center}
\caption{\footnotesize Density of $\widehat{Cov}(u_1,u_2)$ between two uncorrelated standard Normal variables for $T=10$ (a) and $T=100$ (b).}
\label{fig:CovU}
\end{figure}

In Figure \ref{fig:CovU} we report the density of $\widehat{Cov}(u_1,u_2)$ when $u_1$ and $u_2$ are standard Normal in the cases of $T=10$ and 100. Red line shows the density of $N\left(0, \frac{1}{T-1}\right)$. Observations are obtained on 5000 Monte Carlo replications. We observe that the approximation of $\widehat{Cov}(u_1,u_2)$ to $N(0,\frac{1}{T-1})$ holds also when $T$ is small (see Figure \ref{fig:CovU} (a) relative to $T$=10). This analysis corroborate numerically the results in \cite{Glen2004}, which show that if $x$ and $y$ are $N\left(0,1\right)$, then the probability density function of $xy$ is $\frac{K_0(|xy|)}{pi}$, where $K_0(|xy|)$ is the Bessel function of the second kind. 

\end{appendices}

\end{document}